\documentclass[12pt]{amsart}
\usepackage[latin1]{inputenc}
\usepackage{amssymb}
\DeclareMathAlphabet{\mathfr}{U}{euf}{m}{n}
\usepackage{latexsym,amscd,longtable,supertabular,multirow,booktabs,array}
\input xy       
\xyoption{all}

\DeclareFontEncoding{OT2}{}{} 
\newcommand{\textcyr}[1]{{\fontencoding{OT2}\fontfamily{wncyr}\fontseries{m}\fontshape{n}\selectfont #1}}
\newcommand{\Sha}{{\mbox{\textcyr{Sh}}}}

\newtheorem{cor}[equation]{Corollary}

\newtheorem{prop}[equation]{Proposition}

\numberwithin{equation}{section}

\newtheorem{theorem}[equation]{Theorem}
\newtheorem{lemma}[equation]{Lemma}
\newtheorem{conj}[equation]{Conjecture}

\theoremstyle{definition}

\newtheorem{question}[equation]{Question}

\theoremstyle{remark}
\newtheorem{remark}[equation]{Remark}

\newcommand{\Q}{\mathbb{Q}}
\newcommand{\Z}{\mathbb{Z}}
\newcommand{\F}{\mathbb{F}}

\newcommand{\C}{\mathbb{C}}
\newcommand{\D}{{\mathcal{D}}}
\newcommand{\Coeff}{{\operatorname{Coeff}}}
\newcommand{\Gal}{\operatorname{Gal\,}}
\newcommand{\Aut}{{\operatorname{Aut}}} 
\newcommand{\GL}{\operatorname{GL}}
\newcommand{\SL}{\operatorname{SL}}
\newcommand{\End}{\operatorname{End}}
\newcommand{\Frob}{\operatorname{Frob}}

\newcommand{\dv}{\operatorname{div}}
\newcommand{\id}{\operatorname{id}}
\newcommand{\ord}{\operatorname{ord}}

\newcommand{\Tr}{\operatorname{Tr}}

\newcommand{\X}{{\mathcal X}}
\newcommand{\Y}{{\mathcal Y}}
\newcommand{\Cg}{{\mathcal C}_g}
\newcommand{\Vbar}{{\overline{V}}}
\newcommand{\Fbar}{{\overline{\F}}}
\newcommand{\Qisog}{\stackrel{\Q}\sim}
\newcommand{\directsum}{\oplus}
\newcommand{\Directsum}{\bigoplus}
\newcommand{\tensor}{\otimes}
\newcommand{\nichts}{\ensuremath{\left.\right.}}
\newcommand{\ssigma}{{\! \nichts^\sigma \!}}
\newcommand{\ttau}{{\!\! \nichts^\tau \!}}
\newcommand{\isom}{\simeq}
\newcommand{\PGL}{\operatorname{PGL}}
\newcommand{\R}{{\mathbb R}}
\newcommand{\PP}{{\mathbb P}}
\newcommand{\OO}{{\mathcal O}}
\newcommand{\Ohat}{{\hat{\OO}}}

\newcommand{\xtilde}{{\tilde{x}}}
\newcommand{\kbar}{{\overline{k}}}
\newcommand{\Qbar}{{\overline{\Q}}}
\newcommand{\Spec}{{\operatorname{Spec}}}
\newcommand{\Hom}{{\operatorname{Hom}}}
\newcommand{\Jac}{{\operatorname{Jac}}}
\newcommand{\New}{{\operatorname{New}}}
\newcommand{\new}{{\operatorname{new}}}
\newcommand{\old}{{\operatorname{old}}}
\newcommand{\odd}{{\operatorname{odd}}}
\newcommand{\Sparse}{{\operatorname{Sparse}}}
\newcommand{\Smooth}{{\operatorname{Smooth}}}
\newcommand{\cond}{{\operatorname{cond}}}
\newcommand{\abcd}{{\begin{pmatrix} a & b \\ c & d \end{pmatrix}}}
\newcommand{\hh}{{\mathcal H}}
\newcommand{\pp}{{\mathfrak p}}
\newcommand{\im}{{\operatorname{Im}}}
\newcommand{\Edual}{{\E}^\vee}
\newcommand{\ahat}{{\widehat{a}}}

\newcommand{\calE}{{\mathcal E}}

\newcommand{\Wei}{{\operatorname{WP}}}
\newcommand{\E}{\mathbb{E}}

\newcommand{\eps}{\varepsilon}

\newcolumntype{L}[1]{>{\rightskip=0pt plus1fil}p{#1\textwidth}<{}}
\newcolumntype{X}{L{0.10}@{:}L{0.81}}
\newcolumntype{Y}{L{0.15}@{:}L{0.81}}

\begin{document}

\title[Modular curves of genus at least $2$]{Finiteness results for modular
curves of genus at least $2$}
\date{December 20, 2003}

\author[Baker]{Matthew H. Baker}
\address{Department of Mathematics, University of Georgia, Athens, GA 30602-7403, USA}
\email{mbaker@math.uga.edu}
\thanks{The first author was supported in part by an NSF Postdoctoral Fellowship.}

\author[Gonz{\'a}lez-Jim{\'e}nez]{Enrique Gonz{\'a}lez-Jim{\'e}nez}
\address{School of Mathematical Sciences. University of Nottingham. University Park. Nottingham, NG7 2RD, UK} 
\email{enrique.gonzalezj@uah.es}


\thanks{The second author was supported in part by DGICYT Grant BHA2000-0180 and a postdoctoral fellowship from the European Research Training Network "Galois Theory and Explicit Methods in Arithmetic".}

\author[Gonz{\'a}lez]{\hbox{Josep Gonz{\'a}lez}}
\address{Escola Universit\`aria Polit\`ecnica de Vilanova i la
  Geltr{\'u}, Av. V\'{i}ctor Balaguer s/n, E-08800 Vilanova i la Geltr${\cdot}$, Spain}
\email{josepg@mat.upc.es}
\thanks{The third author was supported in part by DGICYT Grant
BFM2003-06768-C02-02.}

\author[Poonen]{\hbox{Bjorn Poonen}}
\address{Department of Mathematics,
        University of California, Berkeley, CA 94720-3840, USA}
\email{poonen@math.berkeley.edu}
\thanks{The fourth author was supported by NSF grant DMS-9801104,
        and a Packard Fellowship.}

\subjclass[2000]{Primary 11G18, 14G35; Secondary 11F11, 11G10, 14H45.}
\keywords{Modular curves, modular forms, de~Franchis-Severi Theorem,
gonality, computability, hyperelliptic curves, automorphisms of curves}

\begin{abstract}
A curve $X$ over $\Q$ is modular if it is dominated by $X_1(N)$
for some $N$; if in addition the image of its jacobian in $J_1(N)$
is contained in the new subvariety of $J_1(N)$, then $X$ is called 
a new modular curve.
We prove that for each $g \ge 2$, the set of new modular curves
over $\Q$ of genus $g$ is finite and computable.
For the computability result, we prove an algorithmic version
of the de~Franchis-Severi Theorem.
Similar finiteness results are proved for new modular curves
of bounded gonality, for new modular curves whose jacobian is
a quotient of $J_0(N)^\new$ with $N$ divisible by a prescribed prime,
and for modular curves (new or not) with levels in a restricted set.
We study new modular hyperelliptic curves in detail.
In particular, we find all new modular curves of genus~2 explicitly, 
and construct what might be the complete list of
all new modular hyperelliptic curves of all genera.
Finally we prove that for each field $k$ of characteristic zero
and $g \ge 2$, the set of genus-$g$ curves over $k$
dominated by a Fermat curve is finite and computable.
\end{abstract}

\maketitle

\section{Introduction}
\label{introduction}

Let $X_1(N)$ be the usual modular curve over $\Q$.
(See Section~\ref{modularformsection} for a definition.)
A curve\footnote{Curves and varieties in this paper are smooth, projective,
and geometrically integral, unless otherwise specified.
When we write an affine equation for a curve, its smooth projective model
is implied.}
 $X$ over $\Q$ will be called {\em modular} if there exists
a nonconstant morphism $\pi:X_1(N) \rightarrow X$ over $\Q$.
If $X$ is modular, then $X(\Q)$ is nonempty,
since it contains the image of the cusp $\infty \in X_1(N)(\Q)$.
The converse, namely that if $X(\Q)$ is nonempty then $X$ is modular,
holds if the genus $g$ of $X$ satisfies $g \le 1$~\cite{breuil2001}.
In particular, there are infinitely many modular curves
over $\Q$ of genus $1$.
On the other hand, we propose the following:

\begin{conj}
\label{fixedgenusconj}
For each $g \ge 2$,
the set of modular curves over $\Q$ of genus $g$ is finite.
\end{conj}

\begin{remark}
$\left.\right.$
\begin{enumerate}
\item[(i)] When we speak of the finiteness
of the set of curves over $\Q$ satisfying some condition,
we mean the finiteness of the set of $\Q$-isomorphism classes of such curves.
\item[(ii)] For any fixed $N$, the {\em de~Franchis-Severi Theorem}
(see Theorem~\ref{defranchis})
implies the finiteness of the set of curves over $\Q$ dominated by $X_1(N)$.
Conjecture~\ref{fixedgenusconj} can be thought of as a version
that is uniform as one ascends the tower of modular curves $X_1(N)$,
provided that one fixes the genus of the dominated curve.
\item[(iii)] Conjecture~\ref{fixedgenusconj} is true if one
restricts the statement to quotients of $X_1(N)$
by subgroups of its group of modular automorphisms.
See Remark~\ref{X^*(N)} for details.
\item[(iv)] If $X_1(N)$ dominates a curve $X$,
then the jacobian $\Jac X$ is a quotient\footnote{
Quotients or subvarieties of varieties,
and morphisms between varieties,
are implicitly assumed to be defined over the same field
as the original varieties.
If $X$ is a curve over $\Q$, 
and we wish to discuss automorphisms over $\C$, for example,
we will write $\Aut(X_\C)$.
Quotients of abelian varieties are assumed to be abelian variety quotients.}
 of $J_1(N):=\Jac X_1(N)$.
The converse, namely that if $X$ is a curve such that $X(\Q)$ is nonempty
and $\Jac X$ is a quotient of $J_1(N)$ then $X$ is dominated by $X_1(N)$,
holds if the genus $g$ of $X$ is $\le 1$,
but can fail for $g \ge 2$.
See Section~\ref{pathologysection} for other ``pathologies.''
\item[(v)]
In contrast with Conjecture~\ref{fixedgenusconj},
there exist infinitely many genus-two curves over $\Q$
whose jacobians are quotients of $J_1(N)$ for some $N$.
See Proposition~\ref{pathologies}(\ref{modularjacobian=modularcurve}).  
\item[(vi)] In Section~\ref{fermatsection}, we use 
a result of Aoki~\cite{Aoki} to prove an analogue
of Conjecture~\ref{fixedgenusconj} in which $X_1(N)$
is replaced by the Fermat curve $x^N+y^N=z^N$ in $\PP^2$.
In fact, such an analogue can be proved over arbitrary fields
of characteristic zero, not just $\Q$.
\end{enumerate}
\end{remark}

We prove many results towards
Conjecture~\ref{fixedgenusconj} in this paper.
Given a variety $X$ over a field $k$,
let $\Omega=\Omega^1_{X/k}$ denote the sheaf of regular $1$-forms.
Call a modular curve $X$ over $\Q$ {\em new of level $N$}
if there exists a nonconstant morphism $\pi: X_1(N) \rightarrow X$
(defined over $\Q$) such that $\pi^* H^0(X,\Omega)$
is contained in the new subspace $H^0(X_1(N),\Omega)_\new$,
or equivalently if the image of the homomorphism
$\pi^\ast: \Jac X \rightarrow J_1(N)$ induced by Picard functoriality
is contained in the new subvariety $J_1(N)_\new$ of $J_1(N)$.
(See Section~\ref{modularformsection} for the definitions of 
$H^0(X_1(N),\Omega)_\new$, $J_1(N)_\new$, $J_1(N)^\new$, 
and so on.)
For example, it is known that every elliptic curve $E$ over $\Q$
is a new modular curve of level $N$, where $N$ is the conductor of $E$.
Here the conductor $\cond(A)$ of an abelian variety $A$ over $\Q$ is
a positive integer $\prod_p p^{f_p}$,
where each exponent $f_p$ is defined in terms of the action of 
an inertia subgroup of $\Gal(\Qbar/\Q)$ at the prime $p$
on a Tate module $T_\ell A$: 
for the full definition, 
see~\cite[\S1]{lockhart1993}, for instance.
If $A$ is a quotient of $J_1(N)^\new$,
then $\cond(A)=N^{\dim A}$, by~\cite{carayol}.

\begin{theorem}
\label{genustheorem}
For each $g \ge 2$, the set of {\em new}
modular curves over $\Q$ of genus $g$ is finite and computable.
\end{theorem}

The first results of this type were proved in~\cite{gogo00},
which showed that the set of new modular curves of genus $2$ over $\Q$
is finite, and that there are exactly $149$ such curves
whose jacobian is $\Q$-simple.
Equations for these $149$ curves are given in
Tables 1 and~2 of~\cite{gogo00}.

\begin{remark}
$\left.\right.$
\begin{enumerate}
\item[(i)] See Section~\ref{computablesection} for the precise definition of
``computable.''
\item[(ii)]
An analysis of the proof of Theorem~\ref{genustheorem}
would show that as $g \rightarrow \infty$,
there are at most $\exp((6+o(1)) g^2)$
new modular curves of genus $g$ over $\Q$.
The truth is probably much smaller, however.
It is conceivable even that there is an upper bound not depending on $g$.
\item[(iii)]
If we consider all genera $g \ge 2$ together,
then there are infinitely many new modular curves.
For example, $X_1(p)$ 
is a new modular curve whenever $p$ is prime,
and its genus tends to infinity with $p$.
\end{enumerate}
\end{remark}

If we drop the assumption that our modular curves are new,
we can still prove results,
but (so far) only if we impose restrictions on the level.
Given $m>0$,
let $\Sparse_m$ denote the set of positive integers $N$
such that if $1=d_1 < d_2 < \dots < d_t = N$ are the positive divisors of $N$,
then $d_{i+1}/d_i > m$ for $i=1,\dots,t-1$.
Define a function $B(g)$ on integers $g \ge 2$
by $B(2)=13$, $B(3)=17$, $B(4)=21$, and $B(g)=6g-5$ for $g \ge 5$.
(For the origin of this function,
see the proofs of
Propositions \ref{curveisdetermined} and~\ref{curveisdetermined2}.)
A positive integer $N$ is called {\em $m$-smooth}
if all primes $p$ dividing $N$ satisfy $p \le m$.
Let $\Smooth_m$ denote the set of $m$-smooth integers.

\begin{theorem}
\label{nonnewtheorem}
Fix $g \ge 2$, and let $S$ be a subset of $\{1,2,\dots\}$.
The set of modular curves over $\Q$ of genus $g$ and of level contained
in $S$ is finite if any of the following hold:
\begin{enumerate}
\item[(i)] $S=\Sparse_{B(g)}$.
\item[(ii)] $S=\Smooth_m$ for some $m>0$.
\item[(iii)] $S$ is the set of prime powers.
\end{enumerate}
\end{theorem}

\begin{remark}
Since $\Sparse_{B(g)} \cup \Smooth_{B(g)}$ contains all prime powers,
parts (i) and~(ii) of Theorem~\ref{nonnewtheorem} imply (iii).
\end{remark}

\begin{remark}
In contrast with Theorem~\ref{genustheorem},
we do not know, even in theory, how to compute the finite sets
of curves in Theorem~\ref{nonnewtheorem}.
The reason for this will be explained in Remark~\ref{nonnewisnotcomputable}.
\end{remark}

If $X$ is a curve over a field $k$, and $L$ is a field extension of $k$,
let $X_{L}$ denote $X \times_k L$.
The {\em gonality} $G$ of a curve $X$ over $\Q$
is the smallest possible degree of a nonconstant morphism
$X_{\C} \rightarrow \PP^1_{\C}$.
(There is also the notion of {\em $\Q$-gonality}, where one only allows
morphisms over $\Q$.
By defining gonality using morphisms
over $\C$ instead of $\Q$, we make the next theorem stronger.)
In Section~\ref{gonalitysection},
we combine Theorem~\ref{genustheorem}
with a known lower bound
on the gonality of $X_1(N)$ to prove the following:

\begin{theorem}
\label{gonalitytheorem}
For each $G \ge 2$, the set of new modular curves over $\Q$
of genus at least $2$ and gonality at most $G$ is finite and computable.
\end{theorem}

\noindent
(We could similarly prove an analogue of Theorem~\ref{nonnewtheorem}
for curves of bounded gonality instead of fixed genus.)

Recall that a curve $X$ of genus $g$ over a field $k$
is called {\em hyperelliptic} if
$g \ge 2$ and the canonical map $X \rightarrow \PP^{g-1}$
is not a closed immersion:
equivalently, $g \ge 2$ and
there exists a degree-$2$ morphism $X \rightarrow Y$ where $Y$ has genus zero.
If moreover $X(k) \not=\emptyset$ then $Y \isom \PP^1_k$,
and if also $k$ is not of characteristic $2$,
then $X$ is birational to a curve of the form $y^2=f(x)$
where $f$ is a separable polynomial in $k[x]$ of degree $2g+1$ or $2g+2$.
Recall that $X(k) \not=\emptyset$ is automatic if $X$ is modular,
because of the cusp $\infty$.

Taking $G=2$ in Theorem~\ref{gonalitytheorem}, we find that
the set of new modular hyperelliptic curves over $\Q$
is finite and computable.
We can say more:

\begin{theorem}
\label{hyperelliptictheorem}
Let $X$ be a new modular hyperelliptic curve over $\Q$ of genus $g \ge 3$
and level $N$.
Then
\begin{enumerate}
\item[(i)] $g \le 16$.
\item[(ii)] If $\Jac X$ is a quotient of $J_0(N)$, then  $g \le 10$. If moreover  $3\mid N$,
then $X$ is the genus-$3$ curve $X_0(39)$.
\item[(iii)] If $\Jac X$ is not a quotient of $J_0(N)$, 
then either $g$ is even or $g \le 9$.
\end{enumerate}
\end{theorem}
\noindent
Further information is given in 
Sections \ref{hypfin} and~\ref{computationalresults}, and in the appendix.
See Section~\ref{modularformsection} for the definitions
of $X_0(N)$ and $J_0(N)$.

\medskip

As we have already remarked, if we consider all genera $g \ge 2$ together,
there are infinitely many new modular curves.
To obtain finiteness results, so far we have needed to restrict
either the genus or the gonality.
The following theorem, proved in Section~\ref{trivialcharsection}, 
gives a different type of restriction that
again implies finiteness.

\begin{theorem}
\label{pinleveltheorem}
For each prime $p$, the set of new modular curves over $\Q$ 
of genus at least $2$ 
whose jacobian is a quotient of $J_0(N)^\new$ for some $N$ divisible by $p$
is finite and computable.
\end{theorem}

\begin{question}
Does Theorem~\ref{pinleveltheorem} remain true if
$J_0(N)^\new$ is replaced by $J_1(N)^\new$?
\end{question}

Call a curve $X$ over a field $k$ of characteristic zero
{\em $k$-modular} if there exists a nonconstant morphism
$X_1(N)_k \rightarrow X$ (over $k$).

\begin{question}
\label{extension}
Is it true that for every field $k$ of characteristic zero,
and every $g \ge 2$,
the set of $k$-modular curves over $k$ of genus $g$ up to $k$-isomorphism
is finite?
\end{question}

\begin{remark}
If $X$ is a $k$-modular curve over $k$,
and we define $k_0 = k \cap \Qbar \subseteq \kbar$,
then $X = X_0 \times_{k_0} k$ for some $k_0$-modular curve $X_0$.
This follows from the de~Franchis-Severi Theorem.
\end{remark}

%

\begin{remark}
If $k$ and $k'$ are fields of characteristic zero
with $[k':k]$ {\em finite},
then a positive answer to Question~\ref{extension} for $k'$
implies a positive answer for $k$,
since Galois cohomology and the finiteness of automorphism
group of curves of genus at least $2$
show that for each $X'$ over $k'$,
there are at most finitely many curves $X$ over $k$
with $X \times_k k' \isom X'$.
But it is {\em not} clear, for instance,
that a positive answer for $\Qbar$ implies
or is implied by a positive answer for $\Q$.
\end{remark}

\section{Recovering curve information from differentials}
\label{recoveringsection}

\subsection{Recovering a curve from partial expansions of its differentials}
\label{recovering}

The goal of this section is to prove the following result, which will
be used frequently in the rest of this paper.

\begin{prop}
\label{curveisdetermined}
Fix an integer $g \ge 2$.
There exists an integer $B>0$ depending on $g$ such that
if $k$ is a field of characteristic zero,
and $w_1$, \dots, $w_g$ are elements of $k[[q]]/(q^B)$,
then up to $k$-isomorphism,
there exists at most one curve $X$ over $k$
such that there exist $P \in X(k)$
and an analytic uniformizing parameter $q$
in the completed local ring $\Ohat_{X,P}$
such that $w_1 \, dq$, \dots, $w_g \, dq$
are the expansions modulo $q^B$ of some basis of $H^0(X,\Omega)$.
\end{prop}

\begin{lemma}
\label{sections}
Let $k$ be a field, and let $X/k$ be a curve of genus $g$.  Let
$P\in X(k)$ be a $k$-rational point, let $F\in k[t_1,\ldots,t_g]$ be a
homogeneous polynomial of degree $d$, and let $q$ be an
analytic uniformizing parameter in $\Ohat_{X,P}$.  Suppose we are given
elements $\omega_1,\ldots,\omega_g \in H^0(X,\Omega)$, and for each
$i=1,\ldots,g$, write $\omega_i = w_i dq$ with $w_i \in k[[q]]$.  Then if
$F(w_1,\dots,w_g) \equiv 0 \pmod{q^B}$ with $B>d(2g-2)$,
then $F(\omega_1,\dots,\omega_g) = 0$ in $H^0(X,\Omega^{\tensor d})$.
\end{lemma}

\begin{proof}
A {\em nonzero} element of $H^0(X,\Omega^{\tensor d})$
has $d(2g-2)$ zeros on $X$,
so its expansion at any given point can vanish to order at most $d(2g-2)$.
\end{proof}

The following is a weak form of a theorem of Petri
appearing as Theorem~2.3 on page 131 of~\cite{arbarello1985},
for example.

\begin{theorem}[Petri, 1923]
\label{petri}
Let $X$ be a nonhyperelliptic curve of genus $g \ge 4$
over a field $k$ of characteristic zero.
Then the image of the canonical map
$X \rightarrow \PP^{g-1}$
is the common zero locus of some set
of homogeneous polynomials of degree $2$ and~$3$.
Moreover, if $X$ is neither trigonal nor a smooth plane quintic,
then degree~$2$ polynomials suffice.
\end{theorem}

\begin{cor}
\label{superpetri}
Let $X$ be a curve of genus $g \ge 2$
over a field $k$ of characteristic zero.
Then the image $X'$ of the canonical map $X \rightarrow \PP^{g-1}$
is the common zero locus of the set of homogeneous polynomials of degree~4
that vanish on $X'$.
\end{cor}

\begin{proof}
We may assume that $k$ is algebraically closed.
If $X$ is hyperelliptic of genus $g$,
say birational to $y^2=f(x)$ where $f$ has distinct roots,
then we may choose $\{\, x^i \, dx/y : 0 \le i \le g-1 \,\}$
as basis of $H^0(X,\Omega)$, and then the image of the canonical
map is the rational normal curve
cut out by $\{\, t_i t_j - t_{i'} t_{j'} : i+j = i'+j' \,\}$
where $t_0$, \dots, $t_{g-1}$ are the homogeneous coordinates
on $\PP^{g-1}$.
If $X$ is nonhyperelliptic of genus~$3$,
its canonical model is a plane quartic.
In all other cases, we use Petri's Theorem.
(The zero locus of a homogeneous polynomial $h$
of degree $d<4$ equals the zero locus of the set
of homogeneous polynomials of degree~$4$ that are multiples of $h$.)
\end{proof}

\begin{lemma}
\label{hyperelliptic}
Let $X$ be a hyperelliptic curve of genus $g$
over a field $k$ of characteristic zero, and suppose $P \in X(k)$.
Let $\{\omega_1,\dots,\omega_g\}$ be a basis of $H^0(X,\Omega)$
such that $\ord_P(\omega_1) < \dots < \ord_P(\omega_g)$.
Then $x:=\omega_{g-1}/\omega_g$ and $y:=dx/\omega_g$
generate the function field $k(X)$,
and there is a unique polynomial $F(x)$ of degree at most $2g+2$
such that $y^2=F(x)$.
Moreover, $F$ is squarefree.
If $P$ is a Weierstrass point,
then $\deg F = 2g+1$ and $\ord_P(\omega_i)=2i-2$ for all $i$;
otherwise $\deg F = 2g+2$ and $\ord_P(\omega_i)=i-1$ for all $i$.
Finally, it is possible to replace each $\omega_i$
by a linear combination of $\omega_i$, $\omega_{i+1}$, \dots, $\omega_g$
to make $\omega_i = x^{g-i} \, dx/y$ for $1 \le i \le g$.
\end{lemma}

\begin{proof}
This follows easily from Lemma 3.6.1, Corollary 3.6.3, and Theorem
3.6.4 of \cite{goldschmidt2002}.
\end{proof}

\begin{proof}[Proof of Proposition~\ref{curveisdetermined}]
Suppose that $X$, $P$, $q$, and the $w_i$
are as in the statement of the proposition.
Let $\omega_i$ be the corresponding elements of $H^0(X,\Omega)$.
We will show that $X$ is determined by the $w_i$
when $B = \max\{8g-7,6g+1\}$.

Since $B > 8g-8$,
Lemma~\ref{sections} implies that the $w_i$ determine
the set of homogeneous polynomial relations of degree~$4$
satisfied by the $\omega_i$,
so by Corollary~\ref{superpetri} the $w_i$ determine the image $X'$
of the canonical map.
In particular, the $w_i$ determine whether $X$ is hyperelliptic,
and they determine $X$ if $X$ is nonhyperelliptic.

Therefore it remains to consider the case where $X$ is hyperelliptic.
Applying Gaussian elimination to the $w_i$,
we may assume $0 = \ord_q(w_1) < \dots <  \ord_q(w_g) \leq 2g-2$
and that the first nonzero coefficient of each $w_i$ is $1$.
We use Lemma~\ref{hyperelliptic} repeatedly in what follows.
The value of $\ord(w_2)$ determines whether $P$ is a Weierstrass point.

Suppose that $P$ is a Weierstrass point.
Then $w_i = q^{2i-2} (1 + \cdots + O(q^{B-2i+2}))$,
where each ``$\cdots$'' here and in the rest of this proof
represents some known linear combination of positive powers of $q$
up to but not including the power in the
big-$O$ term.  (``Known'' means ``determined by the original $w_i$.'')

Define $x = w_{g-1}/w_g = q^{-2} (1 + \cdots + O(q^{B-2g+2}))$.
Define $y = dx/(w_g \, dq) = -2 q^{-(2g+1)} (1 + \cdots + O(q^{B-2g+2}))$.
Then $y^2 = 4 q^{-(4g+2)} (1 + \cdots + O(q^{B-2g+2}))$.
Since $B \ge 6g+1$, we have $-(4g+2) + (B-2g+2) > 0$,
so there is a unique polynomial $F$ (of degree $2g+1$)
such that $y^2 = F(x)$.

A similar calculation shows that in the case where $P$ is not
a Weierstrass point, then $B \ge 3g+2$ is enough.
\end{proof}

\begin{remark}
\label{curveisdeterminedoverkbar}
Let us show that if the hypotheses of Proposition~\ref{curveisdetermined}
are satisfied except that the $w_i$ belong to $\kbar[[q]]/(q^B)$
instead of $k[[q]]/(q^B)$,
and the $\omega_i$ are permitted to lie in $H^0(X_\kbar,\Omega)$,
then the conclusion still holds.
Let $E$ be a finite Galois extension of $k$
containing all the coefficients of the $w_i$.
The $E$-span of the $w_i$ must be stable under $\Gal(E/k)$
if they come from a curve over $\Q$,
and in this case, we can replace the $w_i$
by a $k$-rational basis of this span.
Then Proposition~\ref{curveisdetermined} applies.
\end{remark}

\begin{remark}
\label{curveisdeterminedplus}
We can generalize Proposition~\ref{curveisdetermined}
to the case where $q$ is not a uniformizing parameter on $X$:
\begin{quote}
Fix an integer $g \ge 2$, and let
$k$ be a field of characteristic zero.
Let $B>0$ be the integer appearing
in the statement of Proposition~\ref{curveisdetermined}, and let
$e$ be a positive integer.
Then if $w_1$, \dots, $w_g$ are elements of $k[[q]]/(q^{eB})$,
then up to $k$-isomorphism,
there exists at most one curve $X$ over $k$
such that there exist $P \in X(k)$,
an analytic uniformizing parameter $q' \in \Ohat_{X,P}$
and a relation $q' = c_e q^e + c_{e+1} q^{e+1} + \ldots$ with $c_e \neq 0$,
such that $w_1 \, dq$, \dots, $w_g \, dq$
are the expansions modulo $q^{eB}$ of some basis of $H^0(X,\Omega)$.
\end{quote}
The proof of this statement is similar to the proof
of Proposition~\ref{curveisdetermined}, and is left to the reader.
\end{remark}

The rest of this section is concerned with quantitative improvements
to Proposition~\ref{curveisdetermined}, and is not needed for the
general finiteness and computability results of Sections \ref{finitenesssection} and~\ref{computablesection}.

\begin{prop}
\label{curveisdetermined2}
Proposition~\ref{curveisdetermined} holds with $B=B(g)$,
where $B(2)=13$, $B(3)=17$, $B(4)=21$, and $B(g)=6g-5$ for $g \ge 5$.
Moreover, if we are given that the curve $X$ to be recovered
is hyperelliptic, then we can use $B(g)=4g+5$ or $B(g)=2g+4$,
according as $P$ is a Weierstrass point or not.
\end{prop}

\begin{proof}
For nonhyperelliptic curves of genus $g \ge 4$,
we use Theorem~\ref{petri} instead of Corollary~\ref{superpetri}
to see that $B>6g-6$ can be used in place of $B>8g-8$.

Now suppose that $X$ is hyperelliptic.
As before, assume $\ord_q(w_1) < \dots <  \ord_q(w_g)$
and that the first nonzero coefficient of each $w_i$ is $1$.
The value of $\ord(w_2)$ determines whether $P$ is a Weierstrass point.

Suppose that $P$ is a Weierstrass point.
Then $w_i = q^{2i-2} (1 + \cdots + O(q^{B-2i+2}))$.
(As in the proof of Proposition~\ref{curveisdetermined},
$\cdots$ means a linear combination of positive powers of $q$,
whose coefficients are determined by the $w_i$.)
Define $\xtilde = w_{g-1}/w_g = q^{-2} (1 + \cdots + O(q^{B-2g+2}))$.
For $1 \le i \le g-2$, the expression
$\xtilde^{g-i} w_g = q^{2i-2} (1 + \cdots + O(q^{B-2g+2}))$
is the initial expansion of $w_i + \sum_{j=i+1}^g c_{ij} w_j$
for some $c_{ij} \in k$,
and all the $c_{ij}$ are determined if $2 + (B-2g+2) > 2g-2$,
that is, if $B \ge 4g-5$.
Let $w_i'=w_i + \sum_{j=i+1}^g c_{ij} w_j = q^{2i-2}
(1 + \cdots + O(q^{B-2i+2}))$.
Define $x = w_1'/w_2' = q^{-2} (1 + \cdots + O(q^{B-2}))$.
Define $y = -2 q^{-(2g+1)} (1 + \cdots + O(q^{B-2}))$
as the solution to $w_1' \,dq = x^{g-1} \, dx/y$.
Then $y^2 = 4 q^{-(4g+2)} (1 + \cdots + O(q^{B-2}))$,
and if $-(4g+2) + B-2 > 0$,
we can recover the polynomial $F$ of degree $2g+1$
such that $y^2 = F(x)$.
Hence $B \ge 4g+5$ suffices.
A similar proof shows that $B \ge 2g+4$ suffices in the case
that $P$ is not a Weierstrass point.

Hence $\max\{6g-5,4g+5,2g+4\}$ suffices for all types of curves,
except that the $6g-5$ should be $8g-7$ when $g=3$.
This is the function $B(g)$.
\end{proof}

\begin{remark}
We show here that for each $g \ge 2$,
the bound $B=4g+5$ for the precision needed to recover
a hyperelliptic curve is sharp.
Let $F(x) \in \C[x]$ be a monic polynomial of degree $2g+1$
such that $X: y^2=F(x)$ and $X': y^2=F(x)+1$
are curves of genus $g$ that are not birationally equivalent.
Let $q$ be the uniformizing parameter at the point at infinity on $X$
such that $x=q^{-2}$ and $y=q^{-(2g+1)} + O(q^{-2g})$.
Define $q'$ similarly for $X'$.
A calculation shows that the $q$-expansions of the differentials $x^i \, dx/y$
for $0 \le i \le g-1$ are even power series in $q$ times $dq$,
and modulo $q^{4g+4} \,dq$
they agree with the corresponding $q'$-expansions for $X'$
except for the coefficient of $q^{4g+2} \,dq$ in $x^{g-1} \, dx/y$.
By a change of analytic parameter $q=Q+ \alpha Q^{4g+3}$
for some $\alpha \in \C$, on $X$ only,
we can make even that coefficient agree.

A similar proof shows that in the case that $P$ is not a Weierstrass point,
the bound $2g+4$ cannot be improved.
\end{remark}

\begin{remark}
When studying new modular curves of genus $g$,
we can also use the multiplicativity of Fourier coefficients
of modular forms (see~(\ref{eulerproduct}))
to determine some coefficients from earlier ones.
Hence we can sometimes get away with less than $B(g)$ coefficients
of each modular form.
\end{remark}

\subsection{Descending morphisms}
\label{descendingsection}

The next result will be used a number of times throughout this paper.
In particular, it will be an important ingredient 
in the proof of Theorem~\ref{hyperelliptictheorem}.

\begin{prop}
\label{fieldfunction}
Let $X,Y,Z$ be curves over a field $k$ of characteristic zero,
and assume that the genus of $Y$ is $>1$.
Then:
\begin{itemize}
\item[(i)] Given nonconstant morphisms $\pi: X \rightarrow Z$
and $\phi: X \rightarrow Y$
such that $\phi^{*}H^0(Y,\Omega)\subseteq \pi^{*}H^0(Z,\Omega)$,
there exists a nonconstant morphism $u': Z \rightarrow Y$
making the diagram
$$\xymatrix{
X \ar[d]_{\pi} \ar[rd]^{\phi} & \\
Z \ar@{.>}[r]_{u'}         & Y
}$$
commute.

\item [(ii)] If $\pi: X \rightarrow Y$ is a nonconstant morphism
and $u$ is an automorphism of $X$
such that $u^{*}$ maps $\pi^{*}H^0(Y,\Omega)$ into itself,
then there exists a unique automorphism $u'$ of $Y$
making the diagram
$$\xymatrix{
X \ar[d]_{\pi} \ar[r]^u  & X \ar[d]^{\pi} \\
Y \ar@{.>}[r]^{u'}         & Y
}$$
commute.
\end{itemize}
\end{prop}

\begin{proof} The conclusion of~(i)
is equivalent to the inclusion
$\phi^{*} k(Y)\subseteq \pi^{*} k(Z)$.
It suffices to prove that every function in $\phi^* k(Y)$
is expressible as a ratio of pullbacks of meromorphic differentials on $Z$.
If $Y$ is nonhyperelliptic, then
the field $k(Y)$ is generated by ratios of pairs of differentials in
$H^0(Y,\Omega)$, so the inclusion follows from the hypothesis
$\phi^{*}H^0(Y,\Omega)\subseteq \pi^{*}H^0(Z,\Omega)$.
When $Y$ is hyperelliptic, we must modify this argument slightly.
We have $k(Y)=k(x,y)$,
where $y^2=F(x)$ for some polynomial $F(U)$ in $k[U]$ without double
roots. The field generated by ratios of differentials
in $H^0(Y,\Omega)$ is $k(x)$,
so $\phi^* k(x) \subseteq \pi^* k(Z)$.
To show that $\phi^* y \in \pi^*k(Z)$ too,
write $y=x \,dx/(x \, dx/y)$ and observe that $x \, dx/y\in H^0(Y,\Omega)$.

Now we prove~(ii).
The hypothesis on $u^*$ lets us apply~(i) with $\phi=\pi \circ u$
to construct $u': Y \rightarrow Y$.
Since $Y$ has genus $>1$ and $k$ has characteristic zero,
the Hurwitz formula implies that $u'$ is an automorphism.
Considering function fields proves uniqueness.
\end{proof}

\begin{remark}
Both parts of Proposition~\ref{fieldfunction}
can fail if the genus of $Y$ is $1$.
On the other hand, (ii) remains true under the additional assumption
that $X \rightarrow Y$ is optimal in the sense that it does not
factor nontrivally through any other genus~1 curve.
\end{remark}

\begin{remark}
Proposition~\ref{fieldfunction} remains true if $k$
has finite characteristic, provided that one assumes
that the morphisms are separable.
\end{remark}

\section{Some facts about modular curves}
\label{modularcurvessection}

\subsection{Basic facts about $X_1(N)$}
\label{modularformsection}

We record facts about $X_1(N)$ that we will need for the proof
of finiteness in Theorem~\ref{genustheorem}.
See~\cite{shimurabook} for a detailed introduction.

Let $\hh = \{z \in \C : \im\ z > 0\}$.
The group $\SL_2(\R)$ acts on $\hh$ by
$\displaystyle \abcd z = \frac{az+b}{cz+d}$.
The quotient of $\hh$ by the subgroup
        $$\Gamma_1(N) := \left\{ \abcd \in \SL_2(\Z) \Bigm|
                \begin{pmatrix} a \\ c \end{pmatrix}
                \equiv \begin{pmatrix} 1 \\ 0 \end{pmatrix} \pmod{N}
        \right\}$$
is isomorphic as Riemann surface to the space of complex points
on a smooth affine curve $Y_1(N)$ over $\Q$.
The $\Q$-structure on $Y_1(N)$ is uniquely determined by
the condition that its function field be the subset of the function
field of $Y_1(N)$ whose Fourier expansions in $q$ have coefficients in $\Q$.
Next,
there is a uniquely determined smooth projective curve $X_1(N)$ over $\Q$
having $Y_1(N)$ as a dense open subset.
The difference $X_1(N)(\C) \setminus Y_1(N)(\C)$
is the (finite) set of {\em cusps} on $X_1(N)$.
One such cusp is the point $\infty$, which can be defined
as the limit in $X_1(N)(\C)$ as $t \rightarrow +\infty$
of the image of $it \in \hh$ in $X_1(\C)$.
In fact, $\infty \in X_1(N)(\Q)$.
Let $J_1(N)$ denote the jacobian of $X_1(N)$.
We have the Albanese morphism $X_1(N) \rightarrow J_1(N)$
sending $P$ to the class of the divisor $(P)-(\infty)$.

Let $q=e^{2 \pi i z}$.
Pulling back $1$-forms under
$\hh \rightarrow Y_1(N)(\C) \hookrightarrow X_1(N)(\C) \rightarrow J_1(N)(\C)$
identifies $H^0(J_1(N)_\C,\Omega)$ and $H^0(X_1(N)_\C,\Omega)$
with $S_2(N) \frac{dq}{q}$
for some $g$-dimensional $\C$-subspace $S_2(N)$ of $q \C[[q]]$.
We will not define modular forms in general here,
but $S_2(N)$ is known as the space of weight~2 cusp forms on $\Gamma_1(N)$.

If $M|N$ and $d|\frac{N}{M}$, then $z \mapsto d \cdot z$ on $\hh$
induces a morphism $X_1(N) \rightarrow X_1(M)$,
which in turn induces $S_2(M) \rightarrow S_2(N)$ and
$J_1(M) \rightarrow J_1(N)$.
The {\em old subspace} $S_2(N)_\old$ of $S_2(N)$ is defined as the sum
of the images of all such maps $S_2(M) \rightarrow S_2(N)$
for all $d$ and $M$ such that $M|N$, $M \not=N$, and $d|\frac{N}{M}$.
Similarly define the {\em old subvariety} $J_1(N)_\old$ of $J_1(N)$.
The space $S_2(N)$ has a hermitian inner product called the
{\em Petersson inner product}.
Let $S_2(N)_\new$ denote the orthogonal complement to $S_2(N)_\old$ in $S_2(N)$.
The identifications above also give us new and old subspaces
of $H^0(X_1(N)_\C,\Omega)$ and $H^0(J_1(N)_\C,\Omega)$.
Let $J_1(N)^\new = J_1(N)/J_1(N)_\old$.
There is also an abelian subvariety $J_1(N)_\new$ of $J_1(N)$ that can
be characterized in two ways: either as the abelian subvariety
such that
        $$\ker(H^0(J_1(N)_\C,\Omega)
                \rightarrow H^0((J_1(N)_{\new})_{\C},\Omega))
        = H^0(J_1(N)_\C,\Omega)_\old,$$
or as the abelian subvariety such that
$J_1(N) = J_1(N)_\old + J_1(N)_\new$ with $J_1(N)_\old \cap J_1(N)_\new$
finite.
(The latter description uniquely characterizes $J_1(N)_\new$
because of a theorem that no $\Q$-simple quotient of $J_1(N)_\old$
is isogenous to a $\Q$-simple quotient of $J_1(N)^\new$;
this theorem can be proved by comparing conductors, using~\cite{carayol}.)
The abelian varieties $J_1(N)^\new$ and $J_1(N)_\new$ are $\Q$-isogenous.
We define $X_0(N)$, $J_0(N)$, and $J_0(N)^\new$ similarly,
starting with
        $$\Gamma_0(N) :=
        \left\{\, \begin{pmatrix} a& b\\c& d\end{pmatrix}\in \SL_2(\Z)
        \Bigm| c\equiv 0\pmod N \,\right\}$$
instead of $\Gamma_1(N)$.

For $n \ge 1$,
there exist well-known correspondences $T_n$ on $X_1(N)$,
and they induce endomorphisms of $S_2(N)$ and of $J_1(N)$
known as {\em Hecke operators}, also denoted $T_n$.
There is a unique basis $\New_N$ of $S_2(N)_\new$
consisting of $f = a_1 q+a_2 q^2 + a_3 q^3 + \dots$
such that $a_1=1$ and $T_n f = a_n f$ for all $n \ge 1$.
The elements of $\New_N$ are called the {\em newforms of level $N$}.
(For us, newforms are always normalized: this means that $a_1=1$.)
Each $a_n$ is an algebraic integer,
bounded by $\sigma_0(n)\sqrt{n}$ in each archimedean absolute value,
where $\sigma_0(n)$ is the number of positive integer divisors of $n$.
For each field $k$, let $G_k=\Gal(\kbar/k)$.
The Galois group $G_\Q$ acts on $\New_N$.
For any quotient $A$ of $J_1(N)$, let $S_2(A)$ denote the image of
$H^0(A_\C,\Omega) \rightarrow H^0(J_1(N)_\C,\Omega) \isom S_2(N)$
(the last isomorphism drops the $dq/q$);
similarly for any nonconstant morphism $\pi: X_1(N) \rightarrow X$ of curves,
define $S_2(X) := \pi^\ast H^0(X_\C,\Omega) \frac{q}{dq} \subseteq S_2(N)$.
G.~Shimura~\cite[Theorem~1]{shimura1973} associated to $f \in \New_N$
an abelian variety quotient\footnote{Earlier, 
in Theorem~7.14 of~\cite{shimurabook}, Shimura had attached
to $f$ an abelian {\em subvariety} of $J_1(N)$.} 
$A_f$ of $J_1(N)$
such that $S_2(A_f)$ is spanned by the Galois conjugates of $f$;
this association induces a bijection
\begin{align}
\label{J_1(N)quotients}
        \raisebox{-1 ex}{$G_\Q$} \backslash \raisebox{1 ex}{$\New_N$}
        & \stackrel{\sim}\rightarrow
        \frac{\{\text{$\Q$-simple quotients of $J_1(N)^\new$}\}}
                {\text{$\Q$-isogeny}} \\
\nonumber        f & \mapsto A_f.
\end{align}
Shimura proved that $J_1(N)$ is isogenous to a product of these $A_f$,
and K.~Ribet~\cite{ribet1980} proved that $A_f$ is $\Q$-simple.
This explains the surjectivity of~(\ref{J_1(N)quotients}).
The injectivity is well-known to experts, but we could not find
a suitable reference, so we will prove it, 
as part of Proposition~\ref{A_fisogenies}.

The subfield $E_f=\Q(a_2,a_3,\dots)$ of $\C$ is a number field,
and $\dim A_f = [E_f:\Q]$.
Moreover, $\End(A_f)\otimes\Q$ can be canonically identified with $E_f$,
and under this identification the element
$\lambda \in \End A_f$ acts on $f$ as
multiplication by $\lambda$ (considered as element of $E_f$),
and on each Galois conjugate $\ssigma f$ by multiplication by
$\ssigma \lambda$.
(Shimura~\cite[Theorem~1]{shimura1973} constructed an injection
$\End(A_f)\otimes\Q \hookrightarrow E_f$,
and Ribet~\cite[Corollary~4.2]{ribet1980} proved that it was an isomorphism.)

If $A$ and $B$ are abelian varieties over $\Q$,
let $A \Qisog B$ denote the statement
that $A$ and $B$ are isogenous over $\Q$.

\begin{prop}
\label{A_fisogenies}
Suppose $f \in \New_N$ and $f' \in \New_{N'}$.
Then $A_f \Qisog A_{f'}$
if and only if $N=N'$ and $f=\ttau f'$ for some $\tau \in G_\Q$.
\end{prop}

\begin{proof}[Proof (K.~Ribet)]
The ``if'' part is immediate from Shimura's construction.
Therefore it suffices to show that one can recover $f$,
up to Galois conjugacy, from the isogeny class of $A_f$.
Let $\ell$ be a prime.
Let $V$ be the $\Q_\ell$-Tate module $V_\ell(A_f)$ attached to $A_f$.
Let $\Vbar=V \tensor_{\Q_\ell} \Qbar_\ell$.
The proof of Proposition~4.1 of~\cite{ribet1980}
shows that $\Vbar = \bigoplus_\sigma V_\sigma$,
where $V_\sigma$ is an irreducible $\Qbar_\ell[G_\Q]$-module
of dimension~$2$ over $\Qbar_\ell$,
indexed by embeddings $\sigma: E_f \hookrightarrow \Qbar_\ell$.
Moreover, for $p \nmid \ell N$,
the trace of the $p$-Frobenius automorphism acting on $V_\sigma$
equals $\sigma(a_p)$, where $a_p \in E_f$ is the coefficient of $q^p$
in the Fourier expansion of $f$.

If $f' \in \New_N'$ is another weight~2 newform,
and $A_f \Qisog A_{f'}$,
then (using $'$ in the obvious way to denote objects associated with $f'$),
we have isomorphisms of $G_\Q$-modules
$V \isom V'$ and $\Vbar \isom \Vbar'$.
Fix $\sigma: E_f \hookrightarrow \Qbar_\ell$.
Then $V_\sigma$ is isomorphic to some irreducible component
$V'_{\sigma'}$ of $\Vbar'$,
where $\sigma'$ is some embedding $E_{f'} \hookrightarrow \Qbar_\ell$.
Taking traces of Frobenius elements,
we find that $\sigma(a_p) = \sigma'(a'_p)$ for almost all $p$.
Then Theorem~5 of~\cite{li1975} implies that $f$ and $f'$
have the same level and are Galois conjugate.
\end{proof}

We have parallel decompositions
\begin{align*}
        S_2(N)_\new \quad &=
        \Directsum_{f \in G_\Q \backslash \New_N} \quad
        \Directsum_{\tau: E_f \hookrightarrow \C}
                \C \;\; \ttau f \\
        J_1(N)^\new \quad &\Qisog
        \Directsum_{f \in G_\Q \backslash \New_N}
                A_f
\end{align*}
and parallel decompositions
\begin{align}
\label{S_2decomposition}
        S_2(N) \quad &= \quad
        \Directsum_{M|N} \quad
        \Directsum_{f \in G_\Q \backslash \New_M} \quad
        \Directsum_{d|\frac{N}{M}} \quad
        \Directsum_{\tau: E_f \hookrightarrow \C}
                \C \;\; \ttau f(q^d) \\
\label{J_1(N)decomposition}
        J_1(N) \quad &\Qisog \quad
        \Directsum_{M|N} \quad
        \Directsum_{f \in G_\Q \backslash \New_M}
                {A_f}^{n_f}
\end{align}
where $n_f:=\sigma_0(N/M)$.
(When we write $f \in G_\Q \backslash \New_N$,
we mean that $f$ runs through a set of representatives for
the orbits of $G_\Q$ acting on $\New_M$.)
Because of Proposition~\ref{A_fisogenies}, 
the given decompositions of $J_1(N)^\new$ and $J_1(N)$
are exactly the decompositions up to isogeny into
nonisogenous $\Q$-simple abelian varieties occurring with multiplicity.

\begin{lemma}
\label{nonnewlemmaforA}
Let $A$ be an abelian variety quotient of $J_1(N)$.
Then $S_2(A)$ has a $G_\Q$-stable basis
consisting of cusp forms each of the form
\[
        h(q) = \sum_{d\mid \frac{N}{M}} c_d f(q^d)
\]
for some $M\mid N$ and $f\in \New_M$,
where $c_d \in E_f$ depends on $f$ and $d$.
\end{lemma}

\begin{proof}
By multiplying the quotient map $J_1(N) \rightarrow A$
by a positive integer,
we may assume that it factors through the isogeny
        $$J_1(N) \rightarrow \Directsum_{M \mid N} \quad
        \Directsum_{f \in G_\Q \backslash \New_M} A_f^{n_f}$$
of~(\ref{J_1(N)decomposition}).
We may also assume that $A$ is $\Q$-simple,
and even that $A$ is isomorphic to $A_f$ for some $f$,
so that the quotient map $J_1(N) \rightarrow A$
is the composition of $J_1(N) \rightarrow A_f^{n_f}$
with a homomorphism $A_f^{n_f} \rightarrow A$.
The latter is given by an $n_f$-tuple $c=(c_d)$ of elements of $\End(A_f)$,
indexed by the divisors $d$ of $N/M$.
Under
        $$X_1(N) \hookrightarrow J_1(N) \rightarrow A_f^{n_f}
        \stackrel{c}\rightarrow A,$$
the 1-form on $A_\C \isom (A_f)_\C$ corresponding to $f$ pulls back
to $\sum_{d\mid \frac{N}{M}} c_d f(q^d) dq/q$.
Finally, $H^0(A_\C,\Omega)$ has a basis
consisting of this 1-form and its conjugates,
and the pullbacks of these conjugates are sums of the same form.
\end{proof}

\begin{cor}
\label{nonnewcor}
Let $\pi : X_1(N) \to X$ be a nonconstant morphism of curves over $\Q$.
Then $S_2(X)$ has a $G_\Q$-stable basis $T$
consisting of cusp forms each of the form
        $$h(q) = \sum_{d\mid \frac{N}{M}} c_d f(q^d)$$
for some $M\mid N$ and $f\in \New_M$,
where $c_d \in E_f$ depends on $f$ and $d$.
\end{cor}

\begin{proof}
Apply Lemma~\ref{nonnewlemmaforA} to the Albanese homomorphism
$J_1(N) \rightarrow \Jac X$.
\end{proof}

\subsection{Automorphisms of $X_1(N)$}
\label{automorphismsofX_1(N)}

\subsubsection{Diamonds}
The action on $\hh$ of a matrix
$\begin{pmatrix} a& b\\c& d\end{pmatrix} \in \Gamma_0(N)$ induces an
automorphism of $X_1(N)$ over $\Q$ depending only on $(d \bmod N)$.
This automorphism is called the {\em diamond operator} $\langle d\rangle$.
It induces an automorphism of $S_2(N)$.
Let $\varepsilon$ be a Dirichlet character modulo $N$,
that is, a homomorphism $(\Z/N\Z)^* \rightarrow \C^*$.
Let $S_2(N,\varepsilon)$ be the $\C$-vector space
$\{\, h\in S_2 (N) : h \vert \langle d\rangle=\varepsilon (d) h \,\}$.
A form $h\in S_2(N,\varepsilon)$ is
called a form of Nebentypus $\varepsilon$.
Every newform $f\in \New_N$ is a form of some Nebentypus,
and is therefore an eigenvector for all the diamond operators.
Character theory gives a decomposition
        $$ S_2(N)=\bigoplus_\varepsilon S_2(N,\varepsilon)\,,$$
where $\varepsilon$ runs over all Dirichlet characters modulo $N$.
Define $S_2(N,\varepsilon)_\new = S_2(N,\varepsilon) \cap S_2(N)_\new$.
When we write $\varepsilon=1$, we mean that $\varepsilon$ is the
trivial Dirichlet character modulo $N$, that is,
        $$\varepsilon(n)=\begin{cases} 1 & \text{ if $(n,N)=1$} \\
                0 & \text{ otherwise.} \end{cases}$$
A form of Nebentypus $\varepsilon=1$ is a form on $\Gamma_0(N)$. 

We recall some properties of a newform 
$f = \sum_{n=1}^\infty a_n q^n \in S_2(N,\varepsilon)$.
Let $\cond\ \varepsilon$ denote the smallest integer $M \mid N$
such that $\varepsilon$ is a composition
$(\Z/N\Z)^* \rightarrow (\Z/M\Z)^* \rightarrow \C^*$.
Throughout this paragraph, $p$ denotes a prime,
and $\overline{\phantom{X}}$ denotes complex conjugation.
Let $v_p(n)$ denote the $p$-adic valuation on $\Z$.
If $v_p(\cond\ \varepsilon) < v_p(N)$,
then $\varepsilon$ is a composition
$(\Z/N\Z)^* \rightarrow (\Z/(N/p)\Z)^*
        \stackrel{\varepsilon'}\rightarrow \C^*$
for some $\varepsilon'$.
Then
\begin{align}
\label{eulerproduct}
        \sum_{n=1}^\infty a_n n^{-s}
        &= \prod_p (1-a_p p^{-s}+ p \, \varepsilon(p) \, p^{-2s})^{-1}, \\
\label{p|N}
        p \mid N &\iff \varepsilon(p)=0 \\
\label{zeroap}
        v_p(\cond\ \varepsilon) < v_p(N) \ge 2 &\implies a_p=0 \\
\label{littleap}
        v_p(\cond\ \varepsilon) < v_p(N) = 1
                &\implies a_p^2 = \varepsilon'(p) \\
\label{bigap}
        1 \le v_p(\cond\ \varepsilon) = v_p(N) &\implies |a_p|=\sqrt{p} \\
\label{W_p}
        v_p(N)=1, \; \varepsilon=1
                &\implies f\vert W_p = -a_p f \\
\label{complexconjugate}
        p \nmid N &\implies \overline{a}_p = \overline{\varepsilon(p)} \, a_p.
\end{align}
The equivalence~(\ref{p|N}) is trivial.
For the remaining properties see~\cite{li1975} and~\cite{atkin:li78}.

\subsubsection{Involutions}

For every integer $M \mid N$ such that $(M,N/M)=1$, there is an
automorphism $W_M$ of $X_1(N)_\C$
inducing an isomorphism between
$S_2(N,\varepsilon)_{\operatorname{new}}$ and
$S_2(N, \overline{\varepsilon}_M\varepsilon_{N/M})_{\operatorname{new}}$,
where $\varepsilon_M$ denotes
the character mod $M$ induced by $\varepsilon$.
(For more on this action, see~\cite{atkin:li78}).
Given a newform $f \in S_2(N,\varepsilon)$,
there is a newform
$\overline{\varepsilon}_M \tensor f$
in $S_2(N, \overline{\varepsilon}_M\varepsilon_{N/M})$
whose $p^{\operatorname{th}}$ Fourier coefficient $b_p$ is given by
        $$b_p=\begin{cases}
        \overline{\varepsilon}_M(p) a_p & \text{ if $p\nmid M$,}\\
        \overline{\varepsilon}_{N/M}(p) \overline{a_p} & \text{ if $p\mid M.$}
        \end{cases}
        $$
Then
\begin{equation}
\label{W_Maction}
        f\vert W_M=\lambda_M(f)(\overline{\varepsilon}_M \tensor f)
\end{equation}
for some $\lambda_M(f) \in \C$ with $|\lambda_M(f)|=1$.
Moreover, we have
$$
        f\vert W_M^2=\overline{\varepsilon}_{N/M}(-M)\,f \,, \quad
        f\vert (W_{M'}\, W_M)=\overline{\varepsilon}_{M'}(M)f\vert W_{M'\,M}
$$
whenever $(M,M')=1$.  In the particular case $M=N$, the automorphism
 $W_N$ is an involution defined over the cyclotomic field $\Q(\zeta_N)$,
called the {\em Weil involution}.
For all $\tau \in \Gal(\Q (\zeta_N)/\Q)$, the automorphisms ${}^{\tau}W_N$ are  involutions and are also  called Weil involutions by some authors.
The involution $W_N$ maps each $S_2(A_f)$ into itself and satisfies
\begin{equation}\label{weil}
        \langle d \rangle W_N=W_N\langle d\rangle ^{-1}\,,\quad
         {}^{\tau_d}W_N=W_N\, \langle d \rangle \quad\text{ for all }
d\in(\Z/N\Z)^{*}\,,
\end{equation}
where $\tau_d$ is the element of $\Gal (\Q (\zeta_N)/\Q)$
mapping $\zeta_N$ to $\zeta_N^d$.

Let $\varepsilon$ be a Dirichlet character modulo $N$.
Define the congruence subgroup
 $$ \Gamma(N,\eps):=\left\{\begin{pmatrix} a& b\\c & d\end{pmatrix}
\in\Gamma_0(N) \, | \, \varepsilon (d)=1 \right\}\,.
$$
Let $S_2(\Gamma(N,\eps))$ denote the space of weight $2$ cusp forms
on $\Gamma(N,\eps)$.
Denote by $X(N,\varepsilon)$ the modular curve over $\Q$
such that $X(N,\varepsilon)(\C)$
contains $\Gamma(N,\eps)\backslash\hh$
as a dense open subset.
We can identify $H^0(X(N,\varepsilon)_\C,\Omega)$ with
$S_2(\Gamma(N,\eps)) \frac{dq}{q}$.
It is known that
        $$ S_2(\Gamma(N,\eps))=
                \bigoplus_{k=1}^n S_2(N,\varepsilon^k) \,,$$
where $n$ is the order of the Dirichlet character $\varepsilon$.
 The diamond operators and the Weil
 involution induce automorphisms of $X(N,\varepsilon)_\C$.
 If moreover  $\varepsilon=1$, the
 curve $X(N,1)$ is $X_0(N)$
 and the automorphisms $W_M$ on $X_1(N)_\C$
 induce involutions on $X_0(N)$ over $\Q$
 that are usually called the Atkin-Lehner involutions.

\begin{remark}
\label{X^*(N)}
Define the {\em modular automorphism group} of 
$X_1(N)$ to be the subgroup of $\Aut(X_1(N)_\Qbar)$ generated 
by the $W_M$'s and the diamond automorphisms.
The quotient of $X_1(N)$ by this subgroup equals
the quotient of $X_0(N)$ by its group of Atkin-Lehner involutions,
and is denoted $X^*(N)$.
The genus of $X_0(N)$ is $N^{1+o(1)}$ as $N \to \infty$,
so by Proposition~\ref{gonalitybound}, the gonality of $X_0(N)$
is at least $N^{1+o(1)}$.
On the other hand, the degree of $X_0(N) \rightarrow X^*(N)$
is only $2^{\#\{\text{prime factors of $N$}\}} = N^{o(1)}$.
Hence the gonality of $X^*(N)$ tends to infinity.
Thus the genus of $X^*(N)$ tends to infinity.
In particular, Conjecture~\ref{fixedgenusconj} is true if one restricts the
statement to modular curves $X$ that are {\em quotients} of $X_1(N)$ by
some subgroup of the modular automorphism group.
\end{remark}

\subsubsection{Automorphisms of new modular curves}\label{autnew}

By Proposition~\ref{fieldfunction},
if $X$ is a new modular curve of level $N$ and genus $g \ge 2$,
then the diamond operators and the Weil involution $W_N$
induce automorphisms of $X_\Qbar$,
which we continue to denote by $\langle d\rangle$ and $W_N$ respectively.
Throughout the paper, $\D$ will denote the abelian subgroup of $\Aut(X_\Qbar)$
consisting of diamond automorphisms, 
and $\D_N:=\langle \D,W_N \rangle$ will denote 
the subgroup of $\Aut(X_\Qbar)$ generated by $\D$ and $W_N$. 
If moreover $X$ is hyperelliptic, and $w$ is its hyperelliptic involution, 
then the group generated by $\D$ and $W_N.w$  will be denoted by $\D_N'$. 

Note that $\D$ is a subgroup of $\Aut(X)$,  
and the groups $\D_N$, $\D_N'$ are
 $G_\Q$-stable by~(\ref{weil}).
If $\Jac X \Qisog A_f$ 
for some $f$ with nontrivial Nebentypus, then $\D_N$
 is isomorphic to the dihedral group with $2n$ elements, $D_{2 \cdot n}$,
where $n$ is the order of the Nebentypus of $f$. 

\medskip

For every nonconstant morphism $\pi:X_1(N)\rightarrow X$ 
of curves over $\Q$ such that 
$S_2(X) \subseteq \oplus_{i=1}^n S_2(N,\varepsilon^i)$
for some Nebentypus $\varepsilon$ of order $n$,
there exists a nonconstant morphism 
$\pi(\varepsilon): X(N,\varepsilon)\rightarrow X$ over $\Q$.
This is clear when the genus of $X$ is $\le 1$,
and follows from Proposition~\ref{fieldfunction}(i)
if the genus of $X$ is $>1$.
In particular, for a new modular curve $X$ of genus $>1$, 
there exists a surjective morphism $X_0(N)\rightarrow X$ 
if and only if $\D$ is the trivial group.  More generally, we have the
following result.

\begin{lemma} \label{quotientdiamond}
Let $X$ be a new modular curve of level $N$,
and let $G$ be a $G_{\Q}$-stable subgroup of $\Aut(X_\Qbar)$.
Let $X'=X/G$.
If the genus of $X'$ is at least $2$,
then the group $\D'$ of diamonds of $X'$ is isomorphic to $\D/(G\cap \D)$. 
In particular, if $G=\D$ or $\D_N$,
then there is a nonconstant morphism $\pi': X_0(N)\rightarrow X'$ 
defined over $\Q$.
\end{lemma}

\begin{proof}
By Proposition~\ref{fieldfunction},
each diamond automorphism of $X$ induces an automorphism of $X'$.
Hence we obtain a surjective homomorphism $\rho : \D \to \D'$.
Let $K = \ker(\rho)$.
Since $G \cap \D \subseteq K$, 
it suffices to show that $K \subseteq G$.  
Now $H^0(X',\Omega)$ pulls back
to the space $H^0(X,\Omega)^G$ of $G$-invariant regular differentials of $X$,
so $K$ acts trivially on the latter.
Therefore the result follows from Lemma~\ref{quotientlemma} below.
\end{proof}

\begin{lemma}
\label{quotientlemma}
Let $X$ be a curve over a field of
characteristic zero, and let $G$ be a finite subgroup of $\Aut(X)$.
Assume that the genus $g'$ of $X' := X / G$ is at least $2$.
Let 
        $$\overline{G} := \{\, \phi \in \Aut(X) \; : \; \phi^*
\omega = \omega \text{ for all } \omega\in H^0(X,\Omega)^G \,\} .$$  
Then $\overline{G}=G$.
\end{lemma}

\begin{proof}
It is clear that $G \subseteq \overline{G}$.  Now suppose $\phi \in
\overline{G}$, and let $H := \langle G,\phi \rangle$.  Also, set $X''
:= X/H$, and let $g''$ be the genus of $X''$.  Then there is a natural
dominant morphism $\pi : X' \to X''$.  Moreover, 
\[
g'' = \dim H^0(X,\Omega)^H = \dim H^0(X,\Omega)^{\langle
  G,\phi \rangle} = \dim H^0(X,\Omega)^G = g'.
\]
Since $g' \geq 2$ by hypothesis, 
the Hurwitz formula implies $\deg (\pi) = 1$.
Therefore $\phi \in G$, so that $\overline{G} \subseteq G$ as required.  
\end{proof}

\subsection{Supersingular points} \label{supersingularpoints}

We will use a lemma about curves with good reduction.

\begin{lemma}
\label{goodreductionofcurves}
Let $R$ be a discrete valuation ring with fraction field $K$.
Suppose $f : X \to Y$ is a finite morphism of 
smooth, projective, geometrically integral curves over $K$,
and $X$ extends to a smooth projective model $\X$ over $R$
(in this case we say that $X$ has {\em good reduction}).
If $Y$ has genus $\ge 1$,
then $Y$ extends to a smooth projective model $\Y$ over $R$,
and $f$ extends to a finite morphism $\X \to \Y$ over $R$.
\end{lemma}

\begin{proof}
This result is Corollary~4.10 in~\cite{liu-lorenzini1999}.
See the discussion there also for references to
earlier weaker versions.
\end{proof}

The next two lemmas are well-known (to coding theorists, for example),
but we could not find explicit references, so we supply proofs.

\begin{lemma}
\label{weightedcount}
Let $p$ be a prime.
Let $\Gamma \subseteq \SL_2(\Z)$ denote a congruence subgroup
of level $N$ not divisible by $p$.
Let $X_\Gamma$ be the corresponding 
integral smooth projective curve over $\Qbar$,
and let $\psi$ be the degree of the natural map $X_\Gamma \rightarrow X(1)$.
Then $X_\Gamma$ has good reduction at any place above $p$,
and the number of $\Fbar_p$-points on the reduction
mapping to supersingular points on $X(1)$ is at least $(p-1) \psi/12$.
\end{lemma}

\begin{proof}
By~\cite{igusa1959}, the curve $X(N)$ admits a smooth model over $\Z[1/N]$,
and has a rational point (the cusp $\infty$).
Since $p \nmid N$ and $X(N)$ dominates $X_\Gamma$, 
Lemma~\ref{goodreductionofcurves} implies that
$X_\Gamma$ has good reduction at $p$, 
at least if $X_\Gamma$ has genus $\ge 1$.
If $X_\Gamma$ has genus~0, then the rational point on $X(N)$
gives a rational point on $X_\Gamma$, so $X_\Gamma \isom \PP^1$,
so $X_\Gamma$ has good reduction at $p$ in any case.
Replacing $\Gamma$ by the group generated by $\Gamma$ and $-\id$
does not change $X_\Gamma$, so without loss of generality,
we may assume that $-\id \in \Gamma$.
Then $\psi = (\SL_2(\Z):\Gamma)$.

If $E$ is an elliptic curve,
then $\Gamma$ naturally acts on the finite set of
ordered symplectic bases of $E[N]$.
The curve $Y_\Gamma:=X_\Gamma - \{\text{cusps}\}$
classifies isomorphism classes of pairs $(E,L)$,
where $E$ is an elliptic curve
and $L$ is a $\Gamma$-orbit of symplectic bases of $E[N]$.

Fix $E$.
Since $\SL_2(\Z)$ acts transitively on the symplectic bases of $E[N]$,
the number of $\Gamma$-orbits of symplectic bases is 
$(\SL_2(\Z):\Gamma) = \psi$.
Two such orbits $L$ and $L'$ correspond to the same point of $X_\Gamma$
if and only if $L'=\alpha L$ for some $\alpha \in \Aut(E)$.
Then $\psi$ is the sum of the sizes of the orbits of $\Aut(E)$
acting on the $\Gamma$-orbits,
so
        $$\psi = \sum_{(E,L) \in X_\Gamma} \frac{\#\Aut(E)}{\#\Aut(E,L)},$$
where the sum is over representatives $(E,L)$ of the points
on $X_\Gamma$ above $E$, 
and $\Aut(E,L)$ is the subgroup of $\Aut(E)$ stabilizing $L$.

Dividing by $\#\Aut(E)$ and summing over all supersingular $E$ over $\Fbar_p$,
we obtain
        $$\sum_{\text{supersingular points $(E,L) \in X_\Gamma(\Fbar_p)$}} 
        \frac{1}{\#\Aut(E,L)}
        = \psi \sum_{\text{supersingular $E / \Fbar_p$}} \frac{1}{\# \Aut(E)} 
        = \frac{(p-1)\psi}{24},$$
where the last step is the mass formula of Deuring and Eichler
(see Chapter~13, \S4 of~\cite{husemoller1987} for a proof).
But $[-1] \in \Aut(E,L)$, 
so $\#\Aut(E,L) \ge 2$ at each $(E,L) \in X_\Gamma(\Fbar_p)$.
Therefore the number of supersingular points on $X_\Gamma$ must be at least
$2(p-1)\psi/24 = (p-1) \psi/12$.
\end{proof}

\begin{lemma}
\label{supersingularec}
Let  $p$ be a prime.
Given a supersingular elliptic curve $E$ over $\Fbar_p$,
there exists an elliptic curve $E'$ over $\F_{p^2}$ such that
$E \isom E'_{\Fbar_p}$
and the $p^2$-power Frobenius endomorphism of $E'$ equals $-p$.
\end{lemma}

\begin{proof}
Honda-Tate theory supplies an elliptic curve $\calE$ over $\F_p$ 
such that the characteristic polynomial of the $p$-power Frobenius
endomorphism $\Frob_p$ satisfies $\Frob_p^2=-p$.
All supersingular elliptic curves over $\Fbar_p$ are isogenous,
so there exists an isogeny $\phi: {\calE}_{\Fbar_p} \rightarrow E$.
The inseparable part of this isogeny is a power of $\Frob_p$,
so without loss of generality, we may assume that $\phi$ is separable.
The kernel $K$ of $\phi$ is preserved by $-p = \Frob_p^2$,
so $K$ is defined over $\F_{p^2}$.
Take $E' = {\calE}_{\F_{p^2}}/K$.
\end{proof}

The following is a generalization of inequalities used in~\cite{ogg1974}.

\begin{lemma}
\label{ogglemmaX_1}
Let $X$ be a new modular hyperelliptic curve 
of level $N$ and genus $g$ over $\Q$.
If $p$ is a prime not dividing $N$, then $(p-1)(g-1) < 2(p^2+1)$.
\end{lemma}

\begin{proof}
We may assume $g \ge 2$.
Since $p \nmid N$, Lemma~\ref{goodreductionofcurves}
implies that $X_1(N)$ and $X$ have good reduction at $p$,
and the morphism $\pi: X_1(N) \rightarrow X$
induces a corresponding morphism of curves over $\F_p$.
By Proposition~\ref{fieldfunction}(ii),
the diamond automorphism $\langle -p \rangle$ of $X_1(N)$ induces
an automorphism of $X$, which we also call $\langle -p \rangle$.
These automorphisms induces automorphisms of the corresponding
curves over $\F_p$.
For the rest of this proof, $X_1(N)$, $X$, $\pi$, $\langle -p \rangle$
represent these objects over $\F_p$.
Also denote by $\langle -p \rangle$ 
the induced morphism $\PP^1 \rightarrow \PP^1$
on the canonical image of the hyperelliptic curve $X$.
On each curve over $\F_p$ or $\F_{p^2}$, 
let $F$ denote the $p^2$-power Frobenius morphism.

Let $s$ be the number of $x \in X_1(N)(\Fbar_p)$
satisfying $F x = \langle -p \rangle x$.
We will bound $s$ in two different ways.
On the one hand,
Lemma~\ref{supersingularec}
implies that all supersingular points on $X_1(N)(\Fbar_p)$
have this property.
By Lemma~\ref{weightedcount},
there are at least $(p-1) \psi/12$ such points,
where $\psi$ is the degree of $X_1(N) \rightarrow X(1)$.
Thus $s \ge (p-1) \psi/12$.

On the other hand,
any solution $x$ to $F x = \langle -p \rangle x$ in $X_1(N)(\Fbar_p)$
maps to a solution to the same equation on $X$ or on $\PP^1$.
The number of solutions to $F x = \langle -p \rangle x$ in $\PP^1(\Fbar_p)$
is at most $p^2+1$,
because they can be counted by intersecting 
the graphs of $F$ and $\langle -p \rangle$,
which are curves of type $(1,p^2)$ and $(1,1)$, respectively, 
on $\PP^1 \times \PP^1$.
Therefore $s \le (2 \deg \pi)(p^2+1)$.

Combining the inequalities from the previous
two paragraphs gives $(p-1) \psi/12 \le (2 \deg \pi)$ $(p^2+1)$.
By Proposition~1.40 of~\cite{shimurabook},
the genus $g_1$ of $X_1(N)$ satisfies $g_1-1 < \psi/12$.
By Hurwitz, $\deg \pi \le (g_1-1)/(g-1)$.
Combining the last three inequalities yields $(p-1)(g-1) < 2(p^2+1)$.
\end{proof}

\begin{remark}
One can extend the above arguments to prove the following more general
result:
\begin{quote}
Let $X$ be a new modular curve of genus $g$ and level $N$,
with $p \nmid N$.  Suppose $X$ admits a degree $d$ map (defined over
$\Q$) to a curve $X' / \Q$ of genus $g'$.
Finally, suppose that the diamond automorphisms on $X_1(N)$
are compatible with automorphisms of $X'$.
(This is automatic if $g' \geq 2$, or if $g \ge 2$, $g'=0$, and $d=2$.)  
Then $(g-1)(p-1) \leq d (p^2 + 1 + 2pg')$.
\end{quote}
We omit the proof, since we will not use this result.
\end{remark}

\begin{cor}
\label{6and30}
Let $X$ be a new modular hyperelliptic curve 
of level $N$ and genus $g$ over $\Q$.
If $g>10$, then $6|N$.
If $g>13$, then $30|N$.
\end{cor}

\begin{proof}
Take $p=2$, $3$, and $5$ in Lemma~\ref{ogglemmaX_1}.
\end{proof}

In the case where $X$ is a curve dominated by $X_0(N)$,
$\langle -p \rangle$ is the identity,
so in the proof of Lemma~\ref{ogglemmaX_1} 
we may replace the bound $(2 \deg \pi)(p^2+1)$
by $(\deg \pi) \# X(\F_{p^2})$ to obtain the following:

\begin{lemma}
\label{ogglemma}
Suppose that $X$ is a curve of genus $g \ge 2$ over $\Q$,
and $\pi:X_0(N) \rightarrow X$ is a dominant morphism over $\Q$.
Suppose that $p$ is a prime not dividing $N$.
Then $X$ has good reduction at $p$, and $(p-1)(g-1) < \# X(\F_{p^2})$.
In particular, if $G_\Q$ denotes the $\Q$-gonality of $X$, then 
$g < G_\Q \frac{p^2 + 1}{p-1} + 1$.
\end{lemma}

\subsection{Previous work on curves dominated by modular curves}
\label{previouswork}

Except for the determination of all 
new modular genus-$2$ curves with $\Q$-simple jacobian in~\cite{gogo00},
all work we know of on modular curves of genus $\ge 2$
has focused on $X_0(N)$, $X_1(N)$, or
quotients of these by a subgroup of the
group of Atkin-Lehner involutions or diamonds, respectively.
Here we summarize some of this work.

The $19$ values of $N$ for which $X_0(N)$ is hyperelliptic 
(of genus $\ge 2$)
were determined in~\cite{ogg1974}, and all of their equations
were given in~\cite{gonzalez1991}.
The corresponding determination of the three values of $N$ 
for which $X_1(N)$ is hyperelliptic was carried out in~\cite{mestre1981},
and equations for these (and a few other $X_1(N)$) 
were given in~\cite{reichert84}.
In~\cite{ishii:momose1991}
it is proved that 
a curve strictly between $X_0(N)$ and $X_1(N)$
(that is, a quotient of $X_1(N)$ by a nontrivial proper subgroup
of the diamond group)
is never hyperelliptic.
A series of works 
\cite{kluit1979},~\cite{hasegawa:hashimoto96},~\cite{hasegawa97} 
led up to the determination of all $64$ values of $N$
for which the quotient of $X_0(N)$ by its  Atkin-Lehner group, $X^*(N)$, is hyperelliptic,
and this was generalized in~\cite{furumoto:hasegawa99}
to quotients of $X_0(N)$ by an arbitrary subgroup
of the Atkin-Lehner group.
Similar results determining all {\em trigonal} curves of the form
$X_0(N)$, $X_0(N)/W_d$ for a single Atkin-Lehner involution $W_d$,
and $X^*(N)$, can be found in \cite{hasegawa:shimura1999}, 
\cite{hasegawa:shimura1999-2}, and~\cite{hasegawa:shimura2000}, 
respectively.
Some of these curves are not new, and hence do not appear in our tables.

The method of constructing equations for modular curves one at a time
in terms of a basis of cusp forms has been used by 
many authors: in addition to the papers already mentioned,
we have \cite{murabayashi1992}, \cite{shimura1995}, 
 \cite{galbraith1999}, for instance.
In particular, \cite{shimura1995} gives methods 
for both the hyperelliptic and nonhyperelliptic cases,
using Petri's Theorem in the latter, as we do.
F.~Klein~\cite{klein1879} 
gave an explicit model for $X(p)$ for every prime $p$ in 1879.
See~\cite{velu1978} for a modern construction of $X(N)$ for every $N$. 
An analogous result for $X_1(p)$ can be found 
in~\cite{borisov:gunnells:popescu2001}.

\section{Finiteness theorems}
\label{finitenesssection}

Our goal in this section is to prove the finiteness statements
in Theorems \ref{genustheorem}, \ref{nonnewtheorem}, 
and~\ref{gonalitytheorem}.
Algorithmic and practical issues will be dealt with only in later sections,
because we assume that some readers will be interested only
in the finiteness.
In particular, the computability of the sets in question will be proved
in Section~\ref{computablesection}, and practical algorithms in the
hyperelliptic case will be given in Section~\ref{hyperellipticsection}.

\subsection{New modular curves of fixed genus}

\begin{proof}[Proof of finiteness in Theorem~\ref{genustheorem}]
Fix $g \ge 2$.
Let $X$ be a new modular curve of level $N$ and genus $g$,
given by the nonconstant morphism $\pi: X_1(N) \rightarrow X$.
Since the $\Q$-simple factors $A_f$ of $J_1(N)^{\new}$
are pairwise nonisogenous,
any abelian subquotient of $J_1(N)^{\new}$
is isogenous to a product of some subset of these $A_f$.
Applying this to $J = \Jac X$ and then comparing $1$-forms
shows that there exists a $G_\Q$-stable subset $T \subseteq \New_N$
that is a basis for $S_2(X)$.

From now on, let $f = q + a_2 q^2 + a_3 q^3 + \cdots$
denote an element of $T$.
The map $\pi$ is unramified at $\infty$,
since there exists $\omega \in H^0(X_\C,\Omega)$
such that $\pi^* \omega = f \, dq/q$, which is nonvanishing at $\infty$.
In other words, the analytic uniformizer $q$ at $\infty$ on $X_1(N)$
serves also as analytic uniformizer at $\pi(\infty)$ on $X$.

The field $E_f$ generated by the coefficients of $f$
satisfies $[E_f:\Q] = \dim A_f \le \dim J = g$.
Let $B$ be the integer of Proposition~\ref{curveisdetermined}.
Each $a_n$ for $1 \le n \le B$ is an algebraic integer
of degree at most $g$, bounded in every archimedean absolute value
by $\sigma_0(n) \sqrt{n}$,
so there are only finitely many possibilities for $a_n$.
Thus there are only finitely many possibilities for
$\{\, f \bmod q^B \, | \, f \in T \,\}$, given $g$.
By Proposition~\ref{curveisdetermined} and
Remark~\ref{curveisdeterminedoverkbar},
each such possibility arises for at most one curve over $\Q$.
\end{proof}

\begin{remark}
For this proof, any bound on the absolute values of $|a_n|$
in terms of $n$ would have sufficed.
But when we do computations,
it will be useful to have the strong bound $\sigma_0(n) \sqrt{n}$.
\end{remark}

\subsection{Non-new modular curves of fixed genus}
\label{nonnewsection}

Our goal in this section is to prove the
finiteness assertion of Theorem~\ref{nonnewtheorem}.

\begin{proof}[Proof of finiteness in Theorem~\ref{nonnewtheorem}(i)]
Let $X/\Q$ be a modular curve of genus $g\geq 2$ and level $N$,
so that there exists a map $\pi : X_1(N) \to X$ over $\Q$.
Let $B=B(g)$ be the integer of Proposition~\ref{curveisdetermined}.

\medskip

We will follow the basic strategy of the proof of
Theorem~\ref{genustheorem}, but
two complications arise:
the map $\pi$ need not be unramified at $\infty$,
and the cusp forms whose Fourier coefficients we need to bound
need not be newforms.
The sparseness assumption will allow us to get around both of these
difficulties.

\medskip

Let $T$ be as in Corollary~\ref{nonnewcor}.
Choose $j \in T$ with $\ord_q(j)$ minimal.
Then $j \, dq/q$ corresponds to a regular differential on $X_\C$
not vanishing at $\pi(\infty)$.
Since $q$ is an analytic parameter at $\infty$ on $X_1(N)$,
the ramification index $e$ of $\pi$ at $\infty$ equals $\ord_q(j)$.
In particular, $e | N$.
By Remark~\ref{curveisdeterminedplus},
the curve $X$ is determined by
the set $\{\, h \bmod q^{eB} \, | \, h \in T \,\}$.

Fix $h = \sum_{n=1}^\infty a_n q^n \in T$.
Since $\ord_q(h) \ge e$, Corollary~\ref{nonnewcor} implies that
there exist an integer $M|N$,
a newform $f\in \New_M$, and $c_d \in \C$ such that
        $$h = \sum_{d\mid \frac{N}{M}, d \geq e} c_d f(q^d).$$
Since $N\in\Sparse_B$, we have $d > eB$ for all $d\mid N$ such that $d > e$.
If $c_e = 0$, then $h \bmod q^{eB} = 0$;
otherwise we may scale to assume that $c_e = 1$.
In this case, the sparseness of $N$ implies
that $a_n(h(q)) = a_n(f(q^e))$ for $1\leq n < eB$.
In other words, for $1\leq n < eB$ we have
\[
a_n(h) = \begin{cases}
                a_{n/e}(f)  & \text{ if $e \mid n$,}  \\
                0           & \text{ if $e \nmid n$.}
        \end{cases}
\]
Since $f$ is a {\em newform}, each $a_n(h)$
with $n\leq eB$ and $e|n$ is an algebraic integer satisfying
$|a_n(h)| \leq \sigma_0(n/e)\sqrt{n/e}$
for all archimedean absolute values.
As before, it then follows that there are only
finitely many possibilities for
$\{\, h \bmod q^{eB} \, | \, h \in T \,\}$, and therefore for $X$.
\end{proof}

\begin{proof}[Proof of finiteness in Theorem~\ref{nonnewtheorem}(ii)]
Let $X/\Q$ be a modular curve of genus $g\geq 2$ and level $N \in \Smooth_m$,
with $m>0$.  Then $\Jac X$ is isogenous to a subvariety of $J_1(N)$, so
in particular $\Jac X$ has good reduction outside the finite set
$\Sigma$, where $\Sigma$ is the set of primes $p$ such that $p \leq m$.
By combining the Shafarevich conjecture
(proved by Faltings \cite[Theorem~6]{Faltings}) with a well-known
finiteness result for polarizations
(see \cite[Theorem~18.1]{MilneAV}), it follows that
for any number field $K$, any positive integer $d$,
and any finite set $S$ of places of $K$, there are only finitely
many $K$-isomorphism classes of principally polarized
abelian varieties of dimension $d$ over $K$ with good reduction outside
$S$.
In particular, there are only finitely many possible $\Q$-isomorphism
classes for $\Jac X$ as a principally polarized abelian variety.  By the
Torelli theorem \cite[Corollary~12.2]{MilneJV}, it follows that there
are only finitely many possible $\Q$-isomorphism classes for $X$.
\end{proof}

As mentioned in the introduction, part~(iii) of
Theorem~\ref{nonnewtheorem} (concerning prime power levels)
follows from (i) and~(ii).

\medskip

To close this section, we remark that for modular curves of prime
power level (as opposed to general nonnew modular curves), 
one has good control over the ramification at the cusp infinity.  
More precisely:

\begin{prop}
\label{unramifiedp^r}
Suppose that $N=p^r$ is a prime power.
Suppose that $X$ is modular of level $N$
but not modular of level $M$ for any $M<N$.
Then the map $\pi : X_1(N) \to X$ is unramified at $\infty$.
\end{prop}

\begin{proof}
Corollary~\ref{nonnewcor} says that
$S_2(X)$ has a basis in which each element $h$ has the form
$\sum_{d | \frac{N}{M}} c_d f(q^d)$
for some $M|N$ and $f \in \New_M$.
If $\pi$ is ramified at $\infty$, then $c_1=0$ in each such $h$.
Then $\sum_{d | \frac{N}{M}, \, d>1} c_d f(q^{d/p}) \in S_2(p^{r-1})$,
according to the decomposition~(\ref{S_2decomposition}).
Thus $\pi^* H^0(X_\C,\Omega) \subseteq B_p^* H^0(X_1(p^{r-1})_\C,\Omega)$,
where $B_p: X_1(p^r) \to X_1(p^{r-1})$ denotes the degeneracy map
induced by $q \mapsto q^p$.
By Proposition~\ref{fieldfunction}(i),
$\pi$ factors through $X_1(p^{r-1})$,
contradicting the minimality of $N$.
\end{proof}

\subsection{New modular curves of fixed gonality}
\label{gonalitysection} 

In \cite[p.~1006]{Abramovich},
one finds the following lower bound on the gonality of $X_1(N)$.

\begin{theorem}
\label{Abramovichtheorem}
Let $g'$ and $G'$ be the genus and gonality, respectively, of $X_1(N)$.
Then $G' \ge \frac{21}{200} (g'-1)$.
\end{theorem}

The linearity of the bound in the genus of $X_1(N)$
is what enables us to deduce the following.

\begin{prop}
\label{gonalitybound}
If $X$ is a $\C$-modular curve of genus $g \ge 2$ and gonality $G$,
then $g \le \frac{200}{21} G + 1$.
\end{prop}

\begin{proof}
Let $d$ be the degree of the given morphism $\pi: X_1(N)_\C \rightarrow X$.
Let $g'$ and $G'$ be the genus and gonality of $X_1(N)$, so that
$G' \ge \frac{21}{200} (g'-1)$ by Theorem~\ref{Abramovichtheorem}.
Any morphism $X \rightarrow \PP^1_\C$
can be composed with $\pi$ to obtain a morphism
$X_1(N)_\C \rightarrow \PP^1_\C$,
so $G' \le dG$.
The Hurwitz formula implies $g'-1 \ge d(g-1)$.
Now combine these three inequalities.
\end{proof}

Proposition~\ref{gonalitybound} and
Theorem~\ref{genustheorem} together imply the finiteness
in Theorem~\ref{gonalitytheorem}.

\begin{remark}
\label{genus17}
Abramovich's result used the lower bound $0.21$ for the positive
eigenvalues of the Laplacian on $\Gamma\backslash\hh$
for congruence subgroups $\Gamma$.
Kim and Sarnak~\cite{kim-sarnak2000}
recently improved this to $975/4096 > 0.238$.
This means that $200/21$ can be improved to $2/0.238 < 17/2$.
In particular, taking $G=2$, we find that a $\C$-modular
hyperelliptic curve has genus at most $17$.
\end{remark}

\section{Computability}
\label{computablesection}

\subsection{The meaning of computable}

``Computable'' will mean computable by a Turing machine.
(See~\cite{hopcroftullman1979} for a definition of Turing machine.)
To give a precise sense to each statement in our introduction,
we must specify what the input and output of the Turing machine are to be.
In particular, we will need to choose how to represent various objects,
such as curves over number fields.
In many cases, there exist algorithms for converting between
various possible representations, so then
the particular representation chosen is not important.
A number field $k$ can be represented by $f \in \Z[x]$
such that $k \isom \Q[x]/(f(x))$.
The representation is not unique, but we do not care,
since it is well known that a Turing machine can decide,
given $f_1,f_2 \in \Z[x]$, whether $f_1$ and $f_2$ define isomorphic
number fields (and if so, find an isomorphism).
An element $\alpha \in k$ can then be represented by $g \in \Q[x]$
(of degree at most $\deg f-1$)
whose image in $\Q[x]/(f(x))$ corresponds to $\alpha$.
Turing machines can also handle arithmetic over $\Qbar$,
thought of as a subfield of $\C$,
by representing each $\alpha \in \Qbar$ by its minimal polynomial 
over $\Q$,
together with decimal approximations to its real and imaginary parts
to distinguish $\alpha$ from its conjugates.
If $k_0=\Q$ or $k_0=\Qbar$, then a field finitely generated
over $k_0$ can be represented as the fraction field of
a domain $k_0[t_1,\dots,t_n]/(f_1,\dots,f_m)$,
or alternatively as a finite extension of 
the rational function field $k_0(t_1,\dots,t_n)$
(in the same way that we handled finite extensions of $\Q$).

A curve $X$ over a field $k$ finitely generated over $\Q$ or $\Qbar$
can be represented by $f \in k[x,y]$
such that the (possibly singular) affine curve $f(x,y)=0$
is $k$-birational to $X$.
Meromorphic differentials on $X$ represented by $f \in k[x,y]$
can be expressed as $g(x,y) \,dx$ or $g(x,y) \, dy$,
where $g=g_1/g_2$ with $g_1,g_2 \in k[x,y]$ and $g_2$ not divisible by $f$.
A finite set of $k$-isomorphism classes of curves over $k$
can be represented
by a finite list of curves over $k$ in which each class
is represented exactly once.

A closed subvariety of $\PP^n$ 
over a field $k$ finitely generated over $\Q$ or $\Qbar$
can be represented by a finite set of generators of its homogeneous ideal.
A constructible subset of $\PP^n$ (in the Zariski topology)
can be represented as a Boolean combination of closed subvarieties.
Morphisms or rational maps between quasiprojective varieties
can be defined locally by rational functions on a finite number
of affine open subsets.
Elimination of quantifiers over algebraically closed fields is effective, 
and it follows that the image of a constructible subset of a quasiprojective
variety under a morphism can be computed.
(The elimination of quantifiers is usually attributed to Tarski; 
for a proof, see Section 3.2 of~\cite{marker2002}, 
especially Theorem~3.2.2 and Corollary~3.2.8(ii).)
Irreducible components of a variety can also be computed: this follows from
the primary decomposition algorithm in~\cite{hermann1926}.

When in one of our theorems we claim that some set is computable,
what we really mean is that there exists a Turing machine
that takes as input the various parameters on which the set depends
(such as a ground field, a curve, and/or an integer $g$),
and, after a finite but unspecified amount of computation,
terminates and outputs the set in question.

\subsection{Computability lemmas for curves}

\begin{lemma}
\label{birationaldecision}
A Turing machine can solve the following problems:
Given a field $k$ finitely generated over $\Q$ or $\Qbar$,
and given curves $X$ and $Y$ over $k$,
\begin{enumerate}
\item \label{computinggenus}
Compute the genus $g$ of (the smooth projective model of) $X$.
\item \label{decidinghyperelliptic}
If $g \ge 2$, decide whether or not $X$ is hyperelliptic.
\item \label{decidingbirational}
If any of the following hold:
\begin{itemize}
\item[(i)] $g(X) \geq 2$
\item[(ii)] $g(X)=1$ and $g(Y) \neq 1$
\item[(iii)] $g(X)=g(Y)=0$ and $k$ is a number field or $\Qbar$,
\end{itemize}
decide whether $X$ and $Y$ are $k$-birational.
\end{enumerate}
\end{lemma}

\begin{proof}
\nichts

(\ref{computinggenus})
This is well known.
In fact, a genus-computing algorithm has been implemented,
for example in {\sc Magma}~\cite{magma}.

(\ref{decidinghyperelliptic})
Equations for the image $X'$ of
the canonical map $\phi$
can be computed~\cite{magma}.
We must determine whether $\phi:X \rightarrow X'$ is of degree $1$ or $2$.
We can calculate points in $X'(\kbar)$,
for instance by intersecting $X'$ with a hyperplane.
Also, given $P \in X'(\kbar)$, we can compute $\#\phi^{-1}(P)$
by algebra.
In the hyperelliptic case, this number can be $1$ for at most $2g+2$
points $P$.
Therefore, by testing $2g+3$ points in $X'(\kbar)$,
we can distinguish the cases $\deg \phi=1$ and $\deg \phi=2$.

(\ref{decidingbirational})
By~(\ref{computinggenus}) we may assume that $X$ and $Y$ have the
same genus $g$, and that $g$ is known.
Suppose first that $g=0$.
If $k=\Qbar$, then $X$ and $Y$ are automatically birational.
If $k$ is a number field, using {\sc Magma},
for each curve we can compute the equation for the image
of the anticanonical embedding, which is a plane conic.
We may diagonalize the corresponding quadratic forms
to assume that each curve is in the form $x^2 - a y^2 - b z^2=0$
for some $a,b \in k^*$.
Now the curves are isomorphic over $k$
if and only if the Hilbert symbols $(a,b)_v$
match at every place $v$ of bad reduction for the curves
(including the real places).
It is well known that the Hilbert symbols are computable:
they are given by the function \verb+nfhilbert+ in
the number theory package PARI-GP, for instance 
(see the webpage {\tt http://www.parigp-home.de/}).
Thus we are done when $g=0$.

Now suppose $g \ge 2$.
By Riemann-Roch,
the linear system associated to $3K$, where $K$ is a canonical
divisor, induces an embedding $X \rightarrow \PP^{5g-6}$
called the {\em tricanonical embedding}.
{\sc Magma} has routines for computing equations
for the image of this embedding (or of any other morphism associated
to a complete linear system).
Once we have these equations for $X$ and for $Y$,
it is a matter of algebra
to find all linear automorphisms of $\PP^{5g-6}_\kbar$
that map the embedded $X$ into the embedded $Y$.
There will be finitely many, since a curve of genus $\ge 2$
has at most finitely many automorphisms.
We finish by checking them individually to see if any are defined over $k$.
\end{proof}

\begin{remark}
\label{genus1birationaldecision}
Lemma~\ref{birationaldecision}(\ref{decidingbirational})
holds in the genus $1$ case
if and only if there exists an algorithm to decide,
given a genus $1$ curve $Z$ over $k$, whether $Z(k)=\emptyset$.
Such an algorithm exists trivially if $k=\Qbar$.
If $k$ is a number field,
the existence of such an algorithm is implied by the conjecture that
the Shafarevich-Tate group $\Sha(E)$
is finite for all elliptic curves $E$ over $k$.
\end{remark}

\begin{lemma}
\label{ltorsion}
Let $X$ be a curve over a number field $k$, and let $\ell \ge 1$
be an integer.
Let $A=\Jac X$, 
and let $M=A[\ell]$ be the Galois module of $\ell$-torsion points.
Then we can compute a description of $M$ of the following type:
a finite set equipped with an addition table and an action of $\Gal(L/k)$
for some explicit finite Galois extension $L$ of $k$
over which all points of $M$ are defined.
\end{lemma}

\begin{proof}
Compute the genus $g$ of $X$.
We represent points of $A(\Qbar)$ (nonuniquely) by divisors on $X_\Qbar$
whose support avoids any singularities of the given model.
A solution of the ``Riemann-Roch problem''
decides whether a given divisor is principal: see the references
cited in~\cite{poonenants2}.
Hence we can decide when two given divisors represent
the same element of $A(\Qbar)$.
Now enumerate the countably many divisors on $X_\Qbar$ 
avoiding the singularities,
and continue until finding a set $S$ of $\ell^{2g}$ such divisors $D$
such that $\ell D$ is principal, and such that they represent
distinct elements of $A(\Qbar)$.
For each pair in $S$, we can determine which divisor in $S$ is
linearly equivalent to their sum.
Finally, we can compute a finite Galois extension $L$ of $k$ containing
the coordinates of all points appearing in the divisors in $S$,
and for each $\sigma \in \Gal(L/k)$ and $D \in S$,
we can determine which divisor in $S$ is
linearly equivalent to $\ssigma D$.
\end{proof}

\begin{prop}
\label{computingconductor}
The conductor of the jacobian of a curve $X$ over a number field $k$
is computable.
\end{prop}

\begin{proof}
Let $\ell \in \Z$ be a prime, and let $A=\Jac X$.
Compute $A[\ell]$ and $L$ as in Lemma~\ref{ltorsion}.
We can compute the ramification groups
of $\Gal(L/k)$ at each ramified prime of $k$.
Their action on $A[\ell]$ gives us the exponent of $\pp$ in $\cond(A)$
for each $\pp$ not dividing $\ell$.
Repeating the above with a different $\ell$
lets us compute the remaining exponents.
\end{proof}

\subsection{The de~Franchis-Severi Theorem}

The following result, which we have already mentioned,
is known as the de~Franchis-Severi Theorem;
we show in addition that the finite set it promises
is computable.
We thank Matthias Aschenbrenner, Brian Conrad, Tom Graber, Tom Scanlon,
and Jason Starr for discussions related to the proof.

\begin{theorem}[de~Franchis-Severi, computable version]
\label{defranchis}
Let $k$ be a number field or $\Qbar$.
Let $X$ be a curve over $k$.
Then the set of pairs $(Y,\pi)$
where $Y$ is a curve over $k$ of genus at least $2$
and $\pi:X \rightarrow Y$ is a morphism,
up to $k$-isomorphism, is finite and computable.
\end{theorem}

\begin{remark}
We consider $(Y,\pi)$ and $(Y',\pi')$ to be isomorphic
if and only if
there is an isomorphism $Y \to Y'$
whose composition with $\pi$ gives $\pi'$.
\end{remark}

\begin{proof}
For the finiteness, see pp.~223--224 of~\cite{lang1983}.
We assume first that $k=\Qbar$.
The genus $g$ of each $Y$ is bounded by the genus $g_X$ of $X$,
so it suffices to show that for each fixed $g \ge 2$, 
we can compute the set $\Cg$
of isomorphism classes of pairs $(Y,\pi)$ where $Y$ has genus $g$.

View $X$ as a subvariety of $\PP^n$
using the tricanonical embedding.
Compute equations defining $X$ in $\PP^n$.
If $(Y,\pi) \in \Cg$
and $Y \subseteq \PP^m$ is the tricanonical embedding of $Y$, 
then $3K_X-\pi^*(3K_Y)$ is linearly equivalent to an effective divisor,
so the theory of linear systems implies that 
there is a linear subspace $L \subseteq \PP^n$ of dimension $n-m-1$
such that $\pi$ and the 
linear projection $\pi_L: \PP^n \dashrightarrow \PP^m$
coincide on $X-L$.
Riemann-Roch gives $n=5g_X-6$ and $m=5g-6$.
Let $G$ be the Grassmannian variety whose points
correspond to linear subspaces $L \subseteq \PP^n$ of dimension $n-m-1$.
Then $(Y,\pi) \mapsto L$ defines an injection $\iota: \Cg \to G(\Qbar)$.

Conversely, if we start with a linear subspace $L$,
corresponding to a point $s \in G$,
let $Y_s$ be the Zariski closure of $\pi_L(X-L)$ in $\PP^m$,
and let $\pi_s: X \rightarrow Y_s$ denote the morphism induced by
$\pi_L$.  
The map $s \mapsto (Y_s,\pi_s)$ restricted to $\iota(\Cg)$
is an inverse of $\iota$,
but in general $(Y_s,\pi_s)$ need not be in $\Cg$.
Moreover, the $Y_s$ need not form
the fibers of a smooth (or even flat) family.

\medskip

{\em Claim 1:} 
Each closed subvariety $H \subseteq G$ can be (computably) partitioned
into a finite number of irreducible locally closed subsets $H_i$
such that for each $i$, either
\begin{enumerate}
\item For all $s \in H_i$, the curve $Y_s$ is not smooth 
        over the residue field of $s$, or
\item There is a smooth family $\Y \to H_i$ of curves,
        and an $H_i$-morphism $X \times_k H_i \to \Y$ 
        whose fiber above $s \in H_i$ is $\pi_s: X \to Y_s$.
\end{enumerate}

{\em Proof:} 
We use induction on $\dim H$.
Because irreducible components can be computed, we may reduce to
the case where $H$ is irreducible.

We next use the principle that ``whatever happens at the generic point
also happens over some computable dense Zariski open subset.''
Let $\eta$ be the generic point of $H$, and let $L$ be the corresponding
linear subspace defined over the function field $\kappa$ of $H$.
Choose a dense open affine subset $\Spec A$ of $H$,
and write elements of $\kappa$ as ratios of elements of $A$.
Working over $\kappa$, 
we compute the intersection of $L$ with $X$,
the image of $X-L$ under the projection
$\PP^n \dashrightarrow \PP^m$, its closure $Y_\eta$,
and the morphism $\pi_\eta: X \to Y_\eta$.
Using partial derivatives, we compute also 
whether or not $Y_\eta$ is smooth over $\kappa$.
Compute the localization $A'$ of $A$ obtained
by adjoining to $A$ the inverses of the numerators and denominators of the
finitely many elements of $\kappa$ that appear during these computations.
Then the formulas computed over $\kappa$ make sense over $\Spec A'$,
so that all the constructions can be performed over $\Spec A'$.
Moreover, for each $s \in \Spec A'$, 
the curve $Y_s$ is smooth over the residue field of $s$ 
if and only if $Y_\eta$ is smooth over $\kappa$.

Let $H_1=\Spec A' \subseteq H$.
The complement $H - H_1$ is a closed subvariety of lower dimension
than $H$, so using the inductive hypothesis, 
we can partition $H-H_1$ into $H_2$, \dots, $H_n$
with the desired properties.
This completes the proof of Claim 1.

\medskip

Apply Claim 1 to $G$, and discard the $H_i$ in which the $Y_s$ are
not smooth.
For each remaining $H_i$, 
the genus of $Y_s$ is constant for $s \in H_i$,
so we compute the genus of the generic fiber $Y_\eta$
and discard $H_i$ if this genus is not $g$.

\medskip

{\em Claim 2:}
For each remaining $H_i$, let $J_i$ be the set of $s \in H_i$
for which
the linear subspace $L \subseteq \PP^n$ corresponding to $s$
equals the linear subspace $L' \subseteq \PP^n$
defined as the common zeros of the sections in the image of 
$H^0(Y_s,\omega^{\tensor 3}) \hookrightarrow H^0(X,\omega^{\tensor 3})$.
Then $J_i$ is constructible and computable.

{\em Proof:} 
The strategy of proof is the same as that of Claim~1.
The equality of $L$ and $L'$
can be tested at the generic point $\eta$ of $H_i$,
and the outcome will be the same on some computable dense Zariski open
subset of $H_i$.
We finish the proof of Claim 2 by an induction on the dimension.

\medskip

Let $J$ be the (computable) union of the $J_i$.
By definition, $J=\iota(\Cg)$.
Since $\Cg$ is finite, $J$ is finite.
Therefore $\Cg$ can be computed
by computing $(Y_s,\pi_s)$ for $s \in J$.
This completes the proof of Theorem~\ref{defranchis}
in the case where $k=\Qbar$.

\medskip

Finally, suppose that $k$ is a number field.
Compute the finite subset $J$ for $\kbar=\Qbar$ as above.
Since $\iota$ is $G_k$-equivariant,
it suffices to compute the $G_k$-invariant elements of $J$.
For each $L \in J \subseteq G$, 
we compute the Pl\"ucker coordinates for $L$
relative to a $k$-basis of $H^0(X,\omega^{\tensor 3})$,
and discard those for which the Pl\"ucker embedding does not
map $L$ to a $k$-point of projective space.
The $(Y_s,\pi_s)$ corresponding to the remaining  
elements of $J$ are the ones over $k$, each appearing once.
\end{proof}

\begin{remark}
\label{genus1defranchis}
Suppose that $k$ is a number field.
Then the finiteness of the set of $Y$ (without the morphism) 
in Theorem~\ref{defranchis} holds 
even if we include curves $Y$ of genus $\le 1$!
It is also possible to compute a finite set of curves over $k$
representing all the $k$-isomorphism classes of such $Y$,
but with some classes represented more than once.
Eliminating the redundancy in genus~1
would require an algorithm as in Remark~\ref{genus1birationaldecision}.
\end{remark}

\subsection{Computability of modular curves}

Here we use the results of the previous section to
prove that the finite sets in
Theorems \ref{genustheorem}, \ref{gonalitytheorem}, 
and~\ref{pinleveltheorem} are computable.

\begin{prop}
\label{fixedlevel}
Fix $N \ge 1$ and a number field $k$.
The set of $k$-modular curves of level dividing $N$
up to $k$-isomorphism is finite and computable.
\end{prop}

\begin{proof}
By Theorem~\ref{defranchis},
it suffices to compute $X_1(N)$.
First compute the genus $g$ of $X_1(N)$
(a formula can be found in
Example~9.1.6 on page~77 of~\cite{diamond-im1995},
for example).
Using modular symbols we can compute a basis of $S_2(N)$,
with each $q$-expansion computed up to an error of $O(q^{B(g)+1})$.
(This follows from~\cite{manin1972},
and an explicit algorithm can be found in~\cite{cremona1992}.)
Multiplying each basis element by $dq/q$ results in
expansions for a basis of differentials,
and then Section~\ref{recovering} explains how to
recover either an equation $y^2=f(x)$ if $X_1(N)$ is hyperelliptic,
or equations defining the image of the canonical embedding
if $X_1(N)$ is not hyperelliptic.
In the latter case, we can try various linear projections and
use elimination theory until
we find one that yields a plane curve birational to $X_1(N)$:
we can detect whether a linear projection mapped the canonical
model birationally by checking the genus of the image.
Since we can enumerate all linear projections, we will eventually
find one that will work.
(In practice, almost any projection will work.)
\end{proof}

\begin{question}
The de~Franchis-Severi Theorem lets us prove the finiteness
of modular curves of {\em fixed level} whether or not they are new.
Can the proof of the de~Franchis-Severi Theorem
somehow be combined with our proof of Theorem~\ref{genustheorem}
to prove Conjecture~\ref{fixedgenusconj} in general?
\end{question}

\begin{proof}[Proof of computability in
Theorems \ref{genustheorem}, \ref{gonalitytheorem}, and~\ref{pinleveltheorem}]
The proof of finiteness in Theorem~\ref{genustheorem}
produces a finite list of candidate curves $X$.
For each $X$, we compute $N:=\cond(\Jac X)^{1/g}$,
which will be the level of $X$ if $X$ is a new modular curve of
genus $g$ and level $N$.
(If $N$ is not an integer, discard $X$ immediately.)
By Proposition~\ref{fixedlevel},
we can list all modular curves of level $N$.
Lemma~\ref{birationaldecision}(\ref{decidingbirational})
determines if $X$ is birational to a curve in this list; discard $X$ if not.
If $X$ is in this list, then it is modular of genus $g$ and level $N$,
and it must be new, since otherwise $\cond(\Jac X)$ would be less than $N^g$.
Thus we can obtain a list of all new modular curves of genus $g$.
Finally, we can use Lemma~\ref{birationaldecision}(\ref{decidingbirational})
to eliminate redundancy.

Computability in Theorem~\ref{gonalitytheorem} now follows
from computability in Theorem~\ref{genustheorem} and
Proposition~\ref{gonalitybound}.
Computability in Theorem~\ref{pinleveltheorem} follows
from computability in Theorem~\ref{gonalitytheorem} and
Proposition~\ref{pinlevelprop}.
\end{proof}

\begin{remark}
\label{nonnewisnotcomputable}
We now explain the difficulty in proving computability of the sets
in Theorem~\ref{nonnewtheorem}.
In the proof of finiteness in part~(i) (with $S=\Sparse_{B(g)}$), 
we obtained a finite list of candidates $X$.
The problem is that we do not know how to bound
the level of a given nonnew modular curve;
in fact, we do not even know how to test if a curve is modular.
The levels of the modular forms corresponding to the regular
differentials on $X$ can be bounded by the conductor of the jacobian
of $X$, but this is not enough, since the level of $X$ can be
higher than the levels of these forms: see Section~\ref{pathologysection}.

Our proof of Theorem~\ref{nonnewtheorem}(ii) is ineffective,
because there is no known effective proof of the Shafarevich conjecture.
\end{remark}

\section{New modular hyperelliptic curves}
\label{hyperellipticsection}

By Theorem~\ref{gonalitytheorem}, 
the set of the new modular hyperelliptic curves over $\Q$ is finite.
By Remark~\ref{genus17}, the genus of such a curve is $\le 17$.
The goal of this section is to improve this 
by proving Theorem~\ref{hyperelliptictheorem} 
and other restrictions on these curves.
The main results are summarized in Table~\ref{possiblegenera}
and proved in Section~\ref{hypfin}.
The computational results are given in tables in the appendix,
and summarized in Section~\ref{computationalresults}.

\subsection{Criterion to determine new modular hyperelliptic curves}
\label{criterionsection}

In this subsection we will characterize effectively the class of new
modular hyperelliptic curves. 
From now on, for a hyperelliptic curve
$X$ we denote by $\Wei(X)$ the set of Weiertrass points of $X$.

\begin{lemma}\label{criterion1}
Assume that there exists a nonconstant morphism $\pi:X_1(N)_\C
\rightarrow X$ of curves over $\C$ such that $X$ is
hyperelliptic of genus $g$ and
${\pi}^{*}H^0(X,\Omega)=H^0(A_{\C},\Omega)$, for some abelian
variety $A$ over $\Q$ which is a quotient of $J_1(N)^{\new}$.
We denote by $f^{(j)}= \sum_{n\geq 1}{a_n^{(j)}\,q^n}$, $1\leq j\leq g$,
a basis for $S_2(A)$ consisting of elements of $\New_N$,
and we set $P=\pi (\infty)$.
Then:
\begin{enumerate}
\item
There exists a unique basis $\{h_1,\dots ,h_g\}$ of $S_2(A)$
such that for all $1\leq j\leq g$:
$$
\begin{cases} h_j\equiv q^{j}\phantom{i-11} \qquad\qquad\quad\,\,  \pmod{ q^{
g+1}} & \text{ if
$P\not\in \Wei (X)$ ,}\\ h_j\equiv q^{2j-1}+ \sum_{i=j}^{g-1} C_{j, 2
i} \,q^{2 i} \pmod{ q^{ 2g}} & \text{ if $P\in\Wei (X)$.}
\end{cases}
$$
\item Moreover,
\begin{itemize}
\item [(i)]If $P\not\in \Wei (X)$, then
$
 \det( a_i^{(j)})_{1\leq i,j\leq g} \,\,\, \neq 0\,.
$
\item[(ii)] If $P\in \Wei (X)$, then
$
 \det( a_{2i-1}^{(j)})_{1\leq i,j\leq g} \neq 0\,.
$
\end{itemize}

\end{enumerate}

\end{lemma}
\begin{proof}
 Applying Lemma~\ref{hyperelliptic} and the fact that $\pi$ is
 unramified at the cusp $\infty$, we obtain the existence of a basis
 $\{ h_1^{'},\dots , h_g^{'}\}$ of $S_2(A)$ which satisfies:
$$ \begin{cases} h_j^{'}\equiv q^{2j-1} \pmod{ q^{2 j}} & \text{ if
$P\in\Wei (X) \, $ }\\ h_j^{'}\equiv q^{j} \quad\,\, \pmod{ q^{ j+1}}
& \text{ if $P \not\in \Wei (X)$} \end{cases} $$ for all $1\leq j\leq
g$.
Therefore, (1) is obtained by Gaussian elimination.  To obtain~(2),
it suffices to observe that the matrices in (i) and~(ii) are
the change of basis matrices from $\{f^{(j)}\}$ to the basis
$\{ h_j\}$.
\end{proof}

\begin{remark}
Later, in Proposition \ref{a2n=0}, we prove that if $P\in \Wei (X)$
then $a_{2n}^{(j)}=0$ for all $n\geq 1$ and $j\leq g$. In particular,
$4\mid N$ and the basis $h_j$ satisfies
$$ h_j\equiv q^{2j-1} \pmod{ q^{ 2g}} \,.$$
\end{remark}

\begin{remark}
If moreover $X$ comes from a curve over $\Q$
having jacobian isogenous to $A_f$ with
$f=q+\sum_{n\geq 2}a_nq^n\in \New_N$,
then part~(2) of the lemma implies that $\{1,a_2,\dots,a_g\}$
(resp.\ $\{1,a_3,\dots,a_{2g-1}\}$) is a $\Q$-basis of $E_f$
if $P \notin \Wei (X)$ (resp.\ $P\in \Wei (X)$).
\end{remark}

\begin{remark}
If $A$ is a quotient of $J_1(N)$ (over $\Q$),
and the $\C$-vector space $S_2(A)$
has a basis $\{h_1,\dots, h_g\}$ as in part~(1) of the previous lemma,
then each $h_i$ has Fourier coefficients in $\Q$,
since $S_2(A)$ has a basis contained in $\Q[[q]]$.
\end{remark}

The following proposition provides us with an effective criterion to
determine when a $\Q$-factor of $J_1(N)^{\text{new}}$ is
$\Q$-isogenous to the jacobian of a modular hyperelliptic curve of level $N$.
Let $\langle v_1,\dots,v_n \rangle$ denote the span of
elements $v_1,\dots,v_n$ of a vector space.

\begin{prop}
\label{criterion2}

Let $A$ be a quotient of $J_1(N)^{\new}$ of dimension $g>1$.
The following conditions are equivalent:
\begin{enumerate}
\item There exists a modular hyperelliptic curve $X$ of level $N$
over $\Q$ such that $\Jac X \Qisog A$.

\item There exists a hyperelliptic curve $X'$ over $\C$
and a nonconstant morphism $\pi':X_1(N)_\C \rightarrow X'$
such that ${\pi'}^{*}H^0(X',\Omega)=H^0(A_{\C},\Omega)$.

\item  There exists a basis $\{h_1 ,\dots h_g \}$ of $S_2(A)$
as in part~(1) of Lemma~\ref{criterion1}
such that for every pair $g_1,g_2 \in S_2(A)$ satisfying
$\langle g_1,g_2\rangle= \langle h_{g-1},h_g\rangle $
and $g_2\in \langle h_g\rangle $,
there exists $F(U)\in \mathbb{C}[U]$ of degree
$2g+1$ or $2g+2$ without double roots
such that the functions on $X_1(N)$ given by
$$
        x=\frac{g_1}{g_2}\,,\quad  y=\frac{q\, dx/dq}{g_2}
$$
satisfy the equation $y^2= F(x)$.

\end{enumerate}
\end{prop}

\begin{proof} It is clear that (1) implies~(2).
Also, (3) implies~(1),
because when we apply~(3) with $g_1=h_{g-1}$ and $g_2=h_g$,
the modular functions $x$ and $y$ have rational $q$-expansion,
so the corresponding polynomial $F$ has coefficients in $\Q$.

We now assume~(2) and prove~(3).
By Lemma~\ref{criterion1}, there exists such a basis
$\{h_1,\dots,h_g \}$.
As before, we put $P=\pi'(\infty)$.
Let $u$ and $v$ be nonconstant functions on $X'$ such that:
\begin{itemize}
\item $\dv \,u=(Q)+ (w(Q))-(P)-(w(P))$ for some $Q\in X'$ and where
$w$ denotes the hyperelliptic involution of $X'$.
\item $v^2=G(u)$, where $G(U)$ is a polynomial in $\C [U]$ of degree
$2g+1$ or $2g+2$, without double roots.
\end{itemize}
By looking at $\ord_P du/v$ and $\ord_P u du/v$, and using the fact
that $\pi'$ is unramified at $\infty$, we have:
$$ \langle {\pi'}^{*}(du/v),{ \pi '}^{*}(u du/v)\rangle =
\langle  h_{g-1} dq/q, h_{g} dq/q \rangle \,, \quad
\langle { \pi '}^{*}( du/v)\rangle =
\langle  h_g dq/q \rangle\,.
$$
Thus, for every pair $g_1,g_2$ as in part (3), there exists a matrix
$\begin{pmatrix}a & b\\ 0 & d\end{pmatrix}\in \GL_2(\mathbb{C})$
such that
 $$ \left \{
\begin{array}{l}
{\displaystyle g_1 dq/q= a\, {\pi'}^{*}(udu/v)+b\,{\pi '}^{*}( du/v) }\\
{\displaystyle g_2 dq/q=d\, {\pi '}^{*}( du/v)\,.}
\end{array}\right.
$$
Now,  one can  check easily  that the modular functions
$$
x:=\frac{g_1}{g_2}\,,\quad  y:=\frac{q\, dx/dq}{g_2}
$$
satisfy the equation $y^2= F(x)$, where
$$
F(U)=\frac{a^2}{d^4} G\left ( \frac{d U-b}{a}\right )\,.
$$
\end{proof}

In practice, the previous proposition will be used together 
with the next result.

\begin{lemma}\label{boundqexp}
Let $Y$ be a curve of genus $g_Y>0$ over $\C$.
Let $q$ be an analytic uniformizing parameter in $\Ohat_{Y,P}$ 
for a point $P\in Y(\C)$ 
and let $F\in\C[t]$ be a polynomial of degree $d>0$. 
Suppose we are given $\omega_1, \omega_2 \in H^0(Y,\Omega)$ 
such that the functions $x=\omega_1/\omega_2$ and $y=dx/\omega_2$ satisfy 
$$y^2-F(x)\equiv 0 \pmod{ q^c} \quad \text{for some}\quad 
c\geq (2 g_Y-2) \operatorname{Max} \{6,d\}+1\,.
$$
Then we have $y^2=F(x)$.
\end{lemma}

\begin{proof} Put $d'=\operatorname{Max} \{6,d\}$. Now, 
it suffices to observe that $y\, \omega_2^3\in H^0(Y,\Omega^{\otimes 3})$ and, thus, $(y^2- F(x)) \omega_2^{d '}\in H^0(Y,\Omega^{\otimes d'})$ and has at $P$ a zero of order $d'(2g_Y-2)+1$ at  least. 
\end{proof}

The following proposition improves part~(1) of Lemma~\ref{criterion1}
and is useful for computations.

\begin{prop}
\label{a2n=0}
Let $X$ be a new modular hyperelliptic curve of genus $g$ and level $N$
over $\Q$ and let $\pi:X_1(N)\rightarrow X$ be the corresponding morphism.
If $\pi(\infty) \in \Wei(X)$ and
$f=q+\sum_{n\geq 2} a_n q^n \in \New_N\cap S_2(X)$,
then $a_{2n}=0$ for all $n\geq 1$ and $4\mid N$.
\end{prop}

\begin{proof}
Write
$S_2(X) =\oplus_{i=1}^m S_2(A_{f^{(j)}})$,
where $f^{(j)}=q+\sum_{n\geq 2}a_n^{(j)}q^n$ is a
newform in $S_2(N,\varepsilon_j)$.
The condition that $a_{2n}^{(j)}=0$ for all $n$ and $j$ is equivalent to
the condition that $a_2^{(j)}=0$ for all $j$ and $2\mid N$.
We consider three cases.

\vskip 0.2 cm
\noindent
\underline{Case 1}: $g=2$ and $\Jac X$ is $\Q$-simple.

Put
$f:=f^{(1)}=\sum_{n\geq 0} a_nq^n$.
The coefficients of $f$ generate some quadratic field $\Q(\sqrt{d})$.
Let $\sigma$ denote the nontrivial element of $\Gal(\Q(\sqrt{d})/\Q)$,
so that $f^{(2)}=\ssigma f$.
By parts (1) and~(2)(ii) of Lemma~\ref{criterion1}, respectively,
we have $a_2\in \Q$ and $a_3\notin\Q$.
Write $a_n :=A_n + B_n \sqrt{d}$, where $A_n, B_n\in \Q$.
Define $\varepsilon:=\varepsilon_1$.
By~(\ref{eulerproduct}) the $q$-expansion of $f$ is:
$$ q+ A_2 q^2+ (A_3+ B_3 \sqrt d) q^3+ (A_2^2-2 \, \varepsilon (2)) q^4+
(A_5+ B_5 \sqrt d) q^5+ A_2(A_3+ B_3 \sqrt d) q^6+ O(q^7),
$$
where $B_3 \not=0$.
Put $\varepsilon(2)=C+ D\sqrt d$.
Following Lemma~\ref{criterion1} and Proposition~\ref{criterion2},
we compute
$$ h_1=\frac{f+{}^{\sigma} f}{2}\,, \quad h_2= \frac{f-{}^{\sigma}
f}{2 B_3 \sqrt d}\,, \quad
x=\frac{h_1}{h_2}=\frac{1}{q^2}+\dots\,,\quad y=-\frac{q\, dx/dq}{2
h_2}=\frac{1}{q^5}+\dots.
$$
For $F \in \C((q))$, let $\Coeff[q^n,F]$ denote the coefficient
of $q^n$ in $F$.
Since $y^2-x^5$ is at most a quartic polynomial in $x$,
and $x$ is $O(q^{-2})$, we have $\Coeff[q^{-9},y^2-x^5]=0$.
On the other hand, we compute
$\Coeff[q^{-9},y^2-x^5]=-4(A_2 B_3+ D)/B_3$,
so $D=-A_2 B_3$.
Thus $D=0$ if and only if $A_2=0$.
We claim that $A_2=D=0$.
If $2\mid N$ then $\varepsilon(2)=0$, so $D=0$.
If $2 \nmid N$, then $\overline{a}_2=\overline{\varepsilon(2)} a_2$
by~(\ref{complexconjugate}),
so $A_2 = (C \pm D \sqrt{d}) A_2$ and equating coefficients of $\sqrt{d}$
yields $D A_2 = 0$.
Thus $A_2=D=0$ in both cases.

Now $y^2- x^5 - \Coeff[q^{-8},y^2-x^5] x^4$
is at most a cubic polynomial in $x$, so
        $$0 = \Coeff\left[q^{-7},y^2-x^5-\Coeff[q^{-8},y^2-x^5] x^4\right]
                = 12 C.$$
Thus $C=0$, so $\varepsilon(2)=0$.
By~(\ref{eulerproduct}), $a_2=\varepsilon(2)=0$
implies $a_{2n}=0$ for all $n$.
By (\ref{littleap}) and~(\ref{bigap}), we have $4 \mid N$.

\vskip 0.2 cm
\noindent
\underline{Case 2}: $g=2$ and $\Jac X$ is not $\Q$-simple.

Then
$\Jac X \Qisog  A_{f^{(1)}}\times A_{f^{(2)}}$ with $\dim
A_{f^{(1)}}=\dim A_{f^{(2)}}=1$.
Since $\dim A_{f^{(1)}}$ is odd,
the Nebentypus $\varepsilon$ of $f^{(1)}$ is the trivial character
modulo $N$ (cf. the proof of Lemma \ref{evengenus}).
The same holds for $f^{(2)}$.
By part~(1) of Lemma~\ref{criterion1},
$f^{(1)}$ and $f^{(2)}$ share the same coefficient of $q^2$, say $a$.
By~(\ref{eulerproduct}),
\begin{align*}
        f^{(1)} &= q + a q^2 + a_3 q^3 + (a^2-2\,\varepsilon(2)) q^4 + a_5 q^5
                + a \,a_3 q^6 + O(q^7), \\
        f^{(2)} &= q + a q^2 + b_3 q^3 + (a^2-2\,\varepsilon(2)) q^4 + b_5 q^5
                + a \,b_3 q^6 + O(q^7),
\end{align*}
with $a,a_n,b_n\in\Z$.
By part~(2) of Lemma~\ref{criterion1}, $a_3\neq b_3$.
As in Case~1, we compute
$$ h_1=\frac{f^{(1)}+ f^{(2)}}{2}\,, \,\, h_2=
\frac{f^{(1)}-f^{(2)}}{a_3-b_3}\,, \,\,
x=\frac{h_1}{h_2}=\frac{1}{q^2}+\dots\,,\,\, y=-\frac{q\, dx/dq}{2
h_2}=\frac{1}{q^5}+\dots \,\,.
$$
{}From $0=\Coeff[q^{-9},y^2-x^5]=-4a$ we obtain $a=0$.
{}From
        $$0 = \Coeff\left[q^{-7},y^2-x^5-\Coeff[q^{-8},y^2-x^5] x^4\right]
                = 12 \varepsilon(2)$$
we obtain $\varepsilon(2)=0$.
The result follows from $a=\varepsilon(2)=0$, as in Case~1.

\vskip 0.2 cm
\noindent
\underline{Case 3}: $g>2$.

Put
$$ f^{(j)}=\sum_{n\geq 0} a_n^{(j)}q^n\,,\quad E_j=\Q (\{ a_n^{(j)}
\})\,.
$$
Let $\E$ denote the $\Q$-algebra $E_1\times \dots\times E_m$.
Let $\Edual$ denote the dual vector space $\Hom_\Q(\E,\Q)$,
and let $\phi.\ahat$ denote the evaluation
of a functional $\phi \in \Edual$ at an element
$\ahat=(a^{(j)})$ of $\E$.
For $n\geq 1$, set $\ahat_n=(a_n^{(j)}) \in \E$.
Let $\widehat{\varepsilon(2)}=(\varepsilon_i(2)) \in \E$.
By~(\ref{eulerproduct}), we have
\begin{align}
\label{widehat1}
        \ahat_{nm} &= \ahat_n\ahat_m
                \qquad\text{when $(n,m)=1$,} \\
\label{widehat2}
        \ahat_{2^n} &=
                \ahat_{2^{n-1}}\ahat_2-
                2\, \widehat{\varepsilon(2)} \, \ahat_{2^{n-2}}
                \qquad\text{when $n \ge 2$.}
\end{align}
By part~(2)(ii) of Lemma~\ref{criterion1},
$\{\ahat_1,\ahat_3,\ahat_5,\dots,\ahat_{2g-1}\}$
is a $\Q$-basis of $\E$.
Let $\{\phi_1,\phi_2,\dots,\phi_g\}$ denote the dual basis of $\Edual$.
Let $h_i=\sum_{n=1}^\infty (\phi_i.\ahat_n) q^n$ for $1 \le i \le g$.
Then $h_i$ is a $\Q$-combination of conjugates of the $f^{(j)}$,
and hence $\{h_1,\dots,h_g\}$
is the basis of $S_2(X)$
promised by part~(1) of Lemma~\ref{criterion1}.
Thus the $q$-expansion of $h_i$ has the form
$$      h_i = q^{2 i-1}+\sum_{j=i}^{g-1} C_{i, 2 j} \,q^{2 j}
                + \sum_{j=2g}^{\infty} C_{i,j}q^j.$$
In particular, $\phi_i.\ahat_{2}=0$ for $i>1$,
so $\ahat_2 \in \Q \ahat_1 = \Q$.
That is, $a_2^{(j)}$ has a value $a_2 \in\Q$ independent of $j$.
Let $\gamma$ equal $g$ or $g-1$, depending on whether $g$ is odd or even.
Then by~(\ref{widehat1}),
        $$C_{\gamma,2\gamma} = \phi_\gamma . \ahat_{2\gamma}
                = a_2 \phi_\gamma.\ahat_{\gamma} = 0,$$
since $1<\gamma < 2 \gamma-1$.
The same computation shows $C_{1,2}=a_2$.
Define $x=h_{g-1}/h_g$ and $y = \frac{q \,dx/dq}{-2 h_g}$
so that $X$ has the equation $y^2=F(x)$ for some monic $F$ of
degree $2g+1$.
Then
        $$x = q^{-2} + b q^{-1} + O(q^0) \quad\text{ and }\quad
        y = q^{-(2g+1)} + c q^{-2g} + O(q^{-(2g-1)}),$$
for some $b,c \in \Q$.
Also, $y^2=F(x)$ implies $2c=(2g+1)b$.
Moreover $h_i \, dq/q = P_{g-i}(x) \,dx/(-2y)$ for some monic
$P_{g-i}(x) \in \Q[x]$ of degree $g-i$,
and equating coefficients of $q^{2i-1} \, dq$
(just after the monic leading term)
yields
        $$C_{i,2i} = (g-i) b + (b/2) - c = -ib.$$
Setting $i=\gamma$ and $i=1$ yields
$0 = -\gamma b$ and $a_2 = -b$,
so $a_2=b=0$.

Setting $i=2$ yields
$\phi_2.\ahat_4=C_{2,4}=0$.
But $\phi_i \ahat_4=0$ also for $i>2$,
so $\ahat_4 \in \Q$.
By~(\ref{widehat2}), $\widehat{\varepsilon(2)} \in \Q$.
Induction using~(\ref{widehat2}) shows that
$\ahat_{2^n} \in \Q$ for all $n \ge 0$.
By~(\ref{widehat1}), $\ahat_{n} \in \Q \ahat_{\odd(n)}$,
where $\odd(n)$ is the largest odd divisor of $n$.
Using this, $g \ge 3$, and the definition of $\phi_i$,
we now prove that
$C_{g, 2g}=C_{g,2g+2}=C_{g-1,2g-2}=C_{g-1,2g}= 0$.
Indeed,
\begin{equation}
\label{Cformulas}
\begin{array}{rlll}
        C_{g,2g} &= \phi_g.\ahat_{2g} &= 0,
                &\quad\text{since $\odd(2g) < 2g-1$};\\
        C_{g,2g+2} &= \phi_g.\ahat_{2g+2} &= 0,
                &\quad\text{since $\odd(2g+2) < 2g-1$};\\
        C_{g-1,2g-2} &= \phi_{g-1}.\ahat_{2g-2} &= 0,
                &\quad\text{since $\odd(2g-2) < 2g-3$};\\
        C_{g-1,2g} &= \phi_{g-1}.\ahat_{2g} &= 0,
                &\quad\text{since $\odd(2g) < 2g-3$};
\end{array}
\end{equation}
{\em except} that in the last equation, when $g=3$,
$\odd(2g)<2g-3$ fails so we use $\ahat_{2g}=a_2 \ahat_3 = 0$ instead
to deduce the same result.
Therefore
\begin{align*}
        h_{g-1} &= q^{2g-3}(1+O(q^4)), \\
        h_g &= q^{2g-1}(1+ C_{g,2g+1} q^2 + O(q^4)), \\
        x = h_{g-1}/h_g &= q^{-2}(1-C_{g,2g+1} q^2+ O(q^4)),
\intertext{and the new basis defined by $h'_i:=x^{g-i} h_{g}$ for $1 \le i \le g$ satisfies}
        h'_i &= q^{2i-1}(1+c_i q^2 + O(q^4))
\end{align*}
for some $c_i \in \Q$.
The coefficient of $q^4$ in $h'_i$ is zero for all $i$,
so $\ahat_4=0$.
By~(\ref{widehat2}), $\widehat{\varepsilon(2)}=0$.
The result follows from $\ahat_2=\widehat{\varepsilon(2)}=0$, as in Case~1.
\end{proof}

\subsection{Automorphisms of hyperelliptic curves}

Throughout this subsection we make the following assumptions:
\begin{itemize}
\item $k$ is a field of characteristic zero,
\item $X$ is the smooth projective model of a curve
$y^2=F(x)$ for some squarefree $F \in k[x]$
of degree $2g+1$ or $2g+2$ for some $g \ge 2$
(so in particular $X$ is a hyperelliptic curve of genus $g$), and
\item $w$ is the hyperelliptic involution of $X$.
\end{itemize}

\begin{prop} 
\label{hyperellipticisomorphism}
Let $X'$ be a curve of the same type as $X$,
that is, a genus-$g$ curve ${y'}^2=F'(x')$, 
where $F' \in k[U]$ is squarefree.
Every isomorphism $u: X_\kbar \rightarrow X'_\kbar$ is given by an
expression of the following form:
$$ (x,y)\mapsto  \left( \frac{ax+b}{cx+d},
\frac{ e\, y}{\phantom{c}(cx+d)^{g+1}}\right)\,,
$$
for some $M  = \begin{pmatrix} a & b \\ c & d \end{pmatrix} \in \GL_2(\kbar)$
and $e \in \kbar^{*}$. 
The pair $(M,e)$ is unique up to replacement 
by $(\lambda\,M, e\lambda^{g+1})$ for $\lambda\in \overline{k}^{*}$.
If $u$ is defined over $k$,
then one can take $M \in \GL_2(k)$ and $e \in k^*$.
Moreover, if $u':X'_\kbar \rightarrow X''_\kbar$ is another isomorphism,
given by $(M',e')$, 
then the composition $u' u$ is given by $(M'M, e'e)$.
\end{prop}

\begin{proof}
An isomorphism $X \rightarrow X'$ induces an isomorphism
on the canonical images $\PP^1 \rightarrow \PP^1$.
An automorphism of $\PP^1$ over $\kbar$ is given by some $M \in \GL_2(\kbar)$.
The functions $f \in \kbar(X)$ such that $w^* f =-f$ are
those in $\kbar(x) y$, so $y'$ corresponds to $e y/(cx+d)^{g+1}$
for some $e \in \kbar(x)$.
The image of $(x,y) \in X$ must be on $X'$,
so
        $$\left( \frac{e y}{(cx+d)^{g+1}} \right)^2 = F'\left(\frac{ax+b}{cx+d}\right),$$
or equivalently,
        $$e^2 F(x) = (cx+d)^{2g+2} F'\left(\frac{ax+b}{cx+d}\right).$$
The right hand side is a squarefree polynomial,
since it can be obtained from $F'$ by homogenizing to a polynomial
$z^{2g+2} F'(x/z)$ of degree $2g+2$, performing a linear change of
variable on $x$ and $z$, and dehomogenizing.
But $F(x)$ on the left is a squarefree polynomial too.
Thus $e \in \kbar^*$.
The rest of the statements follow easily.
\end{proof}

Recall that $\Aut(X)$ denotes the group of automorphisms of $X$
defined over $k$.

\begin{lemma}
\label{lemma2} 
Let $u \in \Aut(X)$.
We use the $(M,e)$ notation as in Proposition~\ref{hyperellipticisomorphism}.
Then:
\begin{enumerate}
\item If $g$ is even, $u$ is represented 
by a unique pair $(M, e)\in \GL_2( k)\times k^{*}$ 
satisfying $e=(\det M)^{g/2}$.
Then $u \mapsto M$ defines 
an injective homomorphism $\Aut(X) \rightarrow \GL_2(k)$
mapping $w$ to $-\id$.
\item If $g$ is odd, 
$u$ is represented by exactly two pairs 
$(\pm M, e) \in \GL_2(\kbar) \times \kbar^*$ 
satisfying $e^2=(\det M)^{g}$.
Then $u \mapsto e$ defines
a homomorphism $\Aut(X) \rightarrow \mu(k)$
mapping $w$ to $-1$,
where $\mu(k)$ denote the group of roots of unity in $k$.
\end{enumerate}
\end{lemma}

\begin{proof} 
Suppose $g$ is even.
For a fixed pair $(M_0,e_0)\in \GL_2(k)\times k^{*}$, 
the condition $e_0 \lambda^{g+1}=(\det (\lambda M_0))^{g/2}$ 
on $\lambda \in\kbar^*$
is equivalent to 
$\lambda=(\det M_0)^{g/2}/e_0$. 
Thus $\lambda$ is unique and, in particular,  $\lambda \in k^*$.
Uniqueness also implies that $u \mapsto M$ is a homomorphism.
Since $e=(\det M)^{g/2}$, $M$ determines $(M,e)$,
so the homomorphism is injective.

Now suppose $g$ is odd.
If $g$ is odd, 
the condition
$(e_0 \lambda^{g+1})^2=(\det (\lambda M_0))^{g}$ 
determines $\lambda$ up to sign.
It follows that $\lambda^2\in k^{*}$ and that $e:=e_0 \lambda^{g+1}$ is in $k$.
Then $u \mapsto e$ defines a homomorphism $\Aut(X) \rightarrow k^*$.
Its image is contained in $\mu(k)$, since $\Aut(X)$ is finite.
\end{proof}

\begin{lemma}
\label{lemma3}
Suppose $k=\Q$.
If $g$ is even,
then $\Aut(X)$ is isomorphic to a subgroup of 
$D_{2 \cdot 4}$ or $D_{2 \cdot 6}$.
If $g$ is odd, 
then $\Aut(X)$ is the direct product of $\langle w \rangle$
and a subgroup isomorphic to a subgroup of $D_{2 \cdot 4}$ or $D_{2 \cdot 6}$.
In either case, every element of $\Aut(X)$
has order $1$, $2$, $3$, $4$, or $6$.
\end{lemma}

\begin{proof} 
Suppose $g$ is even.
By Lemma~\ref{lemma2}, 
$\Aut(X)$ is isomorphic to a finite subgroup $G$ of $\GL_2(\Q)$.
By averaging an inner product on $\Q^2$,
we obtain a $G$-invariant inner product,
so $G$ can be embedded in the orthogonal group ${\operatorname{O}}_2(\R)$.
All finite subgroups of ${\operatorname{O}}_2(\R)$
are cyclic or dihedral, so it remains to show that elements of $G$
have order dividing $4$ or $6$.
This follows since the eigenvalues of $G$ are roots of unity
in a number field of degree at most $2$.

Now suppose $g$ is odd.
Let $\overline{\Aut(X)}$ denote the image of 
the homomorphism $\Aut(X) \rightarrow \Aut(\PP^1) \isom \PGL_2(\Q)$
mapping $u \in \Aut(X)$ to its action on the canonical image of $X$.
The first map in the exact sequence
        $$0 \rightarrow \langle w \rangle 
        \rightarrow \Aut(X) \rightarrow \overline{\Aut(X)} \rightarrow 0$$
has a section 
$\Aut(X) \rightarrow \mu(\Q) = \{\pm 1\} \isom \langle w \rangle$
mapping $u$ to $e$ as in Lemma~\ref{lemma2}.
Thus $\Aut(X) \isom \langle w \rangle \times \overline{\Aut(X)}$.
Every finite subgroup of $\PGL_2(\Q)$ is isomorphic to
a subgroup of $D_{2 \cdot 4}$ or $D_{2 \cdot 6}$,
by Proposition~A in~\cite{granville1998}.
\end{proof}

\begin{prop}
\label{brandtprop}
Suppose $k=\kbar$.
Then
\begin{enumerate}
\item $\Aut(X)$ does not contain $(\Z/2\Z)^4$.
\item If $\Aut(X)$ contains $(\Z/2\Z)^3$, then $g$ is odd.
\item If $\Aut(X)$ contains $\Z/2\Z \times D_{2 \cdot 4}$, then $4|(g+1)$.
\item If $\Aut(X)$ contains $\Z/2\Z \times D_{2 \cdot 6}$, then $6|(g+1)$.
\item If $\Aut(X)$ contains $D_{2 \cdot 6}$, then $g \not\equiv 1 \pmod{3}$.
\end{enumerate}
\end{prop}

\begin{proof}
This follows from Proposition~2.1 and Satz~5.1 
of~\cite{brandt-stichtenoth1986}.
See~\cite{bujalance1993} for more,
including a complete classification of the groups
that can be the automorphism group of a hyperelliptic curve
of given genus.
\end{proof}

\begin{lemma}
\label{hyperellipticmorphism}
Let $\pi: X \rightarrow X'$ be a degree-$d$ morphism
between curves of genus $g$ and $g'$, respectively,
over a field of characteristic zero.
Assume that $g \ge 2$ and $X$ is hyperelliptic.
Then
\begin{enumerate}
        \item If $g' \ge 2$, then $X'$ is hyperelliptic.
        \item If $g' \ge 2$, then $\pi(\Wei(X)) \subseteq \Wei(X')$.
        \item 
Suppose that $X'=X/\langle \alpha \rangle$
for some $\alpha \in \Aut(X)$ of order $d$.
Assume that the hyperelliptic involution $w$ of $X$ 
is not in $\langle \alpha \rangle$.
Let $g''$ be the genus of $X'' = X/\langle \alpha w\rangle$.
\begin{enumerate}
\item If $d$ is odd, then $g' = \lfloor g/d \rfloor$.
\item If $d$ is even and $g \not\equiv -1 \pmod d$,
then $g' = \lfloor g/d \rfloor$.
\item If $d$ is even and $g \equiv -1 \pmod d$,
then $g'$ and $g''$ equal $\lfloor g/d \rfloor$
and $\lceil g/d \rceil$ in some order.
\end{enumerate}
\end{enumerate}
\end{lemma}

\begin{proof} 
(1) The image $Y$ of the canonical map $X \rightarrow \PP^{g-1}$ has genus zero and 
dominates the corresponding image  for $X'$, which implies that its genus is zero too.  Therefore, $X'$ is hyperelliptic.

(2)
Let $w$ and $w'$ denote the hyperelliptic involutions on $X$ and $X'$,
respectively.
Since the function field $k(X)$ is of degree~2 over $k(Y)$,
it must equal the compositum of $k(X')$ and $k(Y)$ over $k(Y')$,
and the unique nontrivial element $w^*$ of $\Gal(k(X)/k(Y))$
must restrict to the unique nontrivial element $(w')^*$
of $\Gal(k(X')/k(Y'))$.
Thus $\pi w = w' \pi$.
Hence fixed points of $w$ map to fixed points of $w'$.
(Alternatively, (1) and~(2) could
have been deduced from Proposition~\ref{fieldfunction}(ii).)

(3) We may assume $k=\kbar$.
Since $w \not\in \langle \alpha \rangle$,
the automorphisms $\alpha$ and   $w$  induce an automorphism of $Y$
of degree~$d$ and an involution on $X'$ respectively. Let us denote $Y'=X'/\langle w \rangle$, where $w$ stands for  the involution induced on $X'$.
Up to conjugacy, $\Aut(Y) = \PGL_2(k)$
contains only one cyclic subgroup of order $d$.
Thus we may choose coordinates so that 
the morphism $\PP^1 \isom Y \rightarrow Y' \isom \PP^1$
induced by $X \rightarrow X'$
is $x \mapsto x'=x^d$.
In particular, $Y \rightarrow Y'$
is ramified above two points $0, \infty \in Y'$,
each with ramification index $d$.

Let $r' \in \{0,1,2\}$ be the number of points of $Y'$
among $0$ and $\infty$ that ramify in $X' \rightarrow Y'$.
Concerning the behavior of the points in $X$ above $0$ and $\infty$
in $X \rightarrow X'$:
if $d$ is odd, we have $2(2-r')+r'$ points with ramification index $d$;
if $d$ is even, we have $2(2-r')$ points with ramification index $d$,
and $2r'$ points with ramification index $d/2$.
The Hurwitz formula gives
        $$2g-2 = d(2g'-2) + 
        \begin{cases}
                (4-r')(d-1) & \text{ if $d$ is odd} \\
                (4-2r')(d-1) + (2r')(d/2-1) & \text{ if $d$ is even.}
        \end{cases}$$
This implies $g = dg' + (d-1) - r' \lfloor d/2 \rfloor$.
If $d$ is odd, then we deduce $g' = \lfloor g/d \rfloor$.
If $d$ is even and $g \not\equiv -1 \pmod{d}$, then $r'=1$,
and $g' = \lfloor g/d \rfloor$.

Finally suppose that $d$ is even and $g \equiv -1 \pmod{d}$.
Then $r' \in \{0,2\}$.
Since $X$ is birational to $y^2=f(x^d)$ for some polynomial $f$,
and $d$ is even,
the points $0$ and $\infty$ have a total of $4$ preimages in $X$.
Each of these four preimages is fixed by $\alpha$ or $\alpha w$,
but not both.
There are $2(2-r')$ fixed points of $\alpha$,
so there are $2r'$ fixed points of $\alpha w$.
Applying the same arguments to 
the analogous number $r''$ defined for $X''$,
we find $2(2-r'')$ fixed points of $\alpha w$,
so $r'+r''=2$.
Thus either $(r',r'')=(0,2)$, 
in which case $(g',g'')=(\lfloor g/d \rfloor, \lceil g/d \rceil)$,
or $(r',r'')=(2,0)$, 
in which case $(g',g'')=(\lceil g/d \rceil, \lfloor g/d \rfloor)$.
\end{proof}

\begin{remark}
It is not necessarily true that $\pi(\Wei(X)) = \Wei(X')$.
\end{remark}

\subsection{Restrictions on new modular hyperelliptic curves}
\label{hypfin}

\begin{proof}[Proof of part~(ii) of Theorem~\ref{hyperelliptictheorem}]
If $3 \nmid N$, then Lemma~\ref{ogglemmaX_1} implies $g<11$.
Therefore assume $3 \mid N$ from now on.
We are given that $\Jac\, X$ is a quotient of $J_0(N)$.
By Proposition~\ref{fieldfunction}(i),
the morphism $\pi: X_1(N) \rightarrow X$ factors through
a morphism $\pi_0: X_0(N)\rightarrow X$.
Let $\{f^{(1)},\dots,f^{(g)}\}$ be the basis of newforms of $S_2(X)$.
Write $f^{(j)}=\sum_{n\geq 1}{a_n^{(j)}\,q^n}$.
Then $a_1^{(j)}=1$ for all $j$, and $g \ge 3$,
so Lemma~\ref{criterion1}(2) implies that $a_3^{(j)}$ cannot also be
independent of $j$.
In particular there exists $j$ such that $a_3^{(j)} \not=0$.
By~(\ref{zeroap}), it follows that $9\nmid N$.
Hence we have the Atkin-Lehner involution $W_3$ on $X_0(N)$,
which is defined over $\Q$.
By (\ref{littleap}), $a_3^{(j)} \in \{1,-1\}$.
By~(\ref{W_p}), $f^{(j)}\vert W_3= - a_3^{(j)} f^{(j)}$ for all $j$.
Set $X':=X/\langle W_3\rangle$,
where $W_3$ denotes the automorphism
of $X$ induced from $W_3$ on $X_0(N)$ by Proposition~\ref{fieldfunction}(ii).
Then $X'$ also is new of level $N$ and
$S_2(X')$ is spanned by the $f^{(j)}$ with $a_3^{(j)}=-1$.
There is at least one such $j$,
since otherwise $a_3^{(j)}=1$ for all $j$,
contradicting Lemma~\ref{criterion1}(2).
Applying Lemma~\ref{criterion1}(2) to $X'$
shows that the genus $g'$ of $X'$ is $\le 2$.
Similarly, if $X''=X/\langle w W_3\rangle$
then $S_2(X'')$ is spanned by the $f^{(j)}$ with $a_3^{(j)}=1$,
and the genus $g''$ of $X''$ is $\le 2$.
But $g'+g'' = g \ge 3$, so either $g'=2$ or $g''=2$ (maybe both).
Hence $X$ is a new modular hyperelliptic curve of genus $3$ or $4$
having the same level as some new modular curve of genus $2$.
We know the latter levels,
so we can compute all such curves
using the methods to be discussed in Section~\ref{computationalmethods},
one level at a time.
We find only the curve $C_{39}^{A,B}$ in Table~\ref{table5}.
This must be $X_0(39)$,
because $X_0(39)$ is a new modular genus-$3$ curve
of level divisible by $3$. 
\end{proof}

\begin{lemma} \label{evengenus}
Let $X$ be a new modular curve of genus $g$ and let $\D '$ be a subgroup of $\D$.
If the quotient $X':=X/\D'$ has genus $g'$,
then $g-g'$ is even.
In particular, if $X$ is hyperelliptic
and the hyperelliptic involution $w$ belongs to $\D$,
then $g$ is even.
\end{lemma}

\begin{proof}
Assume $\D'\neq\{1\}$. If a newform $f \in S_2(X)$ lies outside the image of $S_2(X')$
then $f$ has nontrivial Nebentypus,
so the number field $E_f$ is a CM field and then $\dim A_f$ is
even. Indeed, it was proved in \cite{ribet77} (see Proposition 3.2,
Proposition 3.3, and the subsequent remark) that given a newform $f$ of weight $k$ and with
nontrivial Nebentypus 
$\varepsilon$, the number field $E_f$ (which must be either a totally
real field or a CM field) is totally real if and only if
$f$ has complex multiplication by an imaginary quadratic field $K$ and
$\varepsilon$ is the quadratic character $\chi$ attached to $K$. 
By  Theorem 3.4 of \cite{ribet77},
if $f$ has complex multiplication by $K$ then the character
$\eta=\varepsilon \cdot \chi$ satisfies $\eta(-1)=(-1)^{k-1}$. 
In particular, if  $k$ is even then the character  $\varepsilon$ is different from $\chi$ and $E_f$ must be a CM  field.

The quotient of $\Jac X$ by the image of $\Jac X'$ 
is isogenous to a product of such $A_f$,
so it has even dimension.
\end{proof}

\begin{prop} \label{propaut}
Let $X$ be a new modular hyperelliptic curve over $\Q$ of genus $g \ge 2$. 
Then $\D$ is cyclic of order $1$, $2$, $3$, $4$ or~$6$.
If $\#\D=4$ or~$6$, then the hyperelliptic involution $w$ is in $\D$. 
If $\#\D=6$, then $g$ is $2$, $12$ or~$14$.
\end{prop}

\begin{proof} 
We may assume $\D \not=\{1\}$.
If $\D$ is cyclic, let $u$ denote a generator.

Since $\langle \D,w \rangle$ is an abelian subgroup of $\Aut(X)$ 
of even order,
Lemma~\ref{lemma3} implies that it is isomorphic to one of the following:
\begin{equation}
\label{Dwpossibilities}
        \Z/2\Z, \quad \Z/4\Z, \quad \Z/6\Z, 
        \quad (\Z/2\Z)^2, \quad \Z/2\Z \times \Z/4\Z, 
        \quad \Z/2\Z \times \Z/6\Z, 
        \quad (\Z/2\Z)^3.
\end{equation}
We label these Cases 1 through~7, respectively.

By~(\ref{weil}), $W_N \not\in \Aut(X)$,
and the group $\D_N'' = \langle \D,w,W_N \rangle$ is a semidirect product 
\hbox{$\langle \D,w \rangle \rtimes \langle W_N \rangle$} 
with $W_N$ acting on the normal subgroup $\langle \D,w \rangle$ as $-1$.
Hence in Cases 4, 5, 6,~7,
the group $\D_N''$ is isomorphic to
        $$(\Z/2\Z)^3, \quad \Z/2\Z \times D_{2 \cdot 4}, 
        \quad \Z/2\Z \times D_{2 \cdot 6}, 
        \quad (\Z/2\Z)^4,$$
respectively.
Case~7 is impossible by Proposition~\ref{brandtprop}(1).

In Cases 4, 5,~6,
Proposition~\ref{brandtprop} implies $2|(g+1)$, $4|(g+1)$, $6|(g+1)$, 
respectively.
In particular, $g$ is odd,
so Lemma~\ref{evengenus} implies that $w \not\in \D$,
so $\D$ is isomorphic to $\Z/2\Z$, $\Z/4\Z$, $\Z/6\Z$, respectively.
This proves the first statement of Proposition~\ref{propaut}.
To prove the second statement, it suffices to rule out Cases 5 and~6.

\medskip
\noindent \underline{Case~5}: $\D \isom \Z/4\Z$, $w \not\in \D$.

Combining $4|(g+1)$ with the inequality $g \le 17$ of Remark~\ref{genus17}
shows that $g$ is $3$, $7$, $11$, or $15$.
Lemma~\ref{evengenus} implies that the genus $g'$ of $X':=X/\D$ is odd.
Hence Lemma~\ref{hyperellipticmorphism}(3) implies that 
$(g,g')$ is $(3,1)$, $(7,1)$, $(11,3)$, or $(15,3)$.
In the last two cases (with $g'=3$), parts (1) and (3) of
Lemma~\ref{hyperellipticmorphism} imply that 
either $X''=X'/W_N = X/\D_N$ or $X'''=X'/\langle W_n w \rangle = X/\D_N'$
has genus~2 and is dominated by $X_0(N)$.
In the first two cases,
Lemma~\ref{hyperellipticmorphism}(3) implies that 
the curve $X/\langle u^2 w \rangle$ (resp.\ $X/\langle u w \rangle$)
has genus~2 and has a nontrivial diamond of order~4 (resp.\ 2) 
by Lemma~\ref{quotientdiamond} (the nontrivial diamond prevents
the genus from being $1$).
Using the methods of Section~\ref{computationalmethods} 
to check all levels where there are such new genus-$2$ curves, 
we find that no $X$ in Case~5 exists.

\medskip
\noindent \underline{Case~6}: $\D \isom \Z/6\Z$, $w \not\in \D$

Combining $6|(g+1)$ with $g \le 17$ shows that
$g$ is $5$, $11$, or $17$.
In the first (resp.\ second or third) case,
Lemma~\ref{hyperellipticmorphism}(3) and Lemma~\ref{evengenus} 
together imply that 
the curve $X/\langle u^3 w \rangle$ (resp.\ $X/\langle u w \rangle$) 
has genus~2, and 
by Lemma~\ref{quotientdiamond} 
it has a nontrivial diamond of order~6 (resp.~2).
We check as before that such curves do not exist.

\medskip

To finish the proof of Proposition~\ref{propaut},
we must show that if $\D \isom \Z/6\Z$ and $w \in \D$,
then $g$ is $2$, $12$ or~$14$.
By Lemmas \ref{quotientdiamond} and~\ref{evengenus},
$X'':=X/\langle u^2 \rangle$ has diamond group $\langle w \rangle$
and its genus $g''$ is even.
On the other hand, $g$ is even, $g \le 17$,
and $g''=\lfloor g/3 \rfloor$ by Lemma~\ref{hyperellipticmorphism}(3),
so $(g,g'')$ is one of $(2,0)$, $(6,2)$, $(8,2)$, $(12,4)$, or $(14,4)$.
We rule out $(6,2)$ and $(8,2)$ by checking all levels
of genus~2 curves with diamond group of order~2.
\end{proof}

\begin{remark} 
Assume $g\geq 3$ and $\D = \langle u \rangle$.
If $\#\D=2$ and $w\notin \D$ (resp.\ if $\#\D=6$),
then 
$X/\langle u.w\rangle$ (resp.\ $X/\langle u^2 \rangle$)
is a new modular hyperelliptic curve of genus $\ge 2$
of the same level,
and by Lemma~\ref{quotientdiamond}
its group of diamonds has $2$ elements and contains $w$.
Therefore, in order to find all levels of 
new modular hyperelliptic curves with $\#\D>1$, 
it suffices to determine the cases such that 
$\# \D=2,3$ or $4$ with the additional requirement
that $w \in \D$ when $\# \D$ is even.
\end{remark}

\begin{prop} \label{D=2}
Let $X$ be a new modular hyperelliptic curve over $\Q$ of genus $g \ge 2$
such that $2 \mid \#\D$ and $w \notin\D$.
Then $\#\D=2$, and $g$ is $3$, $7$, $8$, or $9$.
If $g=3$, then $X$ is one of the seven curves in Table~\ref{table9}.
If $g>3$, then $3 \nmid N$.
\end{prop} 

\begin{proof}
Proposition~\ref{propaut} implies $\#\D=2$.
Let $u$ be the generator of $\D$.
Let $X'=X/\D$ and $X''=X/\langle uw \rangle$,
and let $g'$ and $g''$ be their genera, respectively.
By Lemma~\ref{hyperellipticmorphism}(3), $g'$ and $g''$
equal $\lfloor g/2 \rfloor$ or $\lceil g/2 \rceil$ in some order.
Moreover $g''=g-g'$ is even (by Lemma~\ref{evengenus})
and positive (since $X \not\isom X/\D$).
In particular $g \not\equiv 2 \pmod{4}$, so $g \not\in \{2,6,10,14\}$.

Suppose $g''=2$.
For each new modular curve of genus~2 with diamond group
of order~2, we compute all new modular curves at that level;
we find only the seven curves in Table~\ref{table9},
all of genus~$3$.
Therefore we may assume $g'' \not=2$ from now on,
so $g'' \ge 4$.
Since $|g'-g''| \le 1$,
we have $g' \ge 3$, and $g \ge 7$.

\medskip
\noindent\underline{Case 1}: $3 \mid N$.

Since $X/\D$ is of genus $g' \ge 3$ with trivial diamond group,
Theorem~\ref{hyperelliptictheorem}(ii) implies $X/\D \isom X_0(39)$,
so $N=39$.
Computations show that $X_0(39)$ is the only new modular
hyperelliptic curve of level $39$ and genus $\ge 3$.

\medskip
\noindent\underline{Case 2}: $3 \nmid N$.

Then $g \le 10$ by Corollary~\ref{6and30}.
The possibility $g=10$ was already ruled out.
\end{proof}

\begin{prop}
\label{D=3}
Let $X$ be a new modular hyperelliptic curve over $\Q$ of genus $g \geq 3$
such that $\#\D=3$.
Then $g=3$ or $g=5$. 
If $g=3$, then $X$ is one of the five curves in Table~\ref{table8}.
If $g=5$ there are at least two such curves, 
which are given in Table~\ref{table10}.
\end{prop}

\begin{proof}
Since $\Aut (X_{\overline{\Q}})$ contains 
$\D_N'' := \langle \D,w,W_N \rangle \isom D_{2 \cdot 6}$,
Proposition~\ref{brandtprop} implies $g \not\equiv 1 \pmod{3}$.
Also $g \le 17$, so $g$ must be 
$17$, $15$, $14$, $12$, $11$, $9$, $8$, $6$, $5$ or $3$.

Suppose $g>5$.
Then we claim that there is a quotient $X'$ of $X$,
such that $X'$ is a new modular curve of genus~2, of the same level $N$, 
and dominated by $X_0(N)$.
Indeed, by repeated application of Lemma~\ref{hyperellipticmorphism}(3),
at least one of the curves $X/\D_N$, $X/\D_N'$ or $X/\D$ 
has genus $2$ in each possible case. 
We now consult the list of new modular curves of genus-2
to check that $X$ does not exist.

Finally we consider the case $g=3$. 
By Lemma~\ref{evengenus}, we have $\Jac X\Qisog A_f\times A_h$, 
where $\dim A_f=2$ and $\dim A_h=1$. 
If ${}^\sigma\! f$ denotes the nontrivial Galois conjugate of $f$,
then $\{f,{}^\sigma\! f,h\}$ is a basis of eigenvectors
of a generator $u$ of $\D$ acting on $S_2(X)$.
Since the $\D$-invariant subspace of $S_2(X)$ is $1$-dimensional,
corresponding to $A_h$, and since $u$ is defined over $\Q$,
the eigenvalues must be $\{\zeta,\zeta^2,1\}$, respectively,
where $\zeta$ is a primitive cube root of $1$.
In particular, the basis of eigenvectors is unique
up to scalar multiples.
Moreover, the field $E_f$ contains the eigenvalue of $u$ acting on $f$,
so $E_f = \Q(\zeta)$.

On the other hand,
all elements of order~3 in $\PGL_2(\C)$ are conjugate,
so we can choose a coordinate function $x$ on 
$\PP^1_\C = X_\C/\langle w \rangle$
so that $u$ induces $x \mapsto \zeta x$.
Then $X_\C$ is the smooth projective model of $y^2=F(x)$
for some squarefree polynomial $F$,
and $\{ dx/y, x\,dx/y, x^2\,dx/y\}$
is a basis of eigenvectors of $u$ acting on 
$H^0(X,\Omega) = S_2(X) \frac{dq}{q}$. In particular, 
there exists a unique eigenvector
whose square equals the product of the other two, 
at least up to a constant factor $c$.

The same must be true of $\{f,{}^\sigma\! f,h\}$ and 
the corresponding constant factor $c$ is in $\Q(\zeta)$. 
For such a relation to be consistent with the action of $\sigma$, 
it must be that $h^2 = c f {}^\sigma\! f$.
Comparing leading coefficients of $q$-expansions shows that $c=1$.
Thus $h^2 = f {}^\sigma\! f$.

Using the methods of Section~\ref{computationalmethods}
we run through all possibilities
with $M=18$, $f=\sum_{i=1}^{M} a_n q^n+O(q^{M+1})$ 
with $a_n \in \Z[\zeta]$,
and $h=\sum_{i=1}^{M} b_n q^n+O(q^{M+1})$ with $b_n\in\Z$
satisfying $h^2 = f {}^\sigma\! f$.
We obtained only the five curves of Table~\ref{table8}.
\end{proof}

\begin{proof}[Proof of parts (i) and~(iii) of Theorem~\ref{hyperelliptictheorem}]
Part~(i) follows from parts (ii) and~(iii)
together with Remark~\ref{genus17}.
Part~(ii) was already proved, so we need only prove~(iii).
Suppose $\#\D > 1$, and $g$ is odd.
Lemma~\ref{evengenus} implies $w \not\in \D$,
and Proposition~\ref{propaut} then implies that $\#\D$ is $2$ or $3$.
Propositions \ref{D=2} and~\ref{D=3}
complete the proof in these two cases, respectively.
\end{proof}

\begin{remark} 
In order to make the process of discarding levels easier 
in the previous propositions, we used some additional information,
such as Corollary~\ref{6and30}. 
In cases where we know that there are 
quotient curves of some smaller genus $\ge 2$,
and we have already determined 
all the hyperelliptic curves of that level and genus,
we have used this information.
\end{remark}

We summarize most of the results of this section 
in Table~\ref{possiblegenera}.
For each possibility for $\D$ and for each possible answer to the question
``Is $w \in \D$?'', we give all integers $\ge 2$
which {\em might} 
be the genus of a new modular hyperelliptic curve over $\Q$
having that $\D$ and such a $w$.
A {\bf bold} number indicates that
our computations have found a curve of that genus.
The other numbers might not actually occur;
in fact, most of them probably do not.

\begin{table}[ht]
\begin{tabular}{|c|c||c|}
      $\D$  & $w \in \D$? & Potential values of $g$ \\
\hline
\hline
$\{1\}$ & no & $\mathbf 2$, $\mathbf 3$, $\mathbf 4$, $\mathbf 5$, $\mathbf 6$, $7$, $8$, $9$, $10$ \\
\hline
\multirow{2}{10mm}{$\Z/2\Z$} & yes & $\mathbf 2$, $\mathbf 4$, $6$, $8$, $10$, $12$, $14$, $16$ \\
      & no & $\mathbf 3$, $7$, $8$, $9$ \\
\hline
$\Z/3\Z$ & no & $\mathbf 3$, $\mathbf 5$ \\
\hline
$\Z/4\Z$ & yes & $\mathbf 2$, $\mathbf 4$, $6$, $8$, $10$, $12$, $14$, $16$ \\
\hline
$\Z/6\Z$ & yes & $\mathbf 2$, $12$, $14$ \\
\hline
\end{tabular}
\medskip
\caption{Possibilities for the diamond group and genus of a new modular hyperelliptic curve.  See Theorem~\ref{hyperelliptictheorem} and Propositions \ref{propaut}, \ref{D=2}, and~\ref{D=3}.}
\label{possiblegenera}
\end{table}

\subsection{Computational methods}
\label{computationalmethods}

Recall that~\cite{gogo00} computed all new modular genus-two curves
with $\Q$-simple jacobian.
Using similar reasoning, and using some of the sieves described there,
we compute all new modular genus-two curves with
jacobian {\em not} $\Q$-simple: there are exactly $64$
of such curves (see Table~\ref{table1}).
In principle, we can also compute
all the equations, levels and newforms for
new modular hyperelliptic curves of genus $g$
for each $g>2$.
But the enormous number of possibilities for the coefficients $a_p$
prevents us from completing these computations in practice.

We demonstrate the method in the case of
new modular hyperelliptic genus-$3$ curves with $\Q$-simple jacobian.
Let $f$ be a corresponding newform.
Then $\dim A_f$ is odd, so $f$ has trivial Nebentypus $\varepsilon$.
The only subfields of $E_f$ are $\Q$ and $E_f$ itself,
so part~(2) of Lemma~\ref{criterion1} implies:
$$
E_f=\begin{cases} \Q( a_2)=\Q (a_3)\,,& \text{ if $\pi (\infty)\not\in\Wei
(X)$,}\\ \Q( a_3)=\Q (a_5)\quad (a_{2n}=0)\,, & \text{ otherwise.}
\end{cases}
$$
For simplicity, we outline the computation
in the case where $a_2 \not=0$
(or equivalently, $\pi(\infty) \not\in \Wei(X)$).
In this case, $E_f=\Q(a_2)$, so $a_2 \not\in \{-1,0,1\}$,
and (\ref{zeroap}) and~(\ref{littleap}) imply that $2 \nmid N$.
\begin{enumerate}
\item
We determine all possible polynomials
$H_2(x)=\prod_{i=1}^3(x-{}^{\sigma_i}a_2)$, that is,
all the monic irreducible cubic polynomials in $\Z[x]$
such that all zeros are real and of absolute value $\le 2 \sqrt{2}$.
In total, there are $80$ such polynomials.

\item For each $H_2(x)$, we fix a zero $a_2$.
Let $M=2g+5=11$,
which is the bound in Proposition~\ref{curveisdetermined2}
plus $1$ (since we multiply newforms by $dq/q$ instead of $dq$).
For every prime $p$ satisfying $3 \le p \le M$,
the possible values of $\varepsilon(p)$ and $a_p$
are those such that
$\varepsilon(p) \in \{0,1\}$ and $a_p$ is an algebraic integer in
$\Q(a_2)$ with  $|a_p| \leq 2\sqrt p$
with respect to every archimedean absolute value.
We may restrict the possibilities by imposing
$\Q(a_3)=\Q(a_2)$ and $\varepsilon(2)=\varepsilon(3)=1$.

\item For each possible $f=q+\sum_{n=2}^M a_n\,q^n +O(q^{M+1})$,
we write $f= g_0 + a_2 g_1 + a_2^2 g_2$ with $g_i \in \Q[[q]]$.
The $g_i$ have the same span as the conjugates of $f$.
Applying linear algebra to the $g_i$,
we compute the basis $\{h_1, h_2, h_3\}$ of part~(1) of
Lemma~\ref{criterion1}.
Next we compute $\xtilde=h_2/h_3$ and $h_1'=\xtilde^2 h_3$.

\item We impose the condition $h_1'\in \langle h_1, h_2,h_3\rangle$
to sieve out some possibilities.
Also we use this condition to extend the precision of $h_1'$,
to determine the first $M$ coefficients of $h_1'$.
Next compute $x=h_2/h_1'$ and $y= (q \, dx/dq)/h_1'$.
\item We impose the condition that $y^2=F(x)$ for a
polynomial $F$ of degree $8$ without double roots.
In the cases that survive, we compute $F(x)$.
\item
Now we have a list of candidate curves.
In principle, we should compute the conductor $\mathcal N$
of each jacobian and keep only those such that
$N := {\mathcal N}^{1/3} \in \Z$
and there exists a newform in $S_2(X_0(N))$
giving rise to $y^2=F(x)$.
In practice, since computing the conductor can be difficult,
it is easier to try to recognize the candidates in a list of
new modular curves calculated from newforms of small level
(as discussed later in this subsection)
and hope that all candidates show up.
\end{enumerate}

\begin{remark}
If $X$ is a new modular curve of level $N$,
and $p$ is a prime not dividing $N$,
then Eichler-Shimura theory shows that
the product $P(x)$ of $x^2-a_p x + p \varepsilon(p)$
over all $f \in S_2(X) \cap \New_N$
must equal the characteristic polynomial of Frobenius
acting on the Tate module of the jacobian of the
mod $p$ reduction of $X$.

In the calculation above, we know $2 \nmid N$,
so we can compute the characteristic polynomials
for all genus-$3$ hyperelliptic curves over $\F_2$,
in order to restrict the possibilities for $H_2(x)$
in step~(1).
Moreover, we can restrict attention to the curves
such that $\#X(\F_2) \ge 1$ (because of $\pi(\infty)$)
and such that $\#X(\F_4) > 2$ (by Lemma~\ref{ogglemma}).

Similar ideas restrict the possibilities for $a_3$,
even when $\pi(\infty) \in \Wei(X)$, i.e., when $N$ is even.
For larger $g$, if it becomes too time-consuming to list all
genus-$g$ hyperelliptic curves over $\F_3$, 
one can at least rule out certain $a_3$ 
by translating the following into conditions on $a_3$:
\begin{itemize}
\item $\#X(\F_3) \ge 1$, because of the image of the cusp $\infty$,
\item $\#X(\F_3) \le 8$ and $\#X(\F_9) \le 20$, since $X$ is hyperelliptic,
\item $\#X(\F_3) \le \#X(\F_9)$,
\item $\#X(\F_9) > 2g-2$, by Lemma~\ref{ogglemma}, 
if $\Jac\, X$ is a quotient of $J_0(N)^\new$.
\end{itemize}
\end{remark}

\medskip

The computation in steps (1) through~(6)
(of new modular hyperelliptic genus-$3$ curves with $a_2 \not=0$
and with $\Q$-simple jacobian)
shows that only two such curves exist.
These are given in Table~\ref{table2} as $C_{41}^A$ and $C_{95}^A$.

In the $a_2=0$ case,
the number of possibilities to analyze
is considerably higher,
because there are more possible values for $a_3$,
hence more possibilities for $E_f$.
Furthermore, in this case, Proposition~\ref{curveisdetermined2}
requires knowledge of $a_p$ for primes $p \le 4g+6$
(which is $18$ for $g=3$)
in order to determine the relation $y^2=F(x)$.

\medskip

Performing all these computations would be extremely time-consuming.
Therefore instead we conducted a search of all new modular hyperelliptic curves $X$ of some bounded level. We used W.~A.~Stein's Modular Symbols package to implement a program in {\sc Magma}~\cite{magma} that detects whether
a set of newforms corresponds to a new modular hyperelliptic curve $X$
and that computes an equation for this curve if so. 
This program, based on Proposition~\ref{criterion2},
was used to determine the $\ge$ entries in Table~\ref{counts}.

\begin{remark}
If there is a nonconstant morphism
$\pi (\varepsilon):X(N,\varepsilon)\rightarrow X$,
then we compute the bound $c$ as in Lemma~\ref{boundqexp}  replacing $g_Y$ by  the genus of $X(N,\varepsilon)$.
If moreover $\varepsilon=1$ and $\Jac\, X$ is
$\Q$-simple,
then for each Atkin-Lehner involution $W_M$ we have
${}^{\sigma}f\vert W_M=\lambda(M) {}^{\sigma}f$ for all $\sigma$,
where $\lambda(M)\in \{-1,1 \}$.
Let $B'(N) = \{W_M\in B(N) \, | \, \lambda(M)=1\}$,
where $B(N)$ is the group of Atkin-Lehner involutions.
Then $X$ is dominated by $X_0(N)/B'(N)$
and the bound $c$ is computed taking $g_Y$ as the 
 genus of $X_0(N)/B'(N)$.
Note that $B'(N)=B(N)$ or $B(N)/B'(N)\simeq \Z/2\Z$,
so $B'(N)=\{ \id\}$ if and only if
$N$ is a power of a prime with $\lambda(N)=-1$.
\end{remark}

\subsection{Computational results}
\label{computationalresults}

Table~\ref{counts} summarizes the results of the computations.
Each entry indicates the number of new modular hyperelliptic curves
over $\Q$ with the genus prescribed by the column,
and with jacobian satisfying the conditions prescribed by the row heading.
In the first (resp.\ second, third, fourth) row,
a number following $\ge$ is the number of curves of that type
of level $\le 3000$ (resp.\ $569$, $2000$, $569$).
We do not know if there are others of higher level.
The zeros in the second row for odd $g$ are from Lemma~\ref{evengenus}.


\begin{table}[ht]
$$
\begin{array}{c||c|c|c|c|c|c}
        & g=2 & g=3 & g=4 & g=5 & g=6 & g \ge 7 \\
\hline
\hline
\text{$\Q$-simple, $\#\D=1$} & 120 \,\, & \ge 14 & \ge 13 & \ge 3 & \ge 0
& \ge 0 \\
\hline
\text{$\Q$-simple, $\#\D>1$} &\phantom{0} 29 \,\, & \phantom{0}0 \,\, & \ge 1 & 0 & \ge 0 & \ge 0 \\
\hline
\text{not $\Q$-simple, $\#\D=1$} & \phantom{0}64 \,\, & \ge 32 & \ge 7& \ge 0 & \ge 1
& \ge 0\\
\hline
\text{not $\Q$-simple, $\#\D>1$} &\phantom{00} 0\,\,  & \ge 12 \,\,  & \geq 3 & \geq 2 & \geq 0 &\geq  0 \\
\hline
\hline
\text{Total} & 213 \,\, & \ge 58 & \ge 23 & \ge 5 & \ge 1 & \ge 0 \\
\hline
\end{array}
$$
\caption{The number of new modular hyperelliptic curves of genus $g$.}
\label{counts}
\end{table}

Equations for the $120+29$ genus-two curves with $\Q$-simple jacobian
have appeared in~\cite{gogo00}.
Equations for all the rest of the curves
are in the tables of the appendix of this paper.

Of the $120$ genus-two curves counted in row~1, only six have level $>3000$: 
their levels are $3159$, $4160$, $7280$, $7424$ (twice), and $7664$. 
Of the $64$ genus-two curves counted in row~3, only two have level $>2000$: 
both are of level $2208$. 
Similarly, of the $29$ genus-two curves counted in rows~2 and 4, 
only four have level $>569$: 
their levels are $768$ (twice), and $928$ (twice).  Moreover, for $g \ge 3$, the largest level we found 
in the curves of row~1 (resp.\ row~3, resp.\ rows~2 and 4) 
was $1664$ (resp.\ $944$, resp.\ $512$ even though we computed up to level $3000$ (resp.\ $2000$, resp.\ $569$). 
This leads us to believe that the lower bounds in Table~\ref{counts} 
are close to the exact numbers.

Finally, we note that all levels obtained when $g>2$ have at most two
different odd prime divisors and only the levels $N=36$, $72$, $144$, $784$ 
are divisible by the square of an odd prime.

\section{New modular curves with trivial character}
\label{trivialcharsection}

\subsection{A finiteness result for curves with $p$ in the level}
\label{pinlevelsection}

The goal of this section is to prove the finiteness statement
in Theorem~\ref{pinleveltheorem}.  
We first give some equivalent characterizations 
of what it means for a curve to have ``trivial character''.

\begin{lemma}
Let $X$ be a new modular curve of genus $g\geq 2$, and let $N$ be a
positive integer.  The following are equivalent:
\begin{itemize}
\item[(i)] There exists a morphism
$\pi_0 : X_0(N) \to X$ with $S_2(X) \subseteq S_2(X_0(N))^{\new}$.
\item[(ii)] There exists a morphism
$\pi_1 : X_1(N) \to X$ with $S_2(X) \subseteq S_2(X_0(N))^{\new}$.
\item[(iii)] The jacobian $J$ of $X$ is a quotient of $J_0(N)^{\new}$.
\end{itemize}
We say that a new modular curve $X$ has {\it trivial character} (and
level $N$) if any of the above conditions are satisfied.
\end{lemma}

\begin{proof}
It is clear that (i) implies (iii).  Furthermore, it follows from
Lemma~\ref{fieldfunction}(i) that (ii) implies (i).  
It remains to show that (iii) implies (ii).  
Since $X$ is a new modular curve,
there exists an integer $M$
and a morphism $\pi : X_1(M) \to X$ 
with $S_2(X) \subseteq S_2(X_1(M)^\new)$.
Then $J$ is $\Q$-isogenous to a product of
abelian varieties of the form $A_f$ with each $f$ a newform of level $M$.  
By assumption $J$ is a quotient of $J_0(N)^{\new}$, 
so $J$ is $\Q$-isogenous also to a
product of abelian varieties of the form $A_{f'}$ with each $f'$ a newform of
level $N$ and trivial Nebentypus character.  
Applying Proposition~\ref{A_fisogenies},
we deduce that $N=M$ and $S_2(X) \subseteq S_2(X_0(N))^\new$.
\end{proof}

The following lemma will be used in the proof of
Proposition~\ref{pinlevelprop} below, and will be used
repeatedly in Section~\ref{smallpinlevelsection}.

\begin{lemma}
\label{RRlemma}
Suppose $X$ is a curve of genus $g>0$ over a field $k$ and that $q$
is an analytic local uniformizing parameter at the point $P\in X(k)$. 
Let $\omega_1,\ldots,\omega_g$ be a basis for $H^0(X,\Omega)$ with 
$\omega_i = (\sum_{j=1}^\infty a_j^{(i)}q^j)\frac{dq}{q}$
and $a_1^{(i)} = 1$ for all $i$.  
Let $m \geq 2$ be an integer, 
and suppose that $a_m^{(i)} = a_m^{(i')}$ for all $1 \leq i,i' \leq g$.
Then there exists a rational function $f$ of degree $m$, defined over $k$,
with poles only at $P$.
\end{lemma}

\begin{proof}
As usual, if $D$ is a divisor on $X$, 
let $l(D)$ denote the dimension of the vector space
$L(D)$ consisting of $0$ and 
the rational functions on $X$ whose divisor is $\ge -D$.
Let $K$ denote a canonical divisor on $X$.  

By hypothesis, if $\omega = (\sum_{j=1}^\infty a_j q^j)\frac{dq}{q}$
is any linear combination of $\omega_1,\ldots,\omega_g$ 
with $a_1 = 0$, then $a_m = 0$. 
Therefore no regular differential on $X$ vanishes to order
exactly $m - 1$ at $P$.  In other words, we have $l(K-(m-1)P) = l(K-mP)$, 
which by Riemann-Roch is equivalent to $l(mP) - l((m-1)P) = 1$.  
Any $f \in L(mP) - L((m-1)P)$ works.
\end{proof}

\begin{cor}
\label{RRcor}
With the same notation as in Lemma~\ref{RRlemma}, we have:
\begin{itemize}
\item[(i)] If $a_2^{(i)} = a_2^{(i')}$ for all $1 \leq i,i' \leq g$,
  then either $g = 1$, 
  or $X$ is hyperelliptic and $P$ is a Weierstrass point.
\item[(ii)] If $a_2^{(i)} = a_2^{(i')}$ and $a_3^{(i)} = a_3^{(i')}$
  for all $1 \leq i,i' \leq g$, then $g = 1$.
\item[(iii)] If $g \geq 2$ and every differential in $H^0(X,\Omega)$ 
  vanishing at $P$ vanishes to order at least $r$ at $P$, 
  then $r \leq 2$, 
  and $r = 2$ if and only if $X$ is hyperelliptic and 
  $P$ is a Weierstrass point.
\end{itemize}
\end{cor}

\begin{proof}
Part~(i) follows from Lemma~\ref{RRlemma},
and (iii) follows from (i) and~(ii).
Therefore we need only prove~(ii).

The argument of Lemma~\ref{RRlemma}
yields rational functions in $L(2P)-L(P)$ and in $L(3P)-L(2P)$.  
Taking products, we find rational functions in $L(mP)-L((m-1)P)$
for all $m \ge 2$.
By induction, we prove $l(mP) \ge m$ for $m \ge 1$.
If we take $m$ large, Riemann-Roch implies $g \le 1$.
\end{proof}

\begin{prop}
\label{pinlevelprop}
Let $X$ be a new modular curve of level $N$ and trivial character.
If $p$ is a prime divisor of $N$, then the $\Q$-gonality
of $X$ is at most $p^2$.
\end{prop}

\begin{proof}
The value of $a_{p^2}$ is the same for each newform
$f = \sum_{n \ge 1} a_n q^n$ in $\New_N$:
it is $0$ or $1$, depending on whether $p^2 \mid N$ or not.
(This follows from the more general 
statements (\ref{eulerproduct}), (\ref{p|N}), 
(\ref{zeroap}), and~(\ref{littleap})
by taking $\varepsilon$ to be the trivial Dirichlet character modulo $N$.)
Also, recall from the proof of Theorem~\ref{genustheorem} that $q$
serves as an analytic uniformizing parameter at the image $P$ of the cusp $\infty$
under $\pi: X_1(N) \rightarrow X$.
The result therefore follows from Lemma~\ref{RRlemma}.
\end{proof}

\begin{proof}[Proof of Theorem~\ref{pinleveltheorem}]
Combine Proposition~\ref{pinlevelprop}
with Theorem~\ref{gonalitytheorem}.
\end{proof}

\subsection{Curves with level divisible by small primes}
\label{smallpinlevelsection}

Theorem~\ref{pinleveltheorem} shows that there are only
finitely many new modular curves $X$ with trivial character and level
divisible by a given prime number $p$.  
In this section, we prove some further results about the levels of new modular
curves with trivial character.  

In addition to relying heavily on Lemma~\ref{RRlemma} and
Corollary~\ref{RRcor}, 
we make use of a classical lemma of Castelnuovo and
Severi.  Before stating it, we make the following definition.
If $f_1:C \to C_1$ and $f_2: C \to C_2$ 
are nonconstant morphisms (of degree $d_1,d_2$, respectively) 
between curves, we say that
$f_1$ and $f_2$ are {\em independent} if 
the product morphism $(f_1,f_2): C \rightarrow C_1 \times C_2$
maps $C$ birationally to its image $C'$.
A consideration of degrees shows that if $\gcd(d_1,d_2)=1$, 
then $f_1$ and $f_2$ are automatically independent.

\begin{lemma}
\label{CSlemma}
Let $C_1,C_2,C$ be curves of genera $g_1,g_2$, and $g$, respectively, 
over the field $k$ of characteristic zero.  
Let $f_i : C \to C_i$ be morphisms of degree $d_i$,
$i=1,2$.  Assume that $f_1$ and $f_2$ are independent.
Then
\[
g\leq (d_1-1)(d_2-1) + d_1g_1 + d_2g_2.
\]
\end{lemma}

\begin{proof}
See Theorem~3.5 of~\cite{accola1994}.  
We remark that Lemma~\ref{CSlemma} is true over any field $k$,
but the proof in \cite{accola1994} assumes that $k$ has
characteristic zero.
\end{proof}

We now investigate some restrictions on the powers of 2 and 3 
dividing the level of new modular curves with trivial character.  

\begin{prop}
\label{2and3levelprop}
Let $X$ be a new modular curve over $\Q$ of genus
$g\geq 2$, level $N$, and trivial character.  Then:
\begin{itemize}
\item [(i)] If $2 \mid N$, then $X$ admits a map of degree 2, defined
  over $\Q$, to a curve of genus at most 1, and $g \leq 16$.
\item[(ii)] $4\mid N$ if and only if $X$ is hyperelliptic and
  $\pi(\infty)$ is a Weierstrass point.
\item[(iii)] If $6 \mid N$, then $g \leq 5$.
\item[(iv)] If $12 \mid N$, then $g=2$.
\item[(v)] If $18 \mid N$, then $g \leq 4 $.
\item[(vi)] $36 \nmid N$.
\end{itemize}
\end{prop}

\begin{proof}

Let $f_1,\ldots,f_g$ be a
basis of newforms for $\Gamma_0(N)$ spanning $S_2(X)$.  

To prove (ii), first note that, as in the proof of Theorem~\ref{genustheorem}, the map $\pi$ is unramified at
$\infty$, so that $q$ serves as an analytic uniformizing parameter at
$P := \pi(\infty)$.  
If $4 \mid N$, then $a_2(f_i) = 0$ for
all $i$ by (\ref{zeroap}), so it follows from Corollary~\ref{RRcor}
that $X$ is hyperelliptic and $P$ is a Weierstrass point.
The other direction of (ii) is a special case of Proposition~\ref{a2n=0}.

\medskip

Before proving (i), note that if $2 \mid N$, then
Proposition~\ref{pinlevelprop} implies that the $\Q$-gonality $G_\Q$ of
$X$ satisfies $G_\Q \leq 4$.  Remark~\ref{genus17} then implies that
$g \leq 34$.  We seek to sharpen both of these statements.

If $X$ is hyperelliptic, then Theorem~\ref{hyperelliptictheorem}(ii) 
implies $g \leq 10$.
Therefore we may suppose that $X$ is not hyperelliptic.  We may then
assume by~(ii) that $4 \nmid N$, so that it makes sense to consider
the Atkin-Lehner involution $W_2$ on $X$.

Let $X' = X/W_2$, let $P'$ be the image in $X'$ of $P$ under the
natural map, and let $f'_1,\ldots,f'_{g'}$ be a basis of newforms
spanning $S_2(X')$.  Then since $f_i \vert W_2 = -a_2(f_i) f_i$ by
(\ref{W_p}), we have $a_2(f'_i) = -1$ for all $i$.  Therefore every regular
differential on $X'$ vanishing at $P'$ vanishes to order at least 2 at $P'$.
We claim that $X'$ has genus 1.  Indeed, if the genus of
$X'$ were at least 2, then Corollary~\ref{RRcor} would then imply that
$X'$ is hyperelliptic and that $P'$ is a Weierstrass point, but then
we find by applying (ii) to $X'$ that $4 \mid N$, contrary to what
we just assumed.  Also, $X'$ cannot have genus 0 or else $X$ would be
hyperelliptic, contrary to assumption.   Thus $X'$ has genus 1 as
claimed, and $X$ is a degree 2 cover (over $\Q$) of the elliptic curve
$X'$.  

If $3\nmid N$, then we have 
$\# X(\F_9) \leq 2 \# X'(\F_9) \le 2(9+1+2\sqrt{9}) = 32$.
But also $(3-1)(g-1) < \# X(\F_9)$ by Lemma~\ref{ogglemma}.  
Thus $g \leq 16$.

Finally, if $3\mid N$, then the conclusion $g \leq 16$ follows from
the stronger conclusion of~(iii) proved below. 

\medskip

To prove~(v), we suppose that $18 \mid N$.  
It follows by Theorem~\ref{hyperelliptictheorem}(ii)
that if $X$ is hyperelliptic then $g < 3$.  Thus by (i), we may
assume that $X$ is a double cover of an elliptic curve $X'$.

Since $9 \mid N$, we have $a_3(f_i) = 0$ for all $i$, so
Lemma~\ref{RRlemma} implies that $X$ is trigonal (i.e., $X$ admits a
degree 3 map to $\PP^1$).  It therefore follows from
Lemma~\ref{CSlemma} that $g \leq 4$.

\medskip

To prove (iii), we may assume by (v) that $3 \mid N$ but $9 \nmid
N$.  As before, we may also assume that $2 \mid N$ but $4 \nmid N$,
that $X$ is not hyperelliptic, and that $X$ is a double cover of the elliptic curve
$X' = X / W_2$.  

Let $w_j(f_i)$ denote the eigenvalue of $W_j$ on the newform $f_i$ for
$j = 2,3,6$ and $i = 1,\ldots,g$. 
We have $w_2(f_i) =-a_2(f_i), w_3(f_i)=-a_3(f_i),$ and
$w_6(f_i)=a_6(f_i)$ for all $i$.  
Let $g_j$ be the genus of $X/W_j$ ($j = 2,3,6$), 
so $g_j = \#\{\,i \mid w_j(f_i)=+1 \,\}$.
Since $X/W_2$ has genus 1, all the $w_2(f_i)$'s are equal to $-1$
except for one, which is equal to $+1$.
Suppose that $f$ is the newform on which $W_2$ acts with eigenvalue $+1$.
The identity $w_2 w_3 = w_6$ and the fact that the $w_2$'s
are $+1, -1, -1, -1, \ldots, -1$
implies that $g_3+g_6$ is equal to either $g-1$ or $g+1$,
depending on whether $f \vert W_3$ equals $-f$ or $+f$.  We consider
these two cases separately.

\medskip
\noindent\underline{Case 1}: $f \vert W_3 = -f$.

If $g_3 \geq 2$, then $X/W_3$ is a curve of genus at least 2 with all $w_2$'s equal to $-1$,
so as before we would have $4 \mid N$, contrary to assumption.  Thus
$g_3 \leq 1$.  Similarly, $g_6 \leq 1$.
But $g_3+g_6=g-1$, so $g \leq 3$.

\medskip
\noindent\underline{Case 2}: $f \vert W_3 = f$.

In this case, we have $g_3 \geq 1$.  Also $f \vert W_6 = f$, so $g_6 \geq 1$.
Assume first that $g_3=1$.  Then $w_2(f_j) = w_3(f_j)$ for all $j$,
so a simple argument using Riemann-Roch (as in the proof of
Lemma~\ref{RRlemma}) 
shows that there is a rational function on $X$ in $L(3P)-L(2P)$.
Since (as in the proof of Proposition~\ref{pinlevelprop}) 
all $a_4$'s are equal to $+1$, there is also a function in
$L(4P)-L(3P)$.  Taking products, we also find functions in
$L(m P)-L((m-1)P)$ for all $m\geq 6$, so that the Weierstrass gap
sequence at $P$ is contained in $\{ 1,2,5 \}$.  
By Riemann-Roch, the gap sequence at any point contains 
exactly $g$ integers, 
so $g\leq 3$.
 
Assume now that $g_6=1$.  Then $w_2(f_j)=w_6(f_j)$ for all $j$, so
$w_3(f_j)=1$ for all $j$.  Corollary~\ref{RRcor} gives a function in $L(3P)-L(2P)$,
and as before we deduce that $g\leq 3$.

In the general case, we may apply the previous reasoning to $X/W_3$,
since $(X/W_3)/W_6$ is of genus 1 (its only newform is $f$).
The conclusion is that $g_3 \leq 3$.  Similarly, since $(X/W_6)/W_3$
has genus 1, we conclude that $g_6 \leq 3$.  Since $g_3+g_6 = g+1$ in Case
II, we have $g \leq 5$ as desired.

\medskip

To prove (vi), suppose that $36 \mid N$.  Then $a_2(f_i) = a_3(f_i) =
0$ for all $i$.  By the same argument as
in the proof of (ii), it follows that every differential on $X$ vanishing at $P$ vanishes to
order at least 3 at $P$, contradicting Corollary~\ref{RRcor}.

\medskip

For (iv), assume $12 \mid N$.  Then $X$ is hyperelliptic by part (ii).
By Theorem~\ref{hyperelliptictheorem}(ii), it follows that $g = 2$ as desired.
\end{proof}


We can deduce stronger results if we assume furthermore that the
jacobian of $X$ is $\Q$-simple.

\begin{prop}
Let $X$ be a new modular curve of genus
$g\geq 2$, level $N$, and trivial character.  Assume furthermore that
the jacobian $J$ of $X$ is $\Q$-simple.  Then:
\begin{itemize}
\item[(i)] If $2 \mid N$, then $4 \mid N$.
\item[(ii)] $6 \nmid N$.
\end{itemize}
\end{prop}

\begin{proof}
Let $f_1,\ldots,f_g$ as before be a basis of newforms for $\Gamma_0(N)$ 
spanning $S_2(X)$.  
Since $J$ is $\Q$-simple, the $f_i$'s are all Galois conjugates of one
another.
To prove (i), suppose $2 \mid N$ but $4 \nmid N$.  
Then $a_2(f_i) = \pm 1$ for all $i$.  Since in particular $a_2(f_i) \in \Q$, all the
$a_2(f_i)$'s must be equal to one another.  
It follows by Corollary~\ref{RRcor} that $X$ is
hyperelliptic and $P$ is a Weierstrass point.  But then $4\mid N$ by
Proposition~\ref{2and3levelprop}(ii), contradicting our assumption.

\medskip

To prove (ii), suppose that $6 \mid N$.  By (i) we have in fact
that $12 \mid N$.  Also, we know by Proposition~\ref{2and3levelprop}(vi)
that $36 \nmid N$.  Therefore we have
$a_2(f_i) = 0$ and $a_3(f_i) = \pm 1$ for all $i$.  As above, it
follows that all the $a_3(f_i)$'s are in fact equal, contradicting Corollary~\ref{RRcor}.

\end{proof}

\section{Examples and Pathologies}
\label{examplesandpathologiessection}

\subsection{Examples of nonnew hyperelliptic curves}
\label{examplesection}

We can construct nonnew modular hyperelliptic curves from 
certain pairs of new modular hyperelliptic curves.
For some levels, there exist at least two new modular hyperelliptic curves 
of level $N$ such that it is possible to take 
the same modular function $x$ in both their equations. 
In other words, $\Q(x,\sqrt{P_1(x)},\sqrt{P_2(x)})$ 
is a subfield of the function field of $X_1(N)$, 
where the two new modular curves are
$$
C_{N,1}: y_1^2= P_1(x)\,, \qquad C_{N,2}: y_2^2=P_2(x).
$$
Then $\Q(x,\sqrt{P_1(x)P_2(x)})$ also is 
a subfield of the function field of $X_1(N)$, 
so $y^2 = P_1(x) P_2(x)$ is a modular hyperelliptic curve $C_N'$ 
dominated by $X_1(N)$. 
(If desired, one may discard square factors of $P_1(x)P_2(x)$ 
to make the right hand side squarefree.) 
In each example we give, $C_N'$ is not new of level $N$, 
although sometimes it is new of a smaller level.

We found $33$ nonnew modular curves $C_N'$ of this type.
Five of these $33$ are not {\em primitive} 
(cf.\ definition~2.2 in~\cite{gogo00}): 
this means that the minimum $N'$ 
such that $\Jac C_N'$ is a quotient of $J_1(N')$ 
is different from the minimum $M$ such that $C_N'$ is dominated by $X_1(M)$.
The following table shows these five curves.

\setlongtables
\begin{longtable}{llll}
\toprule
 $C_N'$ & Decomposition of $\Jac C$ & $M$ & $g'$\\
\midrule
\endfirsthead
\toprule
 $C_N'$ & Decomposition of $\Jac C$ & $M$ & $g'$\\
\midrule
\endhead
\midrule
\endfoot
\bottomrule
\endlastfoot

$ C_{184}^{3,3}$ & $A_{23A}=J_0(23)$ & $46$ & $2$\\[.1cm]

$ C_{248}^{2,3}$ & $A_{31A}=J_0(31)$ & $62$ & $2$\\[.1cm]

$ C_{376}^{2,2}$ & $A_{47A}=J_0(47)$ & $94$ & $4$\\[.1cm]

$ C_{544}^{2,3}$ & $E_{34A}\times A_{68A}$ & $136$ & $3$\\[.1cm]

$ C_{704}^{3,4}$ & $E_{44A}\times E_{88A}\times A_{88B}$ & $176$ & $4$\\[.1cm]
\end{longtable}

In the first column appears the label for each curves. 
For each curve the subscript denotes the level $N$ 
and the superscript denotes the genus of the curves $C_{N,1}$ and $C_{N,2}$.
In the second column appears the decomposition over $\Q$ 
of the jacobians of the curves in the first column 
(throughout this section we use the labelling of modular forms and abelian varieties
described in the appendix). 
In the third column appears the minimum level $M$
such that $X_1(M)$ dominates $C_N'$,  
and in the last column appears the genus of $C_N'$.

\medskip

In the rest of this subsection, we discuss the surprising case $N=376$.  
We have that $J_0(376)^\new$ splits over $\Q$ 
as the product of the $\Q$-simple modular abelian varieties 
$A_{376A}$, $A_{376B}$, $A_{376C}$ and $A_{376D}$. 
Furthermore, each of these abelian varieties is $\Q$-isogenous 
to the jacobian of a new modular hyperelliptic curve. 
Namely, we have
$$
A_{376A}\Qisog\Jac\,C^A_{376}\,,\quad A_{376B}\Qisog\Jac\,C^B_{376}\,,\quad A_{376C}\Qisog\Jac\,C^C_{376}\,,\quad A_{376D}\Qisog\Jac\,C^D_{376}\,,
$$
where the four curves have the following affine models:
$$
\begin{array}{lcl}
C^A_{376}: y_A^2=P_A(x), &  C^C_{376}:y_C^2=P_B (x) Q(x),\\
C^B_{376}: y_B^2=P_B(x), &  C^D_{376}: y_D^2=P_A(x) Q(x),\\
\end{array}
$$
where 
$$
\begin{array}{rcl}
x &= &q^{-2} + q^2 + q^6 + q^8 + q^{10} + \dots \,,\phantom{cc\,\,}\\
P_A(x)&=&x^5 - x^3 + 2 x^2 - 2 x + 1\,, \phantom{ccccccc}\\
P_B(x)&=&x^5 + 4 x^4 + 3 x^3 - 2 x^2 + 2 x + 5\,,\\
Q(x)&=& x^4 - 2 x^3 - 3 x^2 + 4 x - 4\,.\phantom{ccccccc}
\end{array}
$$
Since the modular function $x$ is the same for all four curves, 
$X_0(376)$ also dominates the hyperelliptic curves
$C^{2,2}_{376}: (y_A y_B)^2=P_A(x) P_B(x)$ and 
$C^{2,4}_{376}: (y_B y_D)^2=P_A(x) P_B(x) Q(x)$ 
(and the genus-one curve $(y_B y_C/P_B(x))^2=Q(x)$ 
isomorphic to the curve labelled 94A1 in~\cite{cremona97}). 
In addition, $X_0(376)$ dominates $X_0(47)$, which also is hyperelliptic. 
We will identify the smallest levels of $C^{2,2}_{376}$ and $C^{2,4}_{376}$.

First consider $C^{2,2}_{376}$. 
Let $f$ be the newform labeled by $47A$. It is easy to check that
there is a modular hyperelliptic curve $X$ over $\Q$ of level $94$
attached to the $\C$-vector space spanned by the Galois
conjugates of the eigenform $h(q)= f(q)-2 f(q^2)\in S_2(94)$
and that $X$ is the curve $C^{2,2}_{376}$.
Therefore, $\Jac C^{2,2}_{376}$ and $J_0(47)$ are $\Q$-isogenous.
But one can check that $C^{2,2}_{376}$ and $X_0(47)$ are not
$\Q$-isomorphic.
Hence $94$ is the minimum level
for $C^{2,2}_{376}$ while $47$ is the minimum level for its jacobian.

Now consider $C^{2,4}_{376}$.
We compute that $\Jac C^{2,4}_{376} \Qisog J_0(47) \times J_0(94)^{\new}$.
More precisely, $C^{2,2}_{376}$ is a nonnew modular curve of minimum level $94$:
it is attached to the $\C$-vector space generated by the
Galois conjugates of $f(q)+2 f(q^2)$ and $g$, where $f$ is as above
and $g$ is the newform labeled by $94B$ with $\dim A_{94B}=2$.
The curve $C^{2,4}_{376}$ has genus $6$,
and its jacobian has one factor in the new part and another in the old part.

Finally, we note that the highest genus yet found 
for a modular hyperelliptic curve is $7$:
the curve is given by the equation
$$
C^{2,5}_{1664}:y^2=(x^2-2x+2)(x^2+2x+2)(x^3-x+2)(x^3-x-2)(x^6+2 x^4+ x^2+4)\,,
$$
and has level $N=1664$. 
Its jacobian splits as 
$\Jac C^{2,5}_{1664}\Qisog E_{26A} \times A_{104B} \times A_{416 F}$.

\subsection{``Pathologies''}
\label{pathologysection}

There are many statements about elliptic curve quotients of $J_1(N)$
that fail for modular curves of higher genus or for higher-dimensional
abelian quotients of $J_1(N)$.
We begin with a statement about abelian quotients.

\begin{prop}
\label{abelianpathology}
There exists an abelian variety $A$ over $\Q$ such that
the following statement is {\em false}:
``If $A$ is a quotient of $J_1(N)$ and $N$ is the
smallest integer with this property,
and $M$ is any other integer such that $A$ is a quotient of $J_1(M)$,
then $N$ divides $M$.''
\end{prop}

\begin{proof}
Take $A=J_0(22)$.
One computes that there are no newforms on $\Gamma_0(22)$,
so $A \Qisog J_0(11)^2$.
Now $A$ is not a quotient of $J_1(M)$ for any $M$ prime to $11$
because $A$ has bad reduction at $11$,
but $A$ is also not a quotient of $J_1(11)$, since $\dim J_1(11)=1$.
Hence $N=22$ is the smallest integer for which $A$ is a quotient
of $J_1(N)$.
On the other hand, (\ref{J_1(N)decomposition}) shows that
$A$ is also a quotient of $J_1(11m)$ for any $m \ge 2$.
\end{proof}

Next we list a few false statements about modular curves.

\begin{prop}
\label{pathologies}
For each of the statements below, there is a curve $X$ over $\Q$ of
genus $g \ge 2$ for which the statement is {\em false}.  Let $J=\Jac
X$.
\begin{enumerate}
\item \label{modular=newmodular}
        ``If $X$ is modular, then $X$ is a new modular curve of some level.''
\item \label{nonnew=old}
        ``If $X$ is nonnew modular of level $N$, 
        then $J$ is a quotient of $J_1(N)_\old$.''
\item \label{jacobianlevel=curvelevel}
        ``If $X$ is modular and $N$ is the smallest level such that $J$ is a
        quotient of $J_1(N)$ defined over $\Q$, then $N$ is the smallest level
        of $X$.''
\item \label{primesforjacobian=primesforcurve}
        ``If $X$ is modular,
        then the smallest level of $X$ and $\cond(J)$ are divisible 
        by the same rational primes.''
\item \label{modularjacobian=modularcurve}
        ``If $J$ is a $\Q$-simple quotient of $J_1(N)^\new$ 
        and $X(\Q)$ is nonempty,
        then $X$ is a new modular curve of level $N$.''
\end{enumerate}
In fact,
there exist infinitely many $X$ over $\Q$ of genus $2$
for which~(\ref{modularjacobian=modularcurve}) fails.
\end{prop}

\begin{proof} \nichts

(\ref{modular=newmodular})
Let $X=X_1(p^2)$ where $p \ge 11$ is prime. 
It follows from an explicit formula for the genus of $X_1(N)$ (see Example
9.1.6 of \cite{diamond-im1995}) that $0 < 2g(X_1(p)) < g(X_1(p^2))$.
Therefore, $\dim J_1(p^2)^\old=2 g(X_1(p))>0$ and $\dim
 J_1(p^2)^\new=g(X_1(p^2))- \dim J_1(p^2)^\old>0$. Now, by Proposition
 \ref{A_fisogenies} 
we obtain that if $J_1(p^2)^\new$ is a quotient of $J_1(N)^\new$ over $\Q$ then $N$ must be $p^2$. Hence, $X$ is not new for any level 
because $J_1(p^2)^\old$ is not a quotient of $J_1(p^2)^\new$.

(\ref{nonnew=old})
Let $X=X_1(p^2)$ for prime $p \ge 11$ again.
Then $X$ is modular of level $p^2$ and not new of any level,
but $\Jac X$ is not a quotient of $J_1(p^2)^\old$.

(\ref{jacobianlevel=curvelevel})
Let $X$ be the genus~$4$ curve $C^{2,2}_{376}$ constructed above.
Then $\Jac X$ is a quotient of $J_0(47)$ (they are isogenous),
but the smallest $N$ for which $X$ is modular of level $N$ is $94$.

(\ref{primesforjacobian=primesforcurve})
Again let $X=C^{2,2}_{376}$.
Then $\cond(\Jac X)=47^4$
but the smallest $N$ for which $X$ is modular of level $N$ is $94$.
A simpler example is the genus~2 curve $X=X_0(22)$:
then $\cond(\Jac X)=11^2$,
but the smallest $N$ for which $X$ is modular of level $N$ is $22$.

(\ref{modularjacobian=modularcurve}) 
Let $a \in \Z - \{0\}$, 
and let $X_a$ be the smooth projective model of the affine curve
        $$y^2=x^5+a x^3-4 x$$
over $\Q$ of genus two.
Then $X_a(\Q)\not=\emptyset$, since $(0,0) \in X_a(\Q)$. 
We will show that for infinitely many values of $a$, 
the jacobian $J_a:=\Jac X_a$ is a $\Q$-simple quotient of some $J_1(N)$. 
It can be checked that the group $\Aut (X_a)$ 
is isomorphic to $D_{2 \cdot 4}$ and 
is generated by the hyperelliptic involution 
and the automorphisms $u$, $v$ over $\Q(i)$, 
represented as in Lemma~\ref{lemma2}(1) by 
$M_u=\begin{pmatrix} 0 & 1+i \\ (1-i)/2 & 0 \end{pmatrix}$, 
$M_v=\begin{pmatrix} i & 0 \\ 0 & -i \end{pmatrix}$ (and $e=\det M$).
It can be checked that the Galois action on $\Aut(X_a)$ 
corresponds to the case $C_2^C$ of section~3 of~\cite{cq2002}.
By Proposition 5.3 of \cite{cq2002}, 
we obtain $(\End J_a)\otimes \Q \isom \Q(\sqrt 2)$.
In particular, $J_a$ is $\Q$-simple. 
On the other hand, the quotient $E_a = X_a/u$
is the elliptic curve
        $$ y^2=x^3+ 27(3a-20 i) x-108(28+9a i)(1-i) \,,$$
whose $j$-invariant $j_a$ equals $64(3a-20i)^3(a-4i)/(a^2+16)^2$. 
Thus $J_a$ is isogenous to the Weil restriction of $E_a$ over $\Q$. 
Solving $j_a=\overline{j}_a$ for $a$ shows that $j_a$ is never real,
so $E_a$ does not have CM.
Thus $(\End J_a)\otimes \Q$ is a real quadratic field, 
and the sign of the $2$-cocycle attached to $E_a$ by Ribet in~\cite{ribet94} 
is trivial (see Theorem~5.4 of~\cite{quer}).
In particular, its local component at $3$ is trivial.
For every nonzero $a\in \Z$, 
$E_a$ has good ordinary reduction at $3$, 
so Theorem~5.1 of~\cite{ellenbergskinner2000} 
implies that $E_a$ is modular (over $\Q(i)$),
so $J_a$ is modular.
Since $J_a$ is $\Q$-simple, it must be a quotient of $J_1(N)^\new$
for some $N$ depending on $a$.
Since $j_a$ is nonconstant, the family $X_a$ is not isotrivial.
(In fact, it can be proved that $X_a$ and $X_b$ are isomorphic over $\Qbar$
if and only if $a=\pm b$.) 
Hence by Theorem~\ref{genustheorem},
at most finitely many of the $X_a$ can be 
new modular curves.
(In fact, $X_a$ is a new modular curve exactly 
for $a=\pm 1,\pm 2, \pm 3, \pm 4, \pm 5, \pm 7, \pm 8,\pm 10,\pm 22$.)
 \end{proof}

\begin{question}
Is there a curve analogue of Proposition~\ref{abelianpathology}?
More precisely, does there exist a curve $X$ dominated by
some $X_1(N)$, but such that the set of such $N$ for $X$
are not all multiples of the smallest $N$?
\end{question}

\section{Curves dominated by Fermat curves}
\label{fermatsection}

Here we prove an analogue of Conjecture~\ref{fixedgenusconj} 
for Fermat curves.
Let $k$ be a field of characteristic zero and let $N \ge 1$.
The $N$-th Fermat curve $X_{N,k}$ 
is the smooth plane curve over $k$
given by the homogeneous equation $x^N+y^N=z^N$ in $\PP^2_k$.
Our goal in this section is to prove the following.

\begin{theorem}
\label{fermattheorem}
For each field $k$ of characteristic zero, and integer $g \ge 2$,
the set of genus-$g$ curves over $k$ dominated by $X_{N,k}$
for some $N \ge 1$
is finite.
If $k$ is a number field or $\Qbar$, then the set is also computable.
\end{theorem}

Let $J_N$ denote the jacobian of $X_{N,k}$.
For each $M | N$, there is a morphism 
$\psi_{N,M}:X_{N,k} \rightarrow X_{M,k}$
mapping $(x:y:z)$ to $(x^{N/M}:y^{N/M}:z^{N/M})$.
Define $J_{N,\old} = \sum_{M | N} \psi_{N,M}^* J_M$,
and $J_N^\new = J_N/J_{N,\old}$.
We may consider $X_{M,k}$ as the quotient of $X_{N,k}$
by a $G_k$-stable subgroup $\Gamma_{N,M} \subseteq \Aut(X_{N,\kbar})$.
If we decompose the space $H^0(X_{N,\kbar},\Omega)$ into characters 
of $\Gamma_{N,1}$, and group those for which a particular divisor $M$
is the smallest $M$ for which $\Gamma_{N,M}$ acts trivially,
then we see that $J_N$ is $k$-isogenous to $\prod_{M|N} J_M^\new$.

\begin{lemma}
\label{aokilemma}
If $N>180$, then each $k$-simple quotient of 
$J_N^\new$ has dimension at least $\phi(N)/8$,
where $\phi(N) := \#(\Z/N\Z)^*$.
\end{lemma}

\begin{proof}
We may assume $k=\C$.
According to Theorem~1.3 of \cite{Aoki}, $J_N^\new \sim \prod A_S$,
where each $A_S$ is an abelian variety of dimension $\phi(N)/2$, and
$S$ ranges over the elements of a certain finite set $\Sigma$.
Moreover, if $N$ is not in an explicit finite set of natural numbers
(whose largest element is $180$), then Theorem~0.2 of \cite{Aoki} shows 
that for each $S \in \Sigma$, $A_S = B_S^{W_S}$, where $B_S$ is a simple
abelian variety and $W_S \leq 4$.  The result follows.
\end{proof}

\begin{lemma}
\label{fermatdominated}
There exists a function $M=M(g)$ such that the following holds:
If $k$ is a field of characteristic zero, 
and $X$ is a curve of genus $g \ge 2$ over $k$ dominated by $X_{N,k}$
for some $N \ge 1$,
then $X$ is dominated also by $X_{M',k}$ for some $M' \le M$.
\end{lemma}

\begin{proof}
Choose $m>180$ such that $\phi(n)/8 > g$ for all $n>m$.
Let $M=m!$.
Suppose $X$ is a curve of genus $g \ge 2$ over $k$ dominated by $X_{N,k}$
for some $N \ge 1$.
Let $M'=\gcd(M,N)$.
By Lemma~\ref{aokilemma},
the composition $\Jac X \to J_N \to J_{m'}^\new$
is trivial for each $m'|N$ greater than $m$,
so the image of $\Jac X \to J_N$
is contained in an abelian subvariety of $J_N$ isogenous to
$\prod_{m'|N, m' \le m} J_{m'}^\new$,
which in turn is contained in the image of $J_{M'} \to J_N$.
Considering differentials and applying Proposition~\ref{fieldfunction}(i),
we find that the morphism $X_{N,k} \rightarrow X$
factors through $X_{M',k}$.
\end{proof}

\begin{proof}[Proof of Theorem~\ref{fermattheorem}]
Use Lemma~\ref{fermatdominated},
and then apply Theorem~\ref{defranchis} (our de~Franchis-Severi Theorem)
to $X_{M',k}$ for each $M' \le M(g)$ to obtain a finite list
of curves dominated by Fermat curves over $k$, 
guaranteed to contain all those of genus $g$.
If we now assume that $k$ is a number field or $\Qbar$,
then this list is also computable.
To eliminate curves of genus not $g$,
and to eliminate possible redundancy,
apply parts (\ref{computinggenus}) and~(\ref{decidingbirational}) 
of Lemma~\ref{birationaldecision}.
\end{proof}

\begin{remark}
$\left.\right.$
\begin{enumerate}
\item[(i)] 
If $k$ is a number field, we can obtain finiteness and computability
also for dominated curves of genus $\le 1$ in Theorem~\ref{fermattheorem}.
This follows from Remark~\ref{genus1defranchis}: 
the genus-$1$ isomorphism problem is not a problem here,
because each dominated curve has a rational point, 
thanks to the point $(1:0:1)$ on each $X_{N,k}$.
\item[(ii)] 
The natural analogue of Theorem~\ref{fermattheorem} for bounded
{\em gonality} is false.
In affine coordinates, $(x,y) \mapsto (xy,y^N)$ defines a 
dominant morphism from $x^N+y^N=1$ to the hyperelliptic curve $u^N=v(1-v)$.
This gives infinitely many Fermat-dominated curves of gonality $2$.
These curves are even new, if $N$ is prime.
Perhaps the failure of the gonality result is not surprising,
since in contrast with Theorem~\ref{Abramovichtheorem},
the gonality of high-degree smooth plane curves
is small compared to the genus.
\end{enumerate}
\end{remark}

\begin{question}
Are there finiteness results for curves over $\Q$ of fixed genus $g \ge 2$ 
that are dominated by curves in some other families,
such as the family of Shimura curves associated to orders 
in a fixed or varying indefinite quaternion algebra over $\Q$?
\end{question}


\section*{Appendix}


\subsection*{Labeling}
We use a deterministic procedure to label newforms,
and in particular to fix an ordering of (Galois conjugacy classes of) 
newforms having a given level and Nebentypus.
An ordering was introduced by J.~Cremona in~\cite{cremona97} 
in the case of trivial Nebentypus and weight~$2$ 
and generalized by W.~A.~Stein in~\cite{stein:tesis} 
to the case of weight greater than $2$. We are going to use this last labeling for 
the case of weight $2$ and arbitrary Nebentypus. 

We will define a function mapping each newform 
$f \in \bigcup_{N=1}^\infty \New_N$
to a label of the form $NX_\eps$ (for example, $13A_{\{2\}}$),
where $N$ is the level of $f$,
where $X$ is a letter or string
in $\{A,B,\dots,Z,AA,BB,\dots\}$
and $\eps$ is (some encoding) of the Nebentypus of $f$.
If $\eps=1$, we omit the subscript $\eps$
and use a label of the form $NX$ (for example, $11A$), 
and if the Fourier coefficients of $f$ are integers 
we will use the labeling in~\cite{cremona97}. 
The labeling function will not be injective:
our definition will be such that
if $f$ has label $NX_\eps$ and $\sigma \in G_\Q$,
then $\ssigma f$ will have label $NX_{\sigma(\eps)}$,
which could be the same as $NX_\eps$, even if $\ssigma f \not=f$.

Two things must be explained: how is $X$ constructed from $f$,
and how is $\eps$ encoded?
First we construct $X$.
Fix $N$ and $\eps$.
To $f = \sum a_n q^n \in \New_N \cap S_2(N,\eps)$
associate the infinite sequence of integers 
${\mathbf t}_f = (\Tr_{E_f/\Q} a_1, \Tr_{E_f/\Q} a_2, \dots)$, 
where $E_f$ is the number field $\Q(a_1,a_2,\dots)$.
Choose $X \in \{A,B,\dots,Z,AA,BB,\dots\}$
according to the position of ${\mathbf t}_f$
in the set $\{\, {\mathbf t}_g : g \in \New_N \cap S_2(N,\eps) \,\}$ 
sorted in increasing dictionary order. 
Notice that ${\mathbf t}_f$ determines the Galois conjugacy class of $f$.

The encoding of Dirichlet characters we now describe 
was suggested by J.~Quer. 
Suppose $\eps:(\Z/N\Z)^*\rightarrow \C^*$ 
is a Dirichlet character.
Let $N=\prod{p_n^{\alpha_n}}$ be the prime-ordered factorization.
Then there exist unique $\eps_{p_n}:(\Z/p_n^{\alpha_n}\,\Z)^*\rightarrow \C^*$ 
such that $\eps=\prod{\eps_{p_n}}$. 
If $p$ is an odd prime, let $g_p$ be the smallest positive integer 
that generates $(\Z/p^{\alpha}\,\Z)^*$,
and if $p=2$ and $\alpha \le 2$,
let $g_p=-1$;
in these cases $\eps_{p}$ is determined by the integer 
$e_p \in [0,\varphi(p^\alpha))$
such that $\eps_p(g_p)=e^{2\pi i e_p/\varphi(p^\alpha)}$.
If $p=2$ and $\alpha > 2$,
then $\eps_2$ is determined by the integers 
$e_2',e_2'' \in [0,\varphi(2^\alpha))$
such that $\eps_2(-1)=e^{2\pi i e_2'/\varphi(2^\alpha)}$
and $\eps_2(5)=e^{2\pi i e_2''/\varphi(2^\alpha)}$,
and we write $e_2 = \{e_2',e_2''\}$.
(Note: here and in the next sentence, we use set notation
although we mean sequences.)
Assuming that $N$ is implicit, we denote $\eps$ by $\{e_p\,:\,p|N\}$.

\medskip

If $f\in S_2(N,\eps)$ is a newform with label $NX_{\eps}$, 
then $A_{NX_{\eps}}$ will denote the corresponding modular 
abelian variety $A_f$,
except that when $\dim A_f=1$,
we instead follow the labeling in~\cite{cremona:ftp} 
and use the letter $E$ instead of $A$ 
to denote the modular elliptic curve $A_f$. 

We give an example illustrating the above notation: 
the $2$-dimensional space $S_2(13)$
splits as $S_2(13,\eps) \directsum S_2(13,\eps^{-1})$,
where $\eps=\{2\}$ is a Dirichlet character modulo $13$ (of order $6$), 
and there is only one Galois conjugacy class of newforms of level $13$.
Thus $J_1(13) \Qisog A_{13A_{\{2\}}}$.

\medskip

Each new modular hyperelliptic curve
in our tables will be denoted $C_{{N L_1,\dots,L_n}_\eps}^{M_1,\dots,M_m}$ 
where $N$ is the level,
the $L_i$ and $M_j$ are letters 
(or in one case, the string $AA$)
indicating the simple factors of the jacobian $J$ of the curve,
and $\eps$ is (the label of) a Dirichlet character $(\Z/N\Z)^* \to \C^*$.
More precisely, the notation indicates that $J$ is isogenous
to the product of the $\Q$-simple modular abelian varieties
with the following labels:
$A_{NM_j}$ for each superscript $M_j$ and $A_{{NL_i}_\eps}$
for each subscript $L_i$.
(It turns out that in all known cases, there is at most one
Galois conjugacy class of nontrivial characters $\eps$ involved.)

\subsection*{Tables}

Table~\ref{table1} shows the 64 
new modular curves of genus two with jacobian not $\Q$-simple. 
(The $149$ new modular curves of genus two with $\Q$-simple jacobian
were already listed in~\cite{gogo00}.)
Tables \ref{table2}, \ref{table3} and~\ref{table4} 
show the $30$ new modular hyperelliptic curves of genus $>2$ 
with $\Q$-simple jacobian, trivial character, and level $\leq 3000$, 
grouped by their genus. 
Recall that by Proposition~1.3 of~\cite{shimura72},
a new modular curve with $\Q$-simple jacobian and nontrivial character
must have even genus; this explains why we may omit
the hypothesis $\#\D=1$ in Tables \ref{table2} and~\ref{table4}.

Tables \ref{table5}, \ref{table6} and~\ref{table7} 
show the $40$ new modular hyperelliptic curves that have genus $>2$ 
and level $\leq 2000$ such that the jacobian is a quotient of $J_0(N)^\new$ 
that is not $\Q$-simple. 
Table~\ref{table8} shows the five curves with genus~$3$ and $\# \D=3$. 
Table~\ref{table9} contains all new modular hyperelliptic curves 
with $2\mid \# \D$, $w\notin \D$ and genus $3\le g \le 4$. 
All of them  have genus $3$ and their  jacobians are neither $\Q$-simple nor a $\Q$-factor of $J_0(N)^\new$. 
Table~\ref{table10} contains the remaining curves with $\# \D> 1$ 
and level $\le 569$. 
These curves have $\# \D=2,3$ or $4$.

\begin{remark}
By inspection, the first curve in Table~\ref{table10},
namely $C_{52A,B_{\{0,4\}}}^A$, has a degree-$2$ unramified cover
that dominates $y^2=x^6+4x^5+6x^4+2x^3+x^2+2x+1$,
which is a equation for $X_1(13)$.
The fourth curve, $C^C_{{208A,D}_{\{\{0,0\}, 4 \}}}$, 
similarly has a degree-$2$ unramified cover
that dominates $X_1(13)$.
The second and third curves in Table~\ref{table10},
namely $C_{{160A,E}_{\{\{1,0\},1\}}}$ and $C_{{160B,F}_{\{\{1,0\},1\}}}$,
have degree-$2$ unramified covers dominating
$y^2=p(x)$ and $y^2=-p(-x)$, respectively,
where $p(x)=(x-1)(x^2-2x+2)(x^2-x-1)$;
the latter two curves are the new modular curves 
$C_{160,C}$ and $C_{160,D}$ of genus~2 from~\cite{gogo00}.
We do not have a good explanation for this phenomenon.
\end{remark}

\setlongtables
\begin{longtable}{X}
\caption{Not $\Q$-simple, $g=2$\label{table1}}\\
\toprule
 $C$&$\,\,\,y^2 \,=\, F(x)$\\
\midrule
\endfirsthead
\toprule
 $C$&$\,\,\,y^2 \,=\, F(x)$\\
\midrule
\endhead
\midrule
\endfoot
\bottomrule
\endlastfoot
$C_{26}^{A,B}$ & $\,\,\,y^2 = x^6 + 4 x^5 - 12 x^4 - 114 x^3 - 308 x^2 - 384 x - 191 $\\[.1cm]
$C_{37}^{A,B}$ & $\,\,\,y^2 = x^6 - 4 x^5 - 40 x^4 + 348 x^3 - 1072 x^2 + 1532 x - 860 $\\[.1cm]
$C_{50}^{A,B}$ & $\,\,\,y^2 = x^6 + 2 x^5 - 5 x^4 - 30 x^3 - 55 x^2 - 48 x - 16 $\\[.1cm]
$C_{54}^{A,B}$ & $\,\,\,y^2 = x^6 - 34 x^3 + 1 $\\[.1cm]
$C_{56}^{A,B}$ & $\,\,\,y^2 = x^5 + 6 x^4 - 45 x^3 - 490 x^2 - 1503 x - 1564 $\\[.1cm]
$C_{58}^{A,B}$ & $\,\,\,y^2 = x^6 - 2 x^5 + 11 x^4 - 22 x^3 + 21 x^2 - 12 x + 4 $\\[.1cm]
$C_{66}^{A,B}$ & $\,\,\,y^2 = x^6 + 2 x^5 - 5 x^4 - 22 x^3 - 31 x^2 - 24 x - 8 $\\[.1cm]
$C_{80}^{A,B}$ & $\,\,\,y^2 = x^5 + 2 x^4 - 26 x^3 - 132 x^2 - 231 x - 142 $\\[.1cm]
$C_{84}^{A,B}$ & $\,\,\,y^2 = x^5 + 4 x^4 - 25 x^3 - 172 x^2 - 339 x - 222 $\\[.1cm]
$C_{90}^{A,B}$ & $\,\,\,y^2 = x^6 - 18 x^3 + 1 $\\[.1cm]
$C_{91}^{A,B}$ & $\,\,\,y^2 = x^6 + 2 x^5 - x^4 - 8 x^3 - x^2 + 2 x + 1 $\\[.1cm]
$C_{96}^{A,B}$ & $\,\,\,y^2 = x^5 - 34 x^3 + x $\\[.1cm]
$C_{112}^{A,C}$ & $\,\,\,y^2 = x^5 - 2 x^4 + 10 x^3 - 16 x^2 + 21 x - 14 $\\[.1cm]
$C_{112}^{A,B}$ & $\,\,\,y^2 = x^5 - 6 x^4 - 45 x^3 + 490 x^2 - 1503 x + 1564 $\\[.1cm]
$C_{128}^{B,D}$ & $\,\,\,y^2 = x^5 - 24 x^3 + 16 x $\\[.1cm]
$C_{128}^{A,C}$ & $\,\,\,y^2 = x^5 + 24 x^3 + 16 x $\\[.1cm]
$C_{138}^{A,C}$ & $\,\,\,y^2 = x^6 + 8 x^4 + 6 x^3 + 8 x^2 + 1 $\\[.1cm]
$C_{142}^{B,D}$ & $\,\,\,y^2 = x^6 - 2 x^5 - 5 x^4 + 18 x^3 - 19 x^2 + 12 x - 4 $\\[.1cm]
$C_{160}^{A,B}$ & $\,\,\,y^2 = x^5 + 12 x^3 + 16 x $\\[.1cm]
$C_{162}^{A,D}$ & $\,\,\,y^2 = x^6 + 14 x^3 + 1 $\\[.1cm]
$C_{162}^{B,C}$ & $\,\,\,y^2 = x^6 - 10 x^3 + 1 $\\[.1cm]
$C_{184}^{C,D}$ & $\,\,\,y^2 = x^5 - 10 x^3 - 15 x^2 - 9 x - 7 $\\[.1cm]
$C_{189}^{A,C}$ & $\,\,\,y^2 = x^6 - 12 x^4 + 36 x^3 - 48 x^2 + 36 x - 12 $\\[.1cm]
$C_{189}^{A,B}$ & $\,\,\,y^2 = x^6 - 12 x^4 + 12 x^3 + 24 x^2 - 36 x + 12 $\\[.1cm]
$C_{192}^{C,D}$ & $\,\,\,y^2 = x^5 - 14 x^3 + x $\\[.1cm]
$C_{192}^{B,D}$ & $\,\,\,y^2 = x^5 + 4 x^4 - 6 x^3 - 58 x^2 - 111 x - 70 $\\[.1cm]
$C_{192}^{A,C}$ & $\,\,\,y^2 = x^5 - 4 x^4 - 6 x^3 + 58 x^2 - 111 x + 70 $\\[.1cm]
$C_{192}^{A,B}$ & $\,\,\,y^2 = x^5 + 34 x^3 + x $\\[.1cm]
$C_{200}^{C,E}$ & $\,\,\,y^2 = x^5 - 10 x^3 - 15 x^2 + 8 $\\[.1cm]
$C_{240}^{C,D}$ & $\,\,\,y^2 = x^5 - 2 x^4 + 6 x^3 - 13 x^2 + 12 x - 4 $\\[.1cm]
$C_{256}^{A,D}$ & $\,\,\,y^2 = x^5 + 16 x $\\[.1cm]
$C_{264}^{A,B}$ & $\,\,\,y^2 = x^5 + 2 x^4 - 6 x^3 - 23 x^2 - 24 x - 8 $\\[.1cm]
$C_{312}^{B,C}$ & $\,\,\,y^2 = x^5 - 2 x^4 - x^3 + 8 x^2 - 9 x + 3 $\\[.1cm]
$C_{320}^{A,C}$ & $\,\,\,y^2 = x^5 - 2 x^4 - 2 x^3 - 2 x^2 + x $\\[.1cm]
$C_{320}^{D,E}$ & $\,\,\,y^2 = x^5 - 12 x^3 + 16 x $\\[.1cm]
$C_{320}^{B,F}$ & $\,\,\,y^2 = x^5 + 2 x^4 - 2 x^3 + 2 x^2 + x $\\[.1cm]
$C_{336}^{A,F}$ & $\,\,\,y^2 = x^5 - 4 x^4 - 25 x^3 + 172 x^2 - 339 x + 222 $\\[.1cm]
$C_{368}^{A,G}$ & $\,\,\,y^2 = x^5 - 10 x^3 + 15 x^2 - 9 x + 7 $\\[.1cm]
$C_{384}^{A,D}$ & $\,\,\,y^2 = x^5 + 10 x^3 + x $\\[.1cm]
$C_{384}^{B,C}$ & $\,\,\,y^2 = x^5 - 10 x^3 + x $\\[.1cm]
$C_{400}^{B,E}$ & $\,\,\,y^2 = x^5 - 25 x^2 + 20 x - 4 $\\[.1cm]
$C_{400}^{A,H}$ & $\,\,\,y^2 = x^5 - 10 x^3 + 15 x^2 - 8 $\\[.1cm]
$C_{405}^{A,F}$ & $\,\,\,y^2 = x^6 - 12 x^4 + 28 x^3 - 24 x^2 + 12 x - 4 $\\[.1cm]
$C_{405}^{B,F}$ & $\,\,\,y^2 = x^6 - 12 x^4 + 20 x^3 - 12 x + 4 $\\[.1cm]
$C_{448}^{A,D}$ & $\,\,\,y^2 = x^5 - 2 x^4 + 10 x^3 - 2 x^2 + x $\\[.1cm]
$C_{448}^{B,G}$ & $\,\,\,y^2 = x^5 + 2 x^4 + 10 x^3 + 2 x^2 + x $\\[.1cm]
$C_{480}^{B,C}$ & $\,\,\,y^2 = x^5 + 2 x^4 - 4 x^3 - 17 x^2 - 18 x - 6 $\\[.1cm]
$C_{480}^{B,G}$ & $\,\,\,y^2 = x^5 - 7 x^3 + x $\\[.1cm]
$C_{480}^{A,G}$ & $\,\,\,y^2 = x^5 - 2 x^4 - 4 x^3 + 17 x^2 - 18 x + 6 $\\[.1cm]
$C_{528}^{A,D}$ & $\,\,\,y^2 = x^5 - 2 x^4 - 6 x^3 + 23 x^2 - 24 x + 8 $\\[.1cm]
$C_{544}^{B,C}$ & $\,\,\,y^2 = x^5 - 9 x^3 + 16 x $\\[.1cm]
$C_{624}^{C,D}$ & $\,\,\,y^2 = x^5 + 2 x^4 - x^3 - 8 x^2 - 9 x - 3 $\\[.1cm]
$C_{672}^{A,G}$ & $\,\,\,y^2 = x^5 + 5 x^3 + x $\\[.1cm]
$C_{760}^{A,E}$ & $\,\,\,y^2 = x^5 + 3 x^3 + 14 x^2 + 15 x + 5 $\\[.1cm]
$C_{768}^{D,F}$ & $\,\,\,y^2 = x^5 - 4 x^3 + x $\\[.1cm]
$C_{768}^{B,H}$ & $\,\,\,y^2 = x^5 + 4 x^3 + x $\\[.1cm]
$C_{960}^{A,F}$ & $\,\,\,y^2 = x^5 + 7 x^3 + x $\\[.1cm]
$C_{1088}^{L,N}$ & $\,\,\,y^2 = x^5 + 9 x^3 + 16 x $\\[.1cm]
$C_{1344}^{C,F}$ & $\,\,\,y^2 = x^5 - 5 x^3 + x $\\[.1cm]
$C_{1520}^{B,D}$ & $\,\,\,y^2 = x^5 + 3 x^3 - 14 x^2 + 15 x - 5 $\\[.1cm]
$C_{1664}^{F,G}$ & $\,\,\,y^2 = x^5 - 2 x^4 + x^3 + 2 x - 4 $\\[.1cm]
$C_{1664}^{O,S}$ & $\,\,\,y^2 = x^5 + 2 x^4 + x^3 + 2 x + 4 $\\[.1cm]
$C_{2208}^{A,E}$ & $\,\,\,y^2 = x^5 + 2 x^4 + 8 x^3 + 19 x^2 + 18 x + 6$\\[.1cm]
$C_{2208}^{G,I}$ & $\,\,\,y^2 = x^5 - 2 x^4 + 8 x^3 - 19 x^2 + 18 x - 6$\\[.1cm]
\end{longtable}

\setlongtables
\begin{longtable}{X}
\caption{$\Q$-simple, $g=3$, $N\leq 3000$\label{table2}}\\
\toprule
 $C$  & $\,\,\,y^2 \,=\, F(x)$\\
\midrule
\endfirsthead
\toprule
 $C$  & $\,\,\,y^2 \,=\, F(x)$\\
\midrule
\endhead
\midrule
\endfoot
\bottomrule
\endlastfoot
$C_{41}^A$ & $\,\,\,y^2 = x^8 + 4 x^7 - 8 x^6 - 66 x^5 - 120 x^4 - 56 x^3 + 53 x^2 + 36 x - 16$\\[.1cm]

$C_{95}^A$ & $\,\,\,y^2 = (x^4 + x^3 - 6 x^2 - 10 x - 5)(x^4 + x^3 - 2 x^2 + 2 x - 1)$\\[.1cm]

$C_{152}^C$ & $\,\,\,y^2 = x(x^3 - 2 x^2 - 7 x - 8)(x^3 + 4 x^2 + 4 x + 4)$\\[.1cm]

$C_{248}^E$ & $\,\,\,y^2 = (x^3 + x - 1)(x^4 - 2 x^3 - 3 x^2 - 4 x + 4)$\\[.1cm]

$C_{284}^A$ & $\,\,\,y^2 = x^7 + 4 x^6 + 5 x^5 + x^4 - 3 x^3 - 2 x^2 + 1$\\[.1cm]

$C_{284}^B$ & $\,\,\,y^2 = x^7 - 7 x^5 - 11 x^4 + 5 x^3 + 18 x^2 + 4 x - 11$\\[.1cm]

$C_{304}^G$ & $\,\,\,y^2 = x(x^3 - 4 x^2 + 4 x - 4)(x^3 + 2 x^2 - 7 x + 8)$\\[.1cm]

$C_{496}^J$ & $\,\,\,y^2 = (x^3 + x + 1)(x^4 + 2 x^3 - 3 x^2 + 4 x + 4)$\\[.1cm]

$C_{544}^I$ & $\,\,\,y^2 = (x+1)(x^2 + x - 4)(x^4 - x^2 - 4)$\\[.1cm]

$C_{544}^J$ & $\,\,\,y^2 = (x-1)(x^2 - x - 4)(x^4 - x^2 - 4)$\\[.1cm]

$C_{896}^I$ & $\,\,\,y^2 = (x - 2)(x^2 + 2 x - 1)(x^4 - 2 x^2 - 7)$\\[.1cm]

$C_{896}^K$ & $\,\,\,y^2 = (x + 2)(x^2 - 2 x - 1)(x^4 - 2 x^2 - 7)$\\[.1cm]

$C_{1136}^G$ & $\,\,\,y^2 = x^7 - 7 x^5 + 11 x^4 + 5 x^3 - 18 x^2 + 4 x + 11$\\[.1cm]

$C_{1136}^J$ & $\,\,\,y^2 = x^7 - 4 x^6 + 5 x^5 - x^4 - 3 x^3 + 2 x^2 - 1$\\[.1cm]
\end{longtable}

\setlongtables
\begin{longtable}{X}
\caption{$\Q$-simple, $g=4$, $\#\D=1$, $N\leq 3000$\label{table3}}\\
\toprule
 $C$  & $\,\,\,y^2 \,=\, F(x)$\\
\midrule
\endfirsthead
\toprule
 $C$  & $\,\,\,y^2 \,=\, F(x)$\\
\midrule
\endhead
\midrule
\endfoot
\bottomrule
\endlastfoot
$C_{47}^A$ & $\,\,\,y^2 = (x^5 - 5 x^3 - 20 x^2 - 24 x - 19)(x^5 + 4 x^4 + 7 x^3 + 8 x^2 + 4 x + 1)$\\[.1cm]

$C_{119}^A$ & $\,\,\,y^2 = (x^5 - 2 x^4 + 3 x^3 - 6 x^2 - 7)(x^5 + 2 x^4 + 3 x^3 + 6 x^2 + 4 x + 1)$\\[.1cm]

$C_{164}^A$ & $\,\,\,y^2 = x(x^8 + 4 x^7 - 8 x^6 - 66 x^5 - 120 x^4 - 56 x^3 + 53 x^2 + 36 x - 16)$\\[.1cm]

$C_{376}^C$ & $\,\,\,y^2 = (x^4 - 2 x^3 - 3 x^2 + 4 x - 4)(x^5 + 4 x^4 + 3 x^3 - 2 x^2 + 2 x + 5)$\\[.1cm]

$C_{376}^D$ & $\,\,\,y^2 = (x^4 - 2 x^3 - 3 x^2 + 4 x - 4)(x^5 - x^3 + 2 x^2 - 2 x + 1)$\\[.1cm]

$C_{416}^F$ & $\,\,\,y^2 = x(x^2 + 4)(x^3 - 2 x^2 + x - 4)(x^3 + 2 x^2 + x + 4)$\\[.1cm]

$C_{512}^G$ & $\,\,\,y^2 = x(x^4 - 4 x^2 - 4)(x^4 + 4 x^2 - 4)$\\[.1cm]

$C_{656}^I$ & $\,\,\,y^2 = x(x^8 - 4 x^7 - 8 x^6 + 66 x^5 - 120 x^4 + 56 x^3 + 53 x^2 - 36 x - 16)$\\[.1cm]

$C_{752}^G$ & $\,\,\,y^2 = (x^4 + 2 x^3 - 3 x^2 - 4 x - 4)(x^5 - x^3 - 2 x^2 - 2 x - 1)$\\[.1cm]

$C_{752}^I$ & $\,\,\,y^2 = (x^4 + 2 x^3 - 3 x^2 - 4 x - 4)(x^5 - 4 x^4 + 3 x^3 + 2 x^2 + 2 x - 5)$\\[.1cm]

$C_{832}^P$ & $\,\,\,y^2 = x(x+2)(x-2)(x^6 + 2 x^4 - 15 x^2 + 16)$\\[.1cm]

$C_{1216}^W$ & $\,\,\,y^2 = (x^3 - 2 x + 2)(x^6 + 2 x^4 - 7 x^2 + 8)$\\[.1cm]

$C_{1216}^X$ & $\,\,\,y^2 = (x^3 - 2 x - 2)(x^6 + 2 x^4 - 7 x^2 + 8)$\\[.1cm]
\end{longtable}

\setlongtables
\begin{longtable}{X}
\caption{$\Q$-simple, $g=5$, $N\leq 3000$\label{table4}}\\\toprule
 $C$ & $\,\,\,y^2 \,=\, F(x)$\\
\midrule
\endfirsthead
\toprule
 $C$ & $\,\,\,y^2 \,=\, F(x)$\\
\midrule
\endhead
\midrule
\endfoot
\bottomrule
\endlastfoot
$C_{59}^A$ & $\,\,y^2=(x^9+ 2\, x^8- 4\, x^7- 21\, x^6 - 44\, x^5 -60 \, x^4- 61\, x^3-46\, x^2-24\, x-11)$\\
\multicolumn{2}{r}{$(x^3+2\, x^2+1)$}\\[.1cm]
$C_{1664}^Y$ & $\,\,y^2 = (x^2 + 2 x + 2)(x^3 - x + 2)(x^6 + 2 x^4 + x^2 + 4)$\\[.1cm] 
$C_{1664}^{AA}$ & $\,\,y^2 = (x^2 - 2 x + 2)(x^3 - x - 2)(x^6 + 2 x^4 + x^2 + 4)$\\[.1cm]
\end{longtable}

\setlongtables
\begin{longtable}{X}
\caption{Not $\Q$-simple, $g=3$, $\#\D=1$, $N\leq 2000$\label{table5}}\\
\toprule
 $C$ & $\,\,\,y^2 \,=\, F(x)$\\
\midrule
\endfirsthead
\toprule
 $C$ & $\,\,\,y^2 \,=\, F(x)$\\
\midrule
\endhead
\midrule
\endfoot
\bottomrule
\endlastfoot
$C_{35}^{A,B}$ & $\,\,\,y^2 \,=\,(x^2 + 3 x + 1)(x^6 + x^5 - 10 x^4 - 39 x^3 - 62 x^2 - 51 x - 19)$\\[.1cm]
$C_{39}^{A,B}$ & $\,\,\,y^2 \,=\,(x^4 - 3 x^3 - 4 x^2 - 2 x - 1)(x^4 + 5 x^3 + 8 x^2 + 6 x + 3) $\\[.1cm]
$C_{88}^{A,B}$ & $\,\,\,y^2 \,=\, (x - 2)(x^3 - 2 x^2 + 4 x - 4)(x^3 + 2 x^2 - 4 x + 8)$\\[.1cm]
$C_{104}^{A,B}$ & $\,\,\,y^2 \,=\,(x + 2)(x^6 + 4 x^5 - 12 x^4 - 114 x^3 - 308 x^2 - 384 x - 191)$\\[.1cm]
$C_{116}^{A,B,C}$ & $\,\,\,y^2 \,=\,(x + 2)(x^6 + 2 x^5 - 17 x^4 - 66 x^3 - 83 x^2 - 32 x - 4)$\\[.1cm]
$C_{128}^{A,B,D}$ & $\,\,\,y^2 \,=\,(x - 2)(x^2 - 2 x + 2)(x^4 - 12 x^2 + 32 x - 28)$\\[.1cm]
$C_{128}^{B,C,D}$ & $\,\,\,y^2 \,=\,(x + 2)(x^2 + 2 x + 2)(x^4 - 12 x^2 - 32 x - 28)$\\[.1cm]
$C_{160}^{A,C}$ & $\,\,\,y^2 \,=\,(x - 2)(x^2 + 2 x - 7)(x^4 - 4 x^3 + 10 x^2 - 20 x + 17) $\\[.1cm]
$C_{160}^{B,C}$ & $\,\,\,y^2 \,=\, (x + 2)(x^2 - 2 x - 7)(x^4 + 4 x^3 + 10 x^2 + 20 x + 17)$\\[.1cm]
$C_{176}^{A,D}$ & $\,\,\,y^2 \,=\, (x + 2)(x^3 - 2 x^2 - 4 x - 8)(x^3 + 2 x^2 + 4 x + 4)$\\[.1cm]
$C_{184}^{B,E}$ & $\,\,\,y^2 \,=\, (x - 1)(x^3 - 2 x^2 + 3 x - 1)(x^3 + x^2 - x + 7)$\\[.1cm]
$C_{184}^{A,C,D}$ & $\,\,\,y^2 \,=\, (x - 1)(x^6 - x^5 + 4 x^4 - x^3 + 2 x^2 + 2 x + 1)$\\[.1cm]
$C_{196}^{B,C}$ & $\,\,\,y^2 \,=\, (x^3 + 2 x^2 - x - 1)(x^4 - 2 x^3 - 9 x^2 + 10 x - 3)$\\[.1cm]
$C_{208}^{B,E}$ & $\,\,\,y^2 \,=\,(x - 2)(x^6 - 4 x^5 - 12 x^4 + 114 x^3 - 308 x^2 + 384 x - 191) $\\[.1cm]
$C_{224}^{A,D}$ & $\,\,\,y^2 \,=\, x(x - 1)(x + 1)(x^4 - 6 x^2 + 16 x - 7)$\\[.1cm]
$C_{224}^{B,C}$ & $\,\,\,y^2 \,=\, x(x - 1)(x + 1)(x^4 - 6 x^2 - 16 x - 7)$\\[.1cm]
$C_{248}^{B,D}$ & $\,\,\,y^2 \,=\, (x^3 + 4 x^2 + 5 x + 3)(x^4 - 2 x^3 - 3 x^2 - 4 x + 4)$\\[.1cm]
$C_{256}^{B,E}$ & $\,\,\,y^2 \,=\, x(x^2 + 2)(x^4 + 12 x^2 + 4)$\\[.1cm]
$C_{256}^{C,E}$ & $\,\,\,y^2 \,=\, x(x^2 - 4 x + 2)(x^2 - 2)(x^2 + 4 x + 2)$\\[.1cm]
$C_{280}^{A,D}$ & $\,\,\,y^2 \,=\, (x - 1)(x^6 - x^5 + 7 x^3 - 16 x^2 + 15 x - 5)$\\[.1cm]
$C_{368}^{C,I}$ & $\,\,\,y^2 \,=\, (x + 1)(x^3 - x^2 - x - 7)(x^3 + 2 x^2 + 3 x + 1)$\\[.1cm]
$C_{368}^{A,D,G}$ & $\,\,\,y^2 \,=\, (x + 1)(x^6 + x^5 + 4 x^4 + x^3 + 2 x^2 - 2 x + 1)$\\[.1cm]
$C_{416}^{A,E}$ & $\,\,\,y^2 \,=\, x(x^6 - 2 x^5 - 2 x^4 + 2 x^2 - 2 x - 1)$\\[.1cm]
$C_{416}^{B,C}$ & $\,\,\,y^2 \,=\, x(x^6 + 2 x^5 - 2 x^4 + 2 x^2 + 2 x - 1)$\\[.1cm]
$C_{464}^{D,E,F}$ & $\,\,\,y^2 \,=\, (x - 2)(x^6 - 2 x^5 - 17 x^4 + 66 x^3 - 83 x^2 + 32 x - 4)$\\[.1cm]
$C_{496}^{C,G}$ & $\,\,\,y^2 \,=\, (x^3 - 4 x^2 + 5 x - 3)(x^4 + 2 x^3 - 3 x^2 + 4 x + 4)$\\[.1cm]
$C_{560}^{A,G}$ & $\,\,\,y^2 \,=\, (x + 1)(x^6 + x^5 - 7 x^3 - 16 x^2 - 15 x - 5)$\\[.1cm]
$C_{640}^{C,K}$ & $\,\,\,y^2 \,=\, x(x^2 - 2 x - 1)(x^4 + 2 x^3 - 2 x + 1)$\\[.1cm]
$C_{640}^{G,I}$ & $\,\,\,y^2 \,=\, x(x^2 + 2 x - 1)(x^4 - 2 x^3 + 2 x + 1)$\\[.1cm]
$C_{704}^{B,N}$ & $\,\,\,y^2 \,=\, x(x^3 - 4 x + 4)(x^3 + 2 x^2 - 2)$\\[.1cm]
$C_{704}^{C,O}$ & $\,\,\,y^2 \,=\, x(x^3 - 2 x^2 + 2)(x^3 - 4 x - 4)$\\[.1cm]
$C_{784}^{G,M}$ & $\,\,\,y^2 \,=\, (x^3 - 2 x^2 - x + 1)(x^4 + 2 x^3 - 9 x^2 - 10 x - 3)$\\[.1cm]
\end{longtable}

\setlongtables
\begin{longtable}{X}
\caption{Not $\Q$-simple, $g=4$, $\#\D=1$, $N\leq 2000$\label{table6}}\\
\toprule
 $C$ & $\,\,\,y^2 \,=\, F(x)$\\
\midrule
\endfirsthead
\toprule
 $C$ & $\,\,\,y^2 \,=\, F(x)$\\
\midrule
\endhead
\midrule
\endfoot
\bottomrule
\endlastfoot
$C_{224}^{C, D}$ & $\,\,\,y^2 \,=\, x(x^4-2 x^3-5 x^2-2 x+1)(x^4 + 2 x^3 - 5 x^2 + 2 x + 1)$\\[.1cm]
$C_{236}^{ B, C}$ & $\,\,\,y^2 \,=\, x^9 + 2 x^8 - 4 x^7 - 21 x^6 - 44 x^5 - 60 x^4 - 61 x^3 - 46 x^2 - 24 x - 11$\\[.1cm]
$C_{368}^{B,F,H}$ & $\,\,\,y^2 \,=\, (x^3 - x^2 - x - 7)(x^6 + x^5 + 4 x^4 + x^3 + 2 x^2 - 2 x + 1)$\\[.1cm]
$C_{448}^{I,J}$ & $\,\,\,y^2 \,=\, x(x^8 + 14 x^6 + 19 x^4 + 14 x^2 +1)$\\[.1cm]
$C_{704}^{D,E,M}$ & $\,\,\,y^2 \,=\, (x^3 - 4 x - 4)(x^3 - 4 x + 4)(x^3 + 2 x^2 - 2)$\\[.1cm]
$C_{704}^{F,I,P}$ & $\,\,\,y^2 \,=\, (x^3 - 2 x^2 + 2)(x^3 - 4 x - 4)(x^3 - 4 x + 4)$\\[.1cm]
$C_{944}^{J,L}$ & $\,\,\,y^2 \,=\, x^9 - 2 x^8 - 4 x^7 + 21 x^6 - 44 x^5 + 60 x^4 - 61 x^3 + 46 x^2 - 24 x + 11$\\[.1cm]
\end{longtable}

\setlongtables
\begin{longtable}{X}
\caption{Not $\Q$-simple, $g=6$, $\#\D=1$, $N\leq 2000$\label{table7}}\\
\toprule
$C$ & $\,\,\,y^2 \,=\, F(x)$\\
\midrule
\endfirsthead
\toprule
 $C$ & $\,\,\,y^2 \,=\, F(x)$\\
\midrule
\endhead
\midrule
\endfoot
\bottomrule
\endlastfoot
$C_{71}^{ A, B }$ & $\,\,\,y^2 \,=\,(x^7 - 7 x^5 - 11 x^4 + 5 x^3 + 18 x^2 + 4 x - 11)$\\
\multicolumn{2}{r}{$(x^7 + 4 x^6 + 5 x^5 + x^4 - 3 x^3 - 2 x^2 + 1)$}\\[.1cm]
\end{longtable}

\setlongtables
\begin{longtable}{Y}
\caption{$g=3$, $\#\D=3$\label{table8}}\\
\toprule
$C$ & $\,\,\,y^2 \,=\, F(x)$\\
\midrule
\endfirsthead
\toprule
 $C$ & $\,\,\,y^2 \,=\, F(x)$\\
\midrule
\endhead
\midrule
\endfoot
\bottomrule
\endlastfoot
$C^A_{21A_{\{0,2\}}}$ & $\,\,\,y^2 \,=\, (x^2-x+1)(x^6+ x^5-6 x^4-3 x^3+14 x^2-7 x+1) $\\[.1cm]
$C^A_{36A_{\{0,2\}}}$ & $\,\,\,y^2 \,=\, (x+1)(x+2)(x^2+3 x+3)(x^3-9 x-9) $\\[.1cm]
$C^A_{72A_{\{\{0,0\},2\}}}$ & $\,\,\,y^2 \,=\, x(x+1)( x^2+x+1)( x^3-3 x-1)$\\[.1cm]
$C^A_{144A_{\{\{0,0\},2\}}}$ & $\,\,\,y^2 \,=\, ( x-2)(x-1 )( x^2-3 x+3)(x^3-9 x+9)$\\[.1cm]
$C^B_{144B_{\{\{0,0\},2\}}}$ & $\,\,\,y^2 \,=\, (x-1 )x(x^2-x+1 )( x^3-3 x+1) $\\[.1cm]
\end{longtable}

\setlongtables
\begin{longtable}{l@{:}l}
\caption{$g=3$, $2\mid \# \D$, $w\not\in \D$\label{table9}}\\
\toprule
 $C$  & $\,\,\,y^2 \,=\, F(x)$\\
\midrule
\endfirsthead
\toprule
 $C$  & $\,\,\,y^2 \,=\, F(x)$\\
\midrule
\endhead
\midrule
\endfoot
\bottomrule
\endlastfoot
$C_{40A_{\{\{0,0\},2\}}}^A  $  & $\,\,\, y^2 =   x(x + 1)(x + 2)(x^2 - 2 x - 4)(x^2 + 3 x + 1)$\\[.1cm]
$C_{48A_{\{\{1,0\},1\}}}^A  $  & $\,\,\, y^2 = (x + 1)(x^2 - 2 x - 2)(x^2 + x + 1)(x^2 + 2 x + 2)$\\[.1cm]
$C_{64A_{\{\{0,8\}\}}}^A  $  & $\,\,\, y^2 = x(x-1)(x+1)(x^2-2\, x-1)(x^2+2\, x-1)$\\[.1cm]
$C_{80A_{\{\{0,0\},2\}}}^A$  & $\,\,\, y^2 = x(x-1)(x - 2)(x^2 - 3 x + 1)(x^2 + 2 x - 4)$\\[.1cm]
$C_{80A_{\{\{0,0\},2\}}}^B$  & $\,\,\, y^2 = (x+1)(x^2-x-1)(x^4+4\, x^2+8\, x+4)$\\[.1cm]
$C_{128A_{\{\{0,16\}\}}}^B$  & $\,\,\, y^2 =(x - 2)(x^2 - 2 x - 1)(x^4 - 6 x^2 - 16 x + 41) $\\[.1cm]
$C_{128A_{\{\{0,16\}\}}}^D$  & $\,\,\, y^2 = (x + 2)(x^2 + 2 x - 1)(x^4 - 6 x^2 + 16 x + 41)$\\[.1cm]
\end{longtable}

\setlongtables
\begin{longtable}{l@{:}l}
\caption{$g \ge 4$, $\#\D >1$, $N \leq 569$\label{table10}}\\
\toprule
 $C$  & $\,\,\,y^2 \,=\, F(x)$\\
\midrule
\endfirsthead
\toprule
 $C$  & $\,\,\,y^2 \,=\, F(x)$\\
\midrule
\endhead
\midrule
\endfoot
\bottomrule
\endlastfoot
$C^A_{{52A,B}_{\{0, 4 \}}}$  & $\,\,\, y^2 = x(x+1 )(x^3-x^2-4 x-1 )(x^6+4 x^5+6 x^4+2 x^3+x^2+2 x+ 1) $\\[.1cm] 
$C_{{160A,E}_{\{\{1,0\},1\}}}$  & $\,\,\, y^2 =  (x-1 )(x^2-2 x +2)(x^2-x-1 )(x^4-8 x+8 )$\\[.1cm]
$C_{{160B,F}_{\{\{1,0\},1\}}}$  & $\,\,\, y^2 =  (x+1 )(x^2+2 x +2)(x^2+x-1 )(x^4+8 x+8 )$\\[.1cm]
$C^C_{{208A,D}_{\{\{0,0\}, 4 \}}}$  & $\,\,\, y^2 =  x(x-1 )(x^3+x^2-4 x+1 )(x^6-4 x^5+6 x^4-2 x^3+x^2-2 x+ 1)$\\[.1cm]

$C_{512D_{\{\{0,64\}\}}}$  & $\,\,\, y^2 =  x(x^8+24 x^4+16)$\\[.1cm]
\end{longtable}

\section*{Acknowledgements}

We thank Dinakar Ramakrishnan and Ken Ribet for discussing
Proposition~\ref{A_fisogenies} with us.
The proof of Theorem~\ref{defranchis} grew out of discussions
with Matthias Aschenbrenner, Brian Conrad, Tom Graber, Tom Scanlon,
and Jason Starr.
We thank also 
Paul Gunnells, Sorin Popescu, 
William Stein 
and 
Tonghai Yang 
for suggesting references.


\bibliographystyle{plain}
\bibliography{news}

\end{document}